\documentclass{cpamart1}     
\startingpage{1}                      

\authorheadline{S. Tang and R. Xu}
\titleheadline{A quadratic BSDE approach to the 2D sine-Gordon model}

    \usepackage{amsfonts,amssymb,mathrsfs}
        \usepackage[dvipsnames]{xcolor}
        \usepackage[colorlinks=true,linkcolor=RoyalBlue,citecolor=PineGreen,urlcolor=blue]{hyperref}
        \usepackage[ruled,vlined,linesnumbered]{algorithm2e}
        \usepackage{hyperref}
        \usepackage{tablefootnote}
    
          \usepackage[shortlabels]{enumitem}

\usepackage{newtxtext,newtxmath}

\newtheorem{theorem}{Theorem}[section]
\newtheorem{lemma}[theorem]{Lemma}

\newtheorem{proposition}[theorem]{Proposition}
        
\theoremstyle{definition}
\newtheorem{definition}[theorem]{Definition}

\theoremstyle{remark}
\newtheorem{remark}[theorem]{Remark}

\SetLabelAlign{LeftAlignWithIndent}{\makebox[1.5em][l]{#1}}

        \allowdisplaybreaks

\newcounter{desccount}

\newcommand{\descref}[1]{\hyperref[#1]{#1}}

\begin{document}                        


\title{A quadratic-backward-stochastic-differential-equation approach to normalization for the finite volume 2D sine-Gordon model in the finite ultraviolet regime
}

\author{Shanjian Tang}{School of Mathematical Sciences, Fudan University}
\author{Rundong Xu}{School of Mathematical Sciences, Fudan University}





\begin{abstract}
    This paper is devoted to a new construction of the two-dimensional sine-Gordon model on bounded domains
    by a novel normalization technique in the finite ultraviolet regime. Our methodology involves a family of backward stochastic differential equations (BSDEs for short) driven by a cylindrical Wiener process,
    whose generators are purely quadratic functions of the second unknown variable.
    The terminal conditions of the quadratic BSDEs are uniformly bounded and converge in probability to the real part of imaginary multiplicative chaos tested against an arbitrarily given test function,
    which helps us describe our sine-Gordon measure through some delicate estimates concerning bounded mean oscillation martingales.
    As the ultraviolet cutoffs are vanishing, the quadratic BSDEs converge to a quadratic BSDE that completely characterizes the absolute continuity of our sine-Gordon measure with respect to the law of Gaussian free fields.
    Our approach can also be used effectively to establish the connection between our sine-Gordon measure and the scaling limit of correlation functions of the critical planar XOR-Ising model 
    and to prove the weak convergence of the normalized charge distributions of two-dimensional log-gases.
\end{abstract}

\maketitle 


\vspace{0.1cm}

\noindent {\footnotesize {\textbf{MSC2020 subject classification:}} 60H10; 81S20: 81T08}



\tableofcontents




\section{Introduction}
\label{section-1}

Ever since the original work of Parisi and Wu \cite{Parisi-Wu-1981}, the \textit{stochastic quantization} for the $d$-dimensional Euclidean scalar field theory on a set $\Lambda \subset \mathbb{R}^{d}$ has been a rapid developing research area,
where one advocates to construct a parameterized random scalar field $\{\varphi_{t}(\cdot,\cdot) : \Omega \times \mathbb{R}^{d} \mapsto \mathbb{R}, \ t \geqslant 0 \}$ on a probability space $(\Omega, \mathcal{F}, \mathbb{P})$ such that the invariant measure of $(\varphi_{t})_{t \geqslant 0}$ is formally described by the Gibbs measure
\begin{equation}
    \label{introduction-Gibbs-measure-formal}
    \mu_{\mathrm{Gibbs}}(\mathrm{d}\phi) := \Xi^{-1} \exp\left\{ - \int_{\Lambda} \left[ \frac{1}{2} \left\vert \nabla \phi(x) \right\vert ^{2} + \frac{m^{2}}{2} \phi^{2}(x) - \alpha V(\phi(x)) \right] \mathrm{d}x \right\} \mathcal{D}\phi,
\end{equation}
where $\Xi > 0$ is an appropriate normalizing constant (also called a partition function); $\alpha \neq 0$ is a coupling constant whose magnitude measures the coupling strength between a generic nonlinear differentiable self-interaction $V: \mathbb{R} \mapsto \mathbb{R}$ and the free field;
$m \geqslant 0$ is the mass of the scalar field $\phi$; $\mathcal{D}\phi$ denotes a heuristic ``Lebesgue measure'' on the space $\phi$ lives in.
Such a program is realized by solving the dynamics of $(\varphi_{t})_{t \geqslant 0}$ governed by a nonlinear stochastic partial differential equation (PDE for short) of the Langevin type:
\begin{equation}
    \label{intro-formal-SPDE-Langevin}
    \left\{
        \begin{array}[c]{rcll}
            \mathrm{d}\varphi_{t}(x) & = & - \left[\left( m^{2} - \Delta \right) \varphi_{t}(x) - \alpha V^{\prime} (\varphi_{t}(x)) \right] \mathrm{d}t + \mathrm{d} \mathcal{W}_{t}(x), & (t,x) \in \mathbb{R}_{+} \times \Lambda,\\
            \varphi_{0}(x) & = & \phi(x), & x \in \Lambda
        \end{array}
    \right.
\end{equation}
with initial data $\phi$ and the dynamical noise $\mathcal{W}$ being a cylindrical Wiener process on $L^{2}(\Lambda)$ and satisfying $\mathbb{E}[\mathcal{W}_{t}(x) \mathcal{W}_{s}(y)] = 2 (t \wedge s) \delta_{0} (x - y)$ for $t,s \in [0, +\infty)$ and $x,y \in \mathbb{R}^{d}$.
Actually, (\ref{intro-formal-SPDE-Langevin}) is essentially the stochastic quantization of the classical nonlinear Klein-Gordon equation associated with Gibbs measure (\ref{introduction-Gibbs-measure-formal}), which is derived from the principle of least action. 
Several types of nonlinearity $V$ are concerned, for example, the well-known $\Phi_{d}^{4}$-model for $V(\phi) := \phi^{4}/4$ and the Høegh-Krohn model for $V(\phi) = \exp\{ \beta \phi \}$.

The present paper focuses on the case $d = 2$ and the trigonometric type of nonlinear self-interaction $V(\phi) = \cos (\beta \phi)$ that corresponds to the two-dimensional sine-Gordon model with an inverse temperature $\beta^{2} \in [0,2d]$, 
which has been extensively studied during the past few decades and has gathered significant attention as an interesting toy model due to its importance in constructive quantum field theory as well as to its connection with the two-dimensional gases of charged particles with the Coulomb/Yukawa interaction in the grand canonical ensemble.
(Here, we keep writing $d$ for the dimension in the statements to make the correspondence more transparent with existing/expected results in higher dimensions.)
The basic object of study, Gibbs measure (\ref{introduction-Gibbs-measure-formal}), becomes
\begin{equation}
    \label{introduction-SG-measure-formal}
    \mu_{\mathrm{SG}}(\mathrm{d}\phi) := \Xi^{-1} \exp\left\{ \alpha \int_{\Lambda} \cos(\beta \phi(x)) \mathrm{d}x\right\} \mathbb{P}_{\mathrm{GFF}}(\mathrm{d}\phi),
\end{equation}
where $\mathbb{P}_{\mathrm{GFF}}$ stands for the law of the Gaussian free field (GFF for short, corresponding to the case $V = 0$) on $\Lambda$, being obtained by combining the quadratic terms in (\ref{introduction-Gibbs-measure-formal}) with $\mathcal{D}\phi$ (see \cite{Nelson-1973}). 
(\ref{introduction-SG-measure-formal}) is a prototypical example of non-Gaussian quantum fields and particular interest both as a theory with infinitely many phase transitions for $\beta^{2} \in [0,2d)$ and as a testbed for the understanding of strongly non-polynomial interactions field theories.
(The terminal threshold $\beta^{2} = 2d$ corresponds to the Kosterlitz-Thouless phase transition of the log-gases.)
However, (\ref{introduction-SG-measure-formal}) has no rigorous meaning since either the samples of $\mathbb{P}_{\mathrm{GFF}}$ are known to be genuine distributions on $\Lambda$ whenever $d \geq 2$ so that the cosine potential $\cos(\beta \phi)$ is ill-defined, 
or the sine-Gordon interaction $\int_{\Lambda} \cos(\beta \phi(x)) \mathrm{d}x$ does not exist even in the case $d=1$ where $\Lambda = \mathbb{R}$.

Let us briefly review related works that make sense of (\ref{introduction-SG-measure-formal}) in the constructive literature when $d=2$.
The pioneering constructive treatment of the sine-Gordon model on $\Lambda = \mathbb{R}^{d}$ was given in \cite{Frohlich-Seiler-1976} for $\beta^{2} \in [0,d)$ and small $\alpha$,
and more extensive results are in \cite{Benfatto-Gallavotti-Nicolo-1982,Nicolo-Renn-Steinmann-1986} for the full subcritical regime $\beta^{2} \in [0,2d)$.
Another elegant construction of the sine-Gordon model for $\beta^{2} \in [0,2d)$ relies on the renormalization group method (see \cite{Dimock-Hurd-1993,Dimock-Hurd-2000,Gallavotti-1985} and references therein),
although powerful, which is especially adapted to a translational invariant context for which the renormalization group map is most easily studied. 

Instead of renormalization group techniques, the authors in \cite{Brydges-Kennedy-1987} employed an analysis of the Polchinski flow equations with the Mayer expansions,
and their method was later explored in \cite{Bauerschmidt-Bodineau-CPAM-2021} to prove a logarithmic Sobolev inequality and in \cite{Bauerschmidt-Hofstetter-AOP-2022} to study the maximum of the sine-Gordon quantum field for $\beta^{2} \in [0,3d/2)$.
When $\beta^{2} \in [0,d)$ and $\Lambda = \mathbb{R}^{d}$, utilizing a variational approach \cite{Barashkov-2022}, the Polchinski flow equation is shown to be the Hamilton-Jacobi-Bellman equation
of a stochastic control problem that the minimum of the cost functional is closely related to the Laplace transform of (\ref{introduction-SG-measure-formal}). In terms of the stochastic maximum principle,
this variational approach is extended to $\beta^{2} \in [0, 3d/2)$ in \cite{Gubinelli-Meyer-2024} by introducing a class of forward-backward stochastic differential equations as a counterpart of the Polchinski flow equations.

Besides the above methods, in the spirit of the stochastic quantization method in \cite{Parisi-Wu-1981}, the authors in \cite{Albeverio-Haba-Russo-2001} studied a class of $V$ being the Fourier transform of a complex measure with compact support, which includes the two-dimensional sine-Gordon interaction for sufficiently small $\beta$.
Subsequently, a more natural Langevin dynamics is considered in \cite{Hairer-Shen-2016} as the stochastic quantization equation associated with (\ref{introduction-SG-measure-formal}) for $\beta^{2} \in [0, 4d/3)$ and is extended to the case $\beta^{2} \in [0,2d)$ in \cite{Chandra-Hairer-Shen-2018},
through the regularity structure theory established in \cite{Hairer-Invent-2014} for stochastic nonlinear parabolic PDEs.
Different from the methods just enumerated, a martingale approach is developed in \cite{Lacoin-Rhodes-Vargas-2023-PTRF} for the boundary sine-Gordon model, which covers the full subcritical regime $\beta^{2} \in [0,2d)$ in the case $d=1$ and finite volume $\Lambda$.
Nevertheless, the global understanding of the sine-Gordon model is still very far from being complete from the mathematical perspective.

This paper aims at providing a new probabilistic approach involving the quadratic backward stochastic differential equation (BSDE for short) to construct the two-dimensional sine-Gordon model for $\beta^{2} \in [0,2)$ and bounded domains $\Lambda \subset \mathbb{C}$
and to apply the obtained sine-Gordon measure (with respect to the removal of ultraviolet cutoffs) to specific problems in equilibrium statistical mechanics.
Readers who are not familiar with quadratic BSDEs can refer to the increasingly developed and enriched solvability theory of the scalar- or vector-valued quadratic BSDEs driven by standard Wiener processes in the past two decades (see \cite{Briand-Hu-2008,Kobylanski-AOP-2000,Xing-Zitkovic-AOP-2018} and references therein)
and their wide applications in stochastic linear-quadratic control with random coefficients \cite{Bismut-SICON-1976}, risk-sensitive stochastic optimal control or differential games \cite{Karoui-Hamadene-2003,Hu-Tang-2016}, utility maximization problems \cite{Delbaen-MathFinan-2002,Hu-Imkeller-Muller-AAP-2005}, etc.
To the purpose of this paper, different from the classical quadratic BSDEs just mentioned, we introduce the following quadratic BSDE on the time interval $[0,1]$ driven by a cylindrical Wiener process $W$ such that $W_{t}$ lives in Sobolev space $H_{0}^{-1}(\Lambda)$ for $t \in [0,1]$ and $W_{1}$ is distributed to GFF with the Dirichlet boundary condition (see Section \ref{section-2} for rigorous descriptions):
\begin{equation}
    \label{intro-renormalized-QBSDE}
    \left\{
    \begin{array}[c]{rl}
       \mathrm{d} Y_{t}^{\varepsilon}(\rho) = & - \frac{\alpha}{2} \left\Vert Z_{t}^{\varepsilon}(\rho) \right\Vert _{H_{0}^{-1}(\Lambda)}^{2} \mathrm{d}t + Z_{t}^{\varepsilon}(\rho) \mathrm{d} W_{t}, \quad t \in [0,1),\\
       Y_{1}^{\varepsilon}(\rho) = & \int_{\Lambda_{\varepsilon}} \left[\left[\cos(\beta W_{1}^{\varepsilon}(x))\right]\right] \rho(x) \mathrm{d}x,
    \end{array}
    \right.
\end{equation}
where $\{W_{1}^{\varepsilon}\}_{\varepsilon \in (0,1]}$ is a family of convolution approximations of $W_{1}$; $\left[\left[\cos(\beta W_{1}^{\varepsilon})\right]\right]$ is the Wick-ordered cosine defined on $\Lambda_{\varepsilon}:=\{ x \in \Lambda: \min_{y \in \mathbb{C} \setminus \Lambda} \left\vert y - x \right\vert > 2 \varepsilon \}$ (see Section \ref{section-3});
$\rho$ is an arbitrarily fixed test function. As it is known that only the Wick normalization is needed for $\cos(\beta \phi)$ when $\beta^{2} \in [0,2)$,
the main idea is based on a heuristic observation that, by using the change of variables, (\ref{introduction-SG-measure-formal}) is formally rewritten by $\Xi = \mathbb{E}_{\mathrm{GFF}} \left[ e^{ \alpha \int_{\Lambda} \cos(\beta W_{1}(x)) \mathrm{d}x } \right]$,
where $\cos(\beta W_{1})$ can be viewed as the real part of imaginary multiplicative chaos $e^{i \beta W_{1}}$ that is systematically studied for $\beta^{2} \in [0,d]$ ($d \in \mathbb{N}_{+}$) by Junnila et al. \cite{Webb-AAP-2019} (see also \cite{Lacoin-AAP-2022} and references therein).
In the papers just mentioned, a reasonable normalization of the sine-Gordon interaction $\int_{\Lambda} \cos(\beta W_{1}(x)) \mathrm{d}x$ is $\int_{\Lambda_{\varepsilon}} \left[\left[ \cos(\beta W_{1}^{\varepsilon}(x)) \right]\right] \rho(x) \mathrm{d}x$ 
and the existence of the limit as $\varepsilon \rightarrow 0^{+}$ is proved in the sense of convergence in $\mathbb{P}_{\mathrm{GFF}}$ for any test function $\rho$.
Motivated by this observation, it is natural to consider $\Xi_{\rho, \varepsilon} := \mathbb{E}_{\mathrm{GFF}} \left[ e^{ \alpha \int_{\Lambda_{\varepsilon}} \left[\left[ \cos(\beta W_{1}^{\varepsilon}(x)) \right]\right] \rho(x) \mathrm{d}x } \right]$ as a reasonable normalization of $\Xi$
and then it follows from (\ref{intro-renormalized-QBSDE}) with an application of It\^{o}'s formula (in infinite-dimensional version) to $e^{\alpha Y_{t}^{\varepsilon}(\rho)}$ that
\begin{equation}
    \label{intro-relationship-QBSDE-part-func}
    Y_{0}^{\varepsilon}(\rho) = \frac{1}{\alpha} \log \left( \mathbb{E}_{\mathrm{GFF}} \left[ e^{ \alpha \int_{\Lambda_{\varepsilon}} \left[\left[ \cos(\beta W_{1}^{\varepsilon}(x)) \right]\right] \rho(x) \mathrm{d}x } \right] \right) = \frac{1}{\alpha} \log \Xi_{\rho, \varepsilon},
\end{equation}
provided $\alpha \int_{\Lambda_{\varepsilon}} \left[\left[\cos(\beta W_{1}^{\varepsilon}(x))\right]\right] \rho(x) \mathrm{d}x$ admits an exponential moment. (In fact, it can be shown to be uniformly bounded in Section \ref{section-3}.)
Inspired by (\ref{intro-relationship-QBSDE-part-func}), it is possible to consider a rational approximation $\mu_{\mathrm{SG}}^{\rho, \varepsilon}$ of (\ref{introduction-SG-measure-formal}) by noting that the function
\[
\delta \longmapsto \frac{1}{\alpha} \log \mathbb{E}_{\mathrm{GFF}} \left[  e^{ \alpha \left( \delta F(W_{1}) + \int_{\Lambda_{\varepsilon}} \left[\left[\cos(\beta W_{1}^{\varepsilon}(x))\right]\right] \rho(x) \mathrm{d}x \right) } \right] = Y_{0}^{\varepsilon, \delta F}(\rho)
\]
is differentiable for all bounded and continuous $F$ on $H_{0}^{-1}(\Lambda)$, where $Y^{\varepsilon, \delta F}(\rho)$ is the first component of the solution to (\ref{intro-renormalized-QBSDE}) with the terminal condition being replaced with $\delta F(W_{1}) + \int_{\Lambda_{\varepsilon}} \left[\left[\cos(\beta W_{1}^{\varepsilon}(x))\right]\right] \rho(x) \mathrm{d}x$.
Through investigating the action of any bounded and continuous $F$ on  $\mu_{\mathrm{SG}}^{\rho, \varepsilon}$, we have
\begin{equation*}
    \begin{array}
    [c]{rl}%
    \mu_{\mathrm{SG}}^{\rho, \varepsilon}( F )
    = &  \left.\frac{\mathrm{d}}{\mathrm{d}\delta} \right\vert _{\delta=0} \left( \alpha^{-1} \log \mathbb{E}_{\mathrm{GFF}} \left[  e^{ \alpha \left( \delta F(W_{1}) + \int_{\Lambda_{\varepsilon}} \left[\left[\cos(\beta W_{1}^{\varepsilon}(x))\right]\right] \rho(x) \mathrm{d}x \right) } \right] \right) \\
    = &  \lim\limits_{\delta \rightarrow 0} \frac{\alpha^{-1}}{\delta} \log \left( \frac{\mathbb{E}_{\mathrm{GFF}} \left[  e^{ \alpha \left( \delta F(W_{1}) + \int_{\Lambda_{\varepsilon}} \left[\left[\cos(\beta W_{1}^{\varepsilon}(x))\right]\right] \rho(x) \mathrm{d}x \right) } \right]} {\mathbb{E}_{\mathrm{GFF}} \left[  e^{ \alpha \int_{\Lambda_{\varepsilon}} \left[\left[\cos(\beta W_{1}^{\varepsilon}(x))\right]\right] \rho(x) \mathrm{d}x } \right]} \right) \\
    = &  \lim\limits_{\delta \rightarrow 0} \frac{1}{\delta} \left( Y_{0}^{\varepsilon, \delta F}(\rho) - Y_{0}^{\varepsilon}(\rho) \right) = \widehat{Y}_{0}^{\varepsilon, F}(\rho),
    \end{array}
\end{equation*}
where $\widehat{Y}_{0}^{\varepsilon, F}(\rho)$ is the initial value of a linear BSDE with unbounded coefficients:
\begin{equation}
    \label{intro-variational-eq}
    \left\{
    \begin{array}[c]{rl}
       \mathrm{d} \widehat{Y}_{t}^{\varepsilon, F}(\rho) = & - \alpha \left( Z_{t}^{\varepsilon, 0}(\rho) , \widehat{Z}_{t}^{\varepsilon, F}(\rho) \right) _{H_{0}^{-1}(\Lambda)} \mathrm{d}t + \widehat{Z}_{t}^{\varepsilon, F}(\rho) \mathrm{d} W_{t}, \quad t \in [0,1),\\
       \widehat{Y}_{1}^{\varepsilon, F}(\rho) = & F(W_{1}),
    \end{array}
    \right.
\end{equation}
which we derive from a variational method and whose well-posedness is guaranteed with the aid of some delicate estimates concerning the bounded mean oscillation martingales (see Lemma \ref{sec-3-lemma-3}).
We achieve our ultimate goal by the main result of this paper, Theorem \ref{sec-3-theorem-1}, proving that $\lim_{\varepsilon \rightarrow 0^{+}} \widehat{Y}_{0}^{\varepsilon, F}(\rho) = \mathbb{E}_{\mathrm{GFF}}[ \Gamma(\rho) F(W_{1}) ]$, where $\Gamma(\rho)$ is the stochastic exponential of $\alpha \int_{0}^{1} \overline{Z}_{s}(\rho) \mathrm{d}W_{s}$ such that $\overline{Z}(\rho)$ satisfies the quadratic BSDE
\begin{equation}
    \label{intro-bar-QBSDE}
    \left\{
    \begin{array}[c]{rl}
       \mathrm{d} \overline{Y}_{t}(\rho) = & - \frac{\alpha}{2} \left\Vert \overline{Z}_{t}(\rho) \right\Vert _{H_{0}^{-1}(\Lambda)}^{2} \mathrm{d}t + \overline{Z}_{t}(\rho) \mathrm{d} W_{t}, \quad t \in [0,1),\\
       \overline{Y}_{1}(\rho) = & \langle \cos(\beta W_{1}), \rho \rangle
    \end{array}
    \right.
\end{equation}
with the terminal condition $\langle \cos(\beta W_{1}), \rho \rangle$ being the real part of imaginary multiplicative chaos tested against the given $\rho$, which we obtain by passing the limit as $\varepsilon \rightarrow 0^{+}$ on both sides of (\ref{intro-renormalized-QBSDE}) (see Proposition \ref{sec-3-proposition-2}).
Then the desired sine-Gordon measure is defined by setting $\mu_{\mathrm{SG}}(A) := \mathbb{E}_{\mathrm{GFF}}\left[ \Gamma(\rho) \mathbf{1}_{A} \right] $ for any Borel measurable set $A$ on $H_{0}^{-1}(\Lambda)$.
It should be highlighted that $\Gamma(\rho)$ serves as the Radon-Nikodym derivative since we can prove that $\mathbb{E}_{\mathrm{GFF}}\left[ \Gamma(\rho) \right] = 1$, implying the absolute continuity of $\mu_{\mathrm{SG}}$ with respect to $\mathbb{P}_{\mathrm{GFF}}$ when $\beta^{2} \in [0,d)$, which is consistent with the known result in the existing literature.

To connect our approach with the Parisi-Wu program, a further comparison is necessary between stochastic PDE (\ref{intro-formal-SPDE-Langevin}) and quadratic BSDE (\ref{intro-renormalized-QBSDE}).
Here we compare them in a heuristic way. Consider $\mathcal{W}$ as the canonical process on the path space $\Omega := \{ \omega \in C([0,+\infty); H_{0}^{-1}(\Lambda)): \omega_{0} = 0 \}$ with $\mathfrak{B} (\Omega)$ its Borel $\sigma$-field. One can construct a family of probability measure $\{ \mathbb{P}_{\phi}^{(\mathcal{W})}: \phi \in H_{0}^{-1}(\Lambda) \}$ on $(\Omega, \mathfrak{B} (\Omega))$ such that $\mathcal{W}$ is a cylindrical Wiener process on $L^{2}(\Lambda)$, starting from $\phi$ under $\mathbb{P}_{\phi}^{(\mathcal{W})}$.
This means that a family of linear mappings $(\mathcal{W}_{t})_{t \geqslant 0}$ from $L^{2}(\Lambda)$ to $L^{2}(\Omega)$ satisfying: 
(i) $\forall f \in L^{2}(\Lambda)$, $(\mathcal{W}_{t}(f))_{t \geqslant 0}$ is a real (continuous) Wiener process starting from $0$ under $\mathbb{P}_{0}^{(\mathcal{W})}$; (ii) $\mathbb{E}_{0}^{(\mathcal{W})} [ \mathcal{W}_{t} (f) \cdot \mathcal{W}_{s} (g) ] = 2(t \wedge s) (f,g)_{L^{2}(\Lambda)}$; (iii) $\mathbb{P}_{\phi}^{(\mathcal{W})} (\mathcal{W} \in A) = \mathbb{P}_{0}^{(\mathcal{W})} (\mathcal{W} \in A \mid \mathcal{W}_{0} = \phi)$ for any $A \in \mathfrak{B} (\Omega)$.
When $V$ vanishes, (\ref{intro-formal-SPDE-Langevin}) reduces to a pseudo linear stochastic differential equation:
\begin{equation}
    \label{intro-linear-SPDE}
            \psi_{t} = \phi - \int_{0}^{t} \mathcal{L} \psi_{s} \mathrm{d}s + \mathcal{W}_{t}, \quad  t \in [0,+\infty),
\end{equation}
where the operator $\mathcal{L} := m^{2} \mathrm{Id} - \Delta$ and we omit the space variable $x \in \Lambda$ and the sample $\omega \in \Omega$. Let $\{ e^{- t \mathcal{L}} \}_{t \geqslant 0}$ be the semi-group associated with (\ref{intro-linear-SPDE})
and then (\ref{intro-formal-SPDE-Langevin}) can be also rewritten as an integral equation
\begin{equation}
    \label{intro-integral-eq}
    \varphi_{t} = \psi_{t} + \alpha \int_{0}^{t} e^{- (t - s) \mathcal{L}} V^{\prime} (\varphi_{s}) \mathrm{d}s, \quad t \in [0, +\infty).
\end{equation}
It can be proved that (\ref{intro-linear-SPDE}) has a unique solution $\psi$ being ergodic with $\mathbb{P}_{\mathrm{GFF}}$ as its unique invariant measure.
In the following exposition inspired by the idea in \cite{J-Lasinio-Mitter-1985} concerning the Girsanov transform, we keep writing $V$ for the sine-Gordon interaction because the formal computation is applicable for other nonlinear self-interactions.
Suppose that (\ref{intro-formal-SPDE-Langevin}) admits an ergodic weak solution for any initial condition $\phi \in H_{0}^{-1}(\Lambda)$, namely, 
there is a Markovian family of measures $\{ \mathbb{P}_{\phi}: \phi \in H_{0}^{-1}(\Lambda) \}$ on $(\Omega, \mathfrak{B} (\Omega))$ such that $\mathcal{W}$ is ergodic with $\mathbb{P}_{\mathrm{SG}}(\mathrm{d} \phi) := \Xi_{\mathrm{SG}}^{-1} e^{\alpha \int_{\Lambda} : V( \phi(x) ) : \mathrm{d}x} \mathbb{P}_{\mathrm{GFF}} (\mathrm{d} \phi)$ as its invariant measure and a process defined by
\[
    \widehat{\psi}_{t} := \mathcal{W}_{t} - \alpha \int_{0}^{t} e^{- (t - s) \mathcal{L}} : V^{\prime} (\mathcal{W}_{s}) : \mathrm{d}s, \quad t \in [0, +\infty)
\]
is indistinguishable from the unique solution $\psi$ to (\ref{intro-linear-SPDE}) but not necessarily measurable with respect to it, where $:V:$ (resp. $:V^{\prime}:$) is the Wick-renormalization for $V$ (resp. $V^{\prime}$); $\Xi_{\mathrm{SG}}$ is a normalizing constant.
In this weak sense $(\mathcal{W}, \mathbb{P}_{\phi})$ solves (\ref{intro-integral-eq}) for given $\phi \in H_{0}^{-1}(\Lambda)$ and the corresponding expectation is denoted by $\mathbb{E}_{\phi}$.
For any $t \in [0,+\infty)$, if $\mathbb{E}_{\phi}^{(\mathcal{W})} \left[ e^{ \alpha \int_{0}^{t} : V^{\prime} (\psi_{s}) : \mathrm{d} \mathcal{W}_{s} - \frac{\alpha^{2}}{2} \int_{0}^{t} \left\Vert : V^{\prime} (\psi_{s}) : \right\Vert _{H_{0}^{-1}(\Lambda)}^{2} \mathrm{d}s } \right] = 1$,
then we can define a new probability measure $\overline{\mathbb{P}}_{\phi}^{(\mathcal{W})}$ equivalent to $\mathbb{P}_{\phi}^{(\mathcal{W})}$ with a $H_{0}^{-1}(\Lambda)$-valued Wiener process under $\overline{\mathbb{P}}_{\phi}^{(\mathcal{W})}$ by the Cameron-Martin-Girsanov formula. We deduce from the property of the weak solution that
\begin{equation}
    \label{intro-Girsanov-Fphi}
    \mathbb{E}_{\phi} \left[ F(\mathcal{W}_{t}) \right] = \overline{\mathbb{E}}_{\phi}^{(\mathcal{W})} \left[ F(\psi_{t}) \right] = \mathbb{E}_{\phi}^{(\mathcal{W})} \left[ F(\psi_{t}) e^{ \alpha \int_{0}^{t} : V^{\prime} (\psi_{s}) : \mathrm{d} \mathcal{W}_{s} - \frac{\alpha^{2}}{2} \int_{0}^{t} \left\Vert : V^{\prime} (\psi_{s}) : \right\Vert _{H_{0}^{-1}(\Lambda)}^{2} \mathrm{d}s } \right],
\end{equation}
where $F$ belongs to the space of bounded and continuous functions on $H_{0}^{-1}(\Lambda)$.
We can further rewrite the last line in (\ref{intro-Girsanov-Fphi}) in the Feynman-Kac form $\mathcal{Y}_{0}^{(t,F,\phi)}$ such that $\mathcal{Y}^{(t,F,\phi)}$ is the first unknown variable of the linear BSDE
\begin{equation}
    \label{intro-BSDE-Feynman-Kac}
    \left\{
        \begin{array}[c]{rl}
            \mathrm{d} \mathcal{Y}_{s}^{(t,F,\phi)} = & - \alpha \left( : V^{\prime} (\psi_{s}) :, \mathcal{Z}_{s}^{(t,F,\phi)} \right) _{H_{0}^{-1}(\Lambda)} \mathrm{d}s + \mathcal{Z}_{s}^{(t,F,\phi)} \mathrm{d} \mathcal{W}_{s}, \quad s \in [0,t),\\
            \mathcal{Y}_{t}^{(t,F,\phi)} = & F(\psi_{t}).
         \end{array}
    \right.
\end{equation}
Due to the ergodicity of $( \mathcal{W}, \mathbb{P}_{\phi} )$, we employ (\ref{intro-Girsanov-Fphi})-(\ref{intro-BSDE-Feynman-Kac}) to obtain
\[
    \mathbb{E}_{\mathrm{SG}} (F) = \mathbb{E}_{\mathrm{SG}} \left[ \mathbb{E}_{\phi} \left[ F(B_{t}) \right] \right] =  \mathbb{E}_{\mathrm{SG}} \left[ \mathcal{Y}_{0}^{(t,F,\phi)} \right], \quad \forall t \in [0, +\infty).
\]
On the one hand, it is not hard to see the extreme similarity between (\ref{intro-variational-eq}) and (\ref{intro-BSDE-Feynman-Kac}), whether in their analogous structure or in their role in constructing the limit sine-Gordon measure.
On the other hand, combing (\ref{intro-linear-SPDE}) with (\ref{intro-BSDE-Feynman-Kac}) forms a forward-backward stochastic differential equation that can be regarded as the gradient dynamical system of (\ref{introduction-SG-measure-formal}).
Noting (\ref{intro-relationship-QBSDE-part-func}), the noteworthy difference from (\ref{intro-formal-SPDE-Langevin}) is that (\ref{intro-renormalized-QBSDE}) originates straightforwardly from the normalized log-Laplace transform of Gibbs measure (\ref{introduction-SG-measure-formal})
so that, to some extent, (\ref{intro-formal-SPDE-Langevin}) can be viewed as the primitive dynamics of (\ref{intro-variational-eq}), which is the counterpart of (\ref{intro-BSDE-Feynman-Kac}).

The major novelty of this article is that, without any restriction on the parameter $\alpha$, (\ref{intro-renormalized-QBSDE}) is proved to serve as an effective stochastic quantization equation for (\ref{introduction-SG-measure-formal}) on an arbitrarily simply connected bounded domain $\Lambda \subset \mathbb{C}$ when $\beta^{2} \in [0,2)$,
and the solution of (\ref{intro-bar-QBSDE}) completely characterizes the absolute continuity of our sine-Gordon measure with respect to $\mathbb{P}_{\mathrm{GFF}}$ in such case. Our approach is very different from all the known constructions we enumerated above and brings a new perspective to this extensively investigated classical problem.
It should be emphasized that the authors in \cite{Bauerschmidt-Hofstetter-AOP-2022} employed a distribution-valued BSDE being strongly coupled with a decomposed GFF to construct the sine-Gordon quantum field in the regime $\beta^{2} \in [0, 3d/2)$ when $d=2$ and $\Lambda = \mathbb{T}^{d}$ is the 2D-torus,
where the generator of their BSDE involves the gradient of the solution to the Polchinski flow equation (see \cite{Bauerschmidt-Hofstetter-AOP-2022}, (3.7) on p.477).
Unlike the BSDE in \cite{Bauerschmidt-Hofstetter-AOP-2022}, (\ref{intro-renormalized-QBSDE}) is a scalar-valued BSDE whose generator is purely quadratic growth in the unknown variable $Z^{\varepsilon}(\rho)$, which is much simpler for us to utilize its convex/concave property (resp. $\alpha>0$/$\alpha<0$) to show the well-posedness of (\ref{intro-bar-QBSDE}) as $\varepsilon \rightarrow 0^{+}$ (see Proposition \ref{sec-3-proposition-2}).
Not only that, the coupling relationship with Gaussian field $\{ W_{1}^{\varepsilon}(x): (\varepsilon,x) \in \bigcup _{\varepsilon^{\prime} \in [0,1]} (\{ \varepsilon^{\prime} \} \times \Lambda_{\varepsilon^{\prime}}) \}$ is fully reflected by the terminal condition of (\ref{intro-renormalized-QBSDE}).
Furthermore, compared with the variational approach adopted in \cite{Barashkov-2022,Gubinelli-Meyer-2024},
on the one hand, our method is more flexible in dealing with arbitrarily finite volume $\Lambda$ since many tools in the Fourier analysis are inconvenient to apply in such case; on the other hand, one may not need to use the Bou\'{e}-Dupuis variational formula to transform (\ref{intro-relationship-QBSDE-part-func}) into a stochastic control problem,
which keeps us from solving a complicated forward-backward stochastic Hamiltonian system associated with such kind of stochastic control problem (see \cite{Gubinelli-Meyer-2024}, (2.27) on p.15).

Other contributions of this paper are as follows. For the demand of applications, we may consider more general $\Xi_{\rho, \varepsilon} := \mathbb{E}_{\mathrm{GFF}} \left[ e^{ \alpha \int_{\Lambda_{\varepsilon}} [[\cos(\beta W_{1}^{\varepsilon}(x))]] \rho(x) \mu(\mathrm{d}x) } \right]$, where $\mu(\mathrm{d}x) = \psi(x) \mathrm{d}x$ is a locally finite Borel measure with bounded density $\psi$.
Resorting to several types of the Onsager inequality given in \cite{Webb-AAP-2019}, it is proved by Proposition \ref{sec-4-proposition-3} that the limit $\mathcal{Q}_{\rho, \mu}(\alpha, \beta) := \lim_{\varepsilon \rightarrow 0^{+}} \Xi_{\rho, \varepsilon}$ exists, which plays an pivotal role as a partition function in describing the neutral gas of interacting charged particles with potentials of log-type.
Particularly, when $\psi(x) = \left( \left\vert \varphi^{\prime}(x) \right\vert / 2 \mathrm{Im}\{ \varphi(x) \} \right) ^{\frac{1}{4}}$, by Proposition \ref{sec-4-proposition-1} and Proposition \ref{sec-4-proposition-2} we establish the connection between $\mathcal{Q}_{\rho, \mu} ( 2^{-\frac{1}{2}} \mathcal{C}^{2}, 2^{-\frac{1}{2}} )$ 
and the scaling limit of correlation functions of the spin field for the critical planar XOR-Ising model defined on the dual graph of the square lattice (see \cite{Boutillier-Tiliere-XOR-2014,Webb-AAP-2019}),
where $\mathcal{C}$ is a lattice-dependent constant and $\varphi$ is a conformal bijection from $\Lambda$ to the upper half plane.
Furthermore, since a well-known tool to identify the proper renormalization for log-gases is their sine-Gordon representation,
we utilize the properties of our sine-Gordon measure to obtain the limit of the sine-Gordon representation for the charge distribution under normalization (see Theorem \ref{sec-4-thm-char-func-limit}).
Last but not least, as opposed to the classical quadratic BSDEs driven by ordinary Wiener processes taking values in finite-dimensional Euclidean spaces, the non-triviality of (\ref{intro-renormalized-QBSDE}) itself has a real theoretical appeal and might have other new applications.

The restriction $\beta^{2} \in [0,2)$ is mild enough to be a proper beginning to lay the foundation for our approach.
It is well known that, as $\beta^{2} \rightarrow 2d$, the renormalization theory requires the proof of convergence of more and more auxiliary objects. That is the reason why higher order renormalization schemes need to be introduced, 
such as a sequence of successive thresholds $\{ \beta_{n} := \sqrt{2d(1 - 1/2n)} : n \in \mathbb{N}_{+} \}$ in \cite{Lacoin-Rhodes-Vargas-2023-PTRF} describing the sum of $n$ diverging terms should be subtracted from the Wick-renormalization of the partition function when $\beta \in [\beta_{n}, \beta_{n+1})$ (recall that there $d=1$),
or, in \cite{Hairer-Shen-2016}, the authors build the first-order/second-order auxiliary processes that involve a careful tracking of dipoles/quadrupoles, to implement the regularity structure theory for pursuing the solution to (\ref{intro-formal-SPDE-Langevin}).
We expect our method to be extended to regime $\beta^{2} \in [2,4)$, where the underlying obstacle to overcome is the divergence of $\Xi_{\rho, \varepsilon}$ as $\varepsilon \rightarrow 0^{+}$.
Therefore, a more involved procedure would be required for renormalizing BSDE (\ref{intro-renormalized-QBSDE}) in the spirit of finding the power series expansion of (\ref{intro-relationship-QBSDE-part-func}) with respect to $\alpha$, which attributes to determining certain higher-order variational equations of (\ref{intro-renormalized-QBSDE}) concerning $\alpha$.
It is beyond the objective of this article, so we leave it in our future work.

The outline of this article is as follows. In Section \ref{section-2}, we set up preliminary results that are indispensable for the introduction of quadratic BSDE (\ref{intro-renormalized-QBSDE}).
Section \ref{section-3} sets the stage for approximating (\ref{introduction-SG-measure-formal}) through (\ref{intro-renormalized-QBSDE}) and proving the convergence to our sine-Gordon measure when the ultraviolet cutoffs are removed.
As an application, Section \ref{section-4} discusses the connection between our sine-Gordon measure and the scaling limit of the critical planar XOR-Ising model and investigates the weak convergence of the normalized charge distributions of planar log-gases by its sine-Gordon representation. 

\section{Preliminaries}
\label{section-2}

\subsection{Basic settings and notation}

We start with fixing some general notation, vector spaces, and stochastic processes that will be used in the sequel.
For any $x \in \mathbb{R}^{d}$, $\delta_{x} \in \mathcal{D}^{\prime}(\mathbb{R}^{d})$ denotes the Dirac mass at $x$. The Laplacian operator on $\mathbb{R}^{d}$ is denoted by $\Delta$.
A set $D \subset \mathbb{R}^{d}$ is called a domain if it is open and connected. 
For any domain $D \subset \mathbb{R}^{d}$, $\mathcal{D}_{0}(D)$ denotes the set of compactly supported and $C^{\infty}$ functions in $D$, also known as test functions,
and the support of a test function $\rho$ is denoted by $\mathrm{supp}(\rho) := \overline{\{x \in D: \rho(x) \neq 0 \}}$.
Its dual, the set of distributions on $D$, is denoted by $\mathcal{D}^{\prime}(D)$ and is equipped with the weak-$\ast$ topology on $\mathcal{D}_{0}(D)$. The duality bracket $\langle \cdot , \cdot \rangle: \mathcal{D}^{\prime}(D) \times \mathcal{D}_{0}(D) \mapsto \mathbb{R}$ is a bilinear form.
$L^{2}(D)$ stands for the Hilbert space completion of $\mathcal{D}_{0}(D)$ with respect to the usual inner product $(f,g)_{L^{2}(D)} := \int_{D} f(x) g(x) \mathrm{d}x,\
\forall f,g \in \mathcal{D}_{0}(D)$. 
In particular, we write $L^{2}(D;\mu)$ if the measure of integration concerned on $D$ is a Borel measure $\mu$.
Recall that the Schwartz space $\mathcal{S}(\mathbb{R}^{d})$ consists of infinitely differentiable functions $f: \mathbb{R}^{d} \mapsto \mathbb{R}$ which are rapidly decreasing (also called the Schwartz functions), that is, for all $j \in \mathbb{N}$ and $k = (k_{1},\ldots,k_{d}) \in \mathbb{N}^{d}$,
\[
\left\Vert f \right\Vert _{j,k} := \sup\limits _{x \in \mathbb{R}^{d}} (1 + |x|)^{j} \left\vert D^{k} f(x) \right\vert < +\infty,
\]
where $D^{k} f := \partial _{x_{1}^{k_{1}} \ldots x_{d}^{k_{d}}} f$ denotes the partial derivative of order $k$. It is equipped with the topology generated by the family of semi-norms $\{ \left\Vert \cdot \right\Vert _{j,k}: j \in \mathbb{N}, k \in \mathbb{N}^{d} \}$.
The space of tempered distributions $\mathcal{S}^{\prime}(\mathbb{R}^{d})$ is the topological dual of $\mathcal{S}(\mathbb{R}^{d})$, which is equipped with the weak-$\ast$ topology on $\mathcal{S}(\mathbb{R}^{d})$.
The Borel $\sigma$-field corresponding to this topology on $\mathcal{S}^{\prime}(\mathbb{R}^{d})$ is denoted by $\mathfrak{B}(\mathcal{S}^{\prime}(\mathbb{R}^{d}))$.

Let $C\left([0,1];\mathbb{R}\right)$ be the collection of continuous real-valued functions defined on $[0,1]$
and $\mathfrak{B}(C\left([0,1];\mathbb{R}\right))$ be the Borel $\sigma$-algebra generated by all open sets
under the uniform norm-topology on $C\left([0,1];\mathbb{R}\right)$, i.e. $\left\Vert \varphi \right\Vert _{C\left([0,1];\mathbb{R}\right)} := \sup_{t \in [0,1]} \left\vert \varphi(t) \right\vert$.
Define
\[
    B_{t}(\varphi) := \varphi(t), \quad \forall t \in [0,1], \ \varphi \in C\left([0,1];\mathbb{R}\right)
\]
the canonical process and let $\mathbb{P}_{0}$ be the Wiener measure on $C\left([0,1];\mathbb{R}\right)$ such that $(B_{t})_{t \in [0,1]}$ is a standard Wiener process under $\mathbb{P}_{0}$.
Let $\Omega = \{ \omega = (\omega^{1}, \omega^{2}, \ldots) : \omega^{k} \in C\left([0,1];\mathbb{R}\right), \omega_{0}^{k} = 0, \forall k \in \mathbb{N}_{+} \}$
and $\Pi^{k} : \omega \mapsto \omega^{k}, k \in \mathbb{N}_{+}$ be the projection mappings on $\Omega$.
By $\mathcal{F} := \sigma ( \bigcup _{n=1}^{\infty} \mathcal{A}_{n})$ we denote the $\sigma$-algebra generated by all the measurable rectangle sets in $\Omega$,
where $\mathcal{A}_{n} := \{ \bigcap _{k=1}^{n} \Pi_{k}^{-1} (A_{k}) : A_{n} \in \mathfrak{B}(C\left([0,1];\mathbb{R}\right)) \}$
for any $n \in \mathbb{N}_{+}$. 
According to the Kolmogorov extension theorem, there exists a probability measure $\mathbb{P}$ on $(\Omega, \mathcal{F})$ satisfying
$\mathbb{P} (\bigcap _{k=1}^{n} \Pi_{k}^{-1} (A_{k})) = \mathbb{P}_{0} (A_{1}) \times  \cdots \times \mathbb{P}_{0} (A_{n}), \ \forall \{ A_{k} \}_{k=1}^{n} \subset \mathfrak{B}(C\left([0,1];\mathbb{R}\right))$.
The mathematical expectation corresponding to $\mathbb{P}$ is denoted by $\mathbb{E}$.
Moreover, for any $n \in \mathbb{N}_{+}$ and $\{ k_{1}, \ldots, k_{n} \} \subset \mathbb{N}_{+}$, it can be verified that the $\mathbb{R}^{n}$-valued process $(B_{t}^{(k_{1})},\ldots, B_{t}^{(k_{n})})_{t \in [0,1]}$
is an $n$-dimensional standard Wiener process under $\mathbb{P}$, where $B_{t}^{k_{j}} := B_{t} \circ \Pi^{k_{j}}, \ t \in [0,1], \ j = 1, \ldots, n$.

For any measure space triple $(E, \mathcal{E}, \mu)$ and any $\mu$-integrable function $f$, we make the conventional notation $\mu(F):= \int f \mathrm{d}\mu$.
Let the probability space $(\Omega, \mathcal{F}, \mathbb{P} )$ be equipped with a filtration $\mathbb{F} := \{ \mathcal{F}_{t} \subset \mathcal{F} : t \in [0,1]\}$ satisfying the usual condition.
In the following definitions, all random variables and stochastic processes will be defined 
on the filtered probability space $(\Omega, \mathcal{F}, \mathbb{F}, \mathbb{P} )$,
where all the concepts of measurability for stochastic processes (e.g. predictability etc.) refer to $\mathbb{F}$.
As it will be seen in the rest of this paper, the concrete forms of $\mathbb{F}$ will be indicated in specific cases needed later.
For any separable Banach space $\Theta$ with the norm $\left\Vert \cdot \right\Vert _{\Theta}$ and its dual space $\Theta^{\ast}$, let $\mathfrak{B}(\Theta)$ be the Borel $\sigma$-field generated by the open sets on $\Theta$,
and $\sigma(\Theta^{\ast})$ denotes the smallest $\sigma$-algebra such that all the elements in $\Theta^{\ast}$ are measurable. Then it can be shown that $\mathfrak{B}(\Theta) = \sigma(\Theta^{\ast})$. (The separability of $\Theta$ avoids the occurrence of strict inclusion $\sigma(\Theta^{\ast}) \subsetneq \mathfrak{B}(\Theta)$,
which is crucial to tieing up the measurability of Banach-valued functions with approximation by simple functions, and the application of standard tools in functional analysis.)
For any $p,q \in [1,+\infty)$ and any $\sigma$-field $\mathcal{G} \subset \mathcal{F}$, we introduce:

\begin{itemize}
    \item $\mathcal{B}_{b}(\Theta)$: the space of bounded Borel measurable functions on $\Theta$.
    
    \item $\mathcal{C}_{b}(\Theta)$: the space of bounded continuous functions on $\Theta$.

    \item $L^{p}(\mathcal{G};\Theta)$: the space of $\mathcal{G}$-measurable, $\Theta$-valued random variables $\xi$ such that
    $$
    \left\Vert \xi \right\Vert _{L^{p}(\mathcal{G};\Theta)}:=\left(\mathbb{E}\left[  \left\Vert \xi\right\Vert _{\Theta}^{p}\right]\right)^{\frac{1}{p}}  <+\infty;
    $$

    \item $L^{\infty}(\mathcal{G};\Theta)$: the space of
    $\mathcal{G}$-measurable, $\Theta$-valued random variables $\xi$
    such that 
    \[
        \left\Vert \xi \right\Vert _{L^{\infty}(\mathcal{G};\Theta)}:= \underset{\mathbb{P}}{\mathrm{ess~sup}} \left\Vert \xi  \right\Vert _{\Theta} <+\infty;
    \]

    \item $\mathcal{M}_{\mathbb{F}}^{p,q}([0,1];\Theta)$: the space of predictable, $\Theta$-valued processes $X$ on
    $[0,1]$ such that
    \[
    \left\Vert X \right\Vert _{\mathcal{M}_{\mathbb{F}}^{p,q}([0,1];\Theta)}:= \left\Vert \left( \int_{0}^{1}\left\Vert X_{t}\right\Vert _{\Theta}^{p}dt \right)^{\frac{1}{p}} \right\Vert _{L^{q}(\mathcal{F}_{1};\mathbb{R})} <+\infty,
    \]
    and particularly we write $\mathcal{H}_{\mathbb{F}}^{q}([0,1];\Theta) := \mathcal{M}_{\mathbb{F}}^{2,q}([0,1];\Theta)$;

    \item $L_{\mathbb{F}}^{\infty}([0,1];\Theta)$: the space of predictable, $\Theta$-valued processes $X$ on $[0,1]$ such that
    \[
    \left\Vert X \right\Vert _{L_{\mathbb{F}}^{\infty}([0,1];\Theta)}:=  \underset{\mathbb{P}}{\mathrm{ess~sup}} \left( \underset{\lambda}{\mathrm{ess~sup}} \left\Vert X_{t} \right\Vert _{\Theta} \right) <+\infty,
    \]
    where $\lambda$ denotes the Lebesgue measure on $[0,1]$;

    \item $\mathcal{S}_{\mathbb{F}}^{p}([0,1];\Theta)$: the space of continuous processes $X \in \mathcal{H}_{\mathbb{F}}^{p}([0,1];\Theta)$ such that
    \[
    \left\Vert X \right\Vert _{\mathcal{S}_{\mathbb{F}}^{p}([0,1];\Theta)}:= \left\Vert \sup\limits_{t\in[0,1]}\left\Vert X_{t}\right\Vert _{\Theta} \right\Vert _{L^{p}(\mathcal{F}_{1};\mathbb{R})} <+\infty;
    \]
    \item $\mathcal{S}_{\mathbb{F}}^{\infty}([0,1];\Theta)$: the space of continuous processes $X \in L_{\mathbb{F}}^{\infty}([0,1];\Theta)$ such that
    \[
        \left\Vert X \right\Vert _{\mathcal{S}_{\mathbb{F}}^{\infty}([0,1];\Theta)}:= \left\Vert \sup\limits_{t\in[0,1]}\left\Vert X_{t}\right\Vert _{\Theta} \right\Vert _{L^{\infty}(\mathcal{F}_{1};\mathbb{R})} <+\infty.
    \]
\end{itemize}

Here we list some notations and results of BMO martingales, which will be used in this paper. We refer readers to \cite{BMO-Book} and the references therein for more details.
Let $\mathrm{BMO}$ denote the Banach space completion of continuous, $\mathbb{F}$-local martingales $M \in \mathcal{H}_{\mathbb{F}}^{2}([0,1];\mathbb{R})$ equipped with the norm
    \[
    \left\Vert M \right\Vert _{\mathrm{BMO}} := \sup\limits_{\tau} \left\Vert \sqrt {\mathbb{E} \left[ \ll M \gg _{1} - \ll M \gg _{\tau} \mid \mathcal{F}_{\tau} \right] } \right\Vert _{L^{\infty}(\mathcal{F}_{1};\mathbb{R})} < +\infty,
    \]
where $\ll M \gg$ denotes the quadratic variation of $M$, and the supremum is taken over all $[0,1]$-valued stopping times $\tau$. 
Actually, one can replace $\tau$ with all deterministic times $t \in [0,1]$ in the above definition. Moreover, we write $\mathrm{BMO}(\mathbb{Q})$ for any probability measure $\mathbb{Q}$ defined on $(\Omega, \mathcal{F})$
whenever it is necessary to indicate the underlying probability.

The very important feature of BMO martingales is the following: the Dol\'{e}ans-Dade exponential of a continuous, $\mathbb{F}$-local martingale $M$, denoted by $\mathcal{E}(M):=\exp \{ M - \ll M \gg /2 \}$,
is a uniformly integrable martingale if $M$ belongs to $\mathrm{BMO}$ (see Theorem 2.3 in the monograph of Kazamaki \cite{BMO-Book}).
Moreover, $\mathcal{E}(M)$ satisfies a reverse H\"{o}lder inequality. Let $\kappa$ be the function defined on $(1, +\infty)$ by setting
\begin{equation}
    \label{sec-2-func-Psi}
    \kappa(x) := \left( 1 + \frac{1}{p^{2}} \log \frac{2p-1}{2(p-1)} \right) ^{\frac{1}{2}} - 1;
\end{equation}
$\kappa$ is non-increasing with $\lim_{p \rightarrow 1^{+}} \kappa(p) = +\infty$ and $\lim_{p \rightarrow +\infty} \kappa(p) = 0$. Let $p_{a}$ be such that $\kappa(p_{a}) = a$.
Then, for each $p \in (1,p_{a})$, and for all $[0,1]$-valued stopping times $\tau$,
\begin{equation}
    \label{sec-2-reverse-Holder-ineq}
    \mathbb{E} \left[ \mathcal{E}(M_{1}) \mid \mathcal{F}_{\tau} \right] \leqslant K \left( p, a \right) \mathcal{E}(M_{\tau}), \quad \mathbb{P}-a.s.,
\end{equation}
where the constant $K \left( p, a \right)$ can be chosen depending only on $p$ and $a = \left\Vert M \right\Vert _{\mathrm{BMO}}$, e.g.
\begin{equation}
    \label{sec-2-reverse-Holder-coeff}
    K \left( p, a \right) = 2 \left( 1 - \frac{2(p-1)}{2p-1} \exp \{ p^{2} (a^{2} + 2a) \} \right)^{-1}.
\end{equation}
   
\subsection{Cylindrical Wiener processes and 2D Dirichlet Gaussian free fields}

Now we focus on the two-dimensional case by fixing $d = 2$. Let $\Lambda \subset \mathbb{C}$ be a bounded and simply connected domain with its boundary denoted by $\partial \Lambda$.
It is well known that the following spectral problem with the Dirichlet boundary condition:
\[
    \left\{
        \begin{array}[c]{rcll}
            - \Delta u(x) & = & \lambda^{-1} u(x), & x \in \Lambda,\\
            u(x) & = & 0, & x \in \partial \Lambda
        \end{array}
    \right.
\]
defines a self-adjoint and compact operator $(- \Delta)^{-1}$ being strictly positive on $L^{2}(\Lambda)$. 
The Hilbert-Schmidt theorem yields a sequence of non-increasing and strictly positive real numbers $\{ \lambda_{k} \}_{k \in \mathbb{N}_{+}}$
and an orthonormal basis $\{ e_{k} \}_{k \in \mathbb{N}_{+}}$ for $L^{2}(\Lambda)$ such that, for each $k \in \mathbb{N}_{+}$, $\lambda_{k}$ is the eigenvalue of $(- \Delta)^{-1}$ 
and $e_{k}$ is the eigenfunction corresponding to $\lambda_{k}$. 
Without any assumption of smoothness on $\partial \Lambda$, it follows from the Weyl asymptotic formula that
\begin{equation}
    \label{sec-2-Weyl-law}
\lim\limits _{k \rightarrow \infty} \lambda_{k} k = 4\pi \left\vert \Lambda \right\vert ^{-1},
\end{equation}
where $\left\vert \Lambda \right\vert$ denotes the volume of $\Lambda$ (References and much more precise estimates on the growth of $\lambda_{k}$ are given in \cite{Netrusov-Safarov-2005}).

\begin{remark}
    If one considers the massive sine-Gordon model with a mass $m>0$ (which is closely related to the 2D-Yukawa gas), then the Dirichlet Laplacian $-\Delta$ on $\Lambda$ can be replaced with the operator $m^{2}\mathrm{Id} - \Delta$
    in the above spectral problem. By Theorem 3.1 and Remark 3.2 in \cite{Birman-Solomyak-1971}, the estimates on the growth of the corresponding eigenvalues $\{ \lambda_{k}^{(m)} \}_{k \in \mathbb{N}_{+}}$ can be obtained such that
    $\lambda_{k}^{(m)} \leqslant \lambda_{k} \leqslant C k^{-1}$ for each $k \in \mathbb{N}_{+}$, where the constant $C>0$ depends only on $\left\vert \Lambda \right\vert$.
\end{remark}

The above leads us to define Sobolev spaces of real index on $\Lambda$, denoted by $H_{0}^{s}(\Lambda)$ for $s \in \mathbb{R}$, to be the Hilbert space completion of
\[
(f,g)_{H_{0}^{s}(\Lambda)} := \sum\limits _{k=1}^{\infty} \lambda_{k}^{-s} (f,e_{k})_{L^{2}(\Lambda)} (g,e_{k})_{L^{2}(\Lambda)}, \quad \forall f,g \in \mathcal{D}_{0}(\Lambda).
\]
It is easy to check that $\{ \lambda_{k}^{\frac{s}{2}} e_{k} \} _{k=1}^{\infty}$ forms an orthonormal basis of $H_{0}^{s}(\Lambda)$.

\begin{remark}
    \label{sec-2-remark-2}
The above series does converge for any $f,g \in \mathcal{D}_{0}(\Lambda)$: this can be seen by applying Cauchy-Schwartz, utilizing that $\mathcal{D}_{0}(\Lambda) \subset L^{2}(\Lambda)$, 
and that all derivatives of test functions in $\mathcal{D}_{0}(\Lambda)$ are again elements of $\mathcal{D}_{0}(\Lambda)$, with $(\Delta f,e_{k})_{L^{2}(\Lambda)} = - \lambda_{k}^{-1} (f,e_{k})_{L^{2}(\Lambda)}$
for $f \in \mathcal{D}_{0}(\Lambda)$ and $k \in \mathbb{N}_{+}$.
\end{remark}

Here some basic facts about $H_{0}^{s}(\Lambda)$ are listed below. We refer readers to \cite{Berestycki-GFF-2024, Sheffield-Sun-surveys-2016} for more details.
\begin{itemize}
\item When $s = 0$ the above space is equivalent, by definition, to $L^{2}(\Lambda)$. 

\item When $s = 1$, $H_{0}^{1}(\Lambda)$ agrees with the standard Sobolev space which is the completion of $\mathcal{D}_{0}(\Lambda)$ with respect to the Dirichlet inner product
$(f,g)_{\nabla} := (\nabla f,\nabla g)_{L^{2}(\Lambda)}$.

\item When $s > 0$, it is simple to check that $H_{0}^{s}(\Lambda) \subset L^{2}(\Lambda)$, and that $f \in L^{2}(\Lambda)$
is an element of $H_{0}^{s}(\Lambda)$ if and only if $\sum\limits _{k=1}^{\infty} \lambda_{k}^{-s} (f,e_{k})_{L^{2}(\Lambda)}^{2} < +\infty$.

\item When $s < 0$, $H_{0}^{s}(\Lambda)$ can be identified with a subspace of $\mathcal{D}^{\prime}(\Lambda)$, and is the dual space of $H_{0}^{-s}(\Lambda)$.
It can be proved that $\left\vert \langle \phi, \rho \rangle \right\vert \leqslant \left\Vert \phi \right\Vert _{H_{0}^{s}(\Lambda)} \left\Vert \rho \right\Vert _{H_{0}^{-s}(\Lambda)}$ for any $\phi \in H_{0}^{s}(\Lambda), \rho \in \mathcal{D}_{0}(\Lambda)$.
\end{itemize}

The separable Hilbert space $H_{0}^{1}(\Lambda)$ plays a key role in the definition of a cylindrical Wiener process $W$ on $L^{2}(\Lambda)$ and the corresponding stochastic integral with respect to $W$.
Recall the definition of $\{(B_{t}^{k})_{0 \leqslant t \leqslant 1}: k \in \mathbb{N}_{+} \}$ in Subsection 2.1. 
For each $t \in [0,1]$, let $\mathcal{F}_{t}^{k} := \sigma(B_{s}^{k}: 0 \leqslant s \leqslant t) \vee \mathcal{N}$ and $ \mathcal{F}_{t}^{B} := \sigma( \bigcup _{k=1}^{\infty} \mathcal{F}_{t}^{k} )$,
where $\mathcal{N}$ is smallest $\sigma$-algebra that contains all $\mathbb{P}$-null sets.
The filtration $\mathbb{F}^{B} := \{ \mathcal{F}_{t}^{B} : t \in [0,1]\}$ is continuous since each $\mathcal{F}_{t}^{k}$ is continuous with respect to $t$.

\begin{proposition}
    \label{sec-2-proposition-1}
    For any $s < 0$, there exists a $H_{0}^{s}(\Lambda)$-valued, $\mathbb{F}^{B}$-adapted process $W$ such that
    \begin{equation}
        \label{sec-2-W-L2-def}
    W_{t}(f) := \sum\limits _{k=1}^{\infty} (2\pi \lambda_{k})^{\frac{1}{2}} B_{t}^{k} (f, e_{k})_{L^{2}(\Lambda)}
    \end{equation}
    is a real-valued Wiener process for arbitrarily $f \in L^{2}(\Lambda)$, and
    \begin{equation}
        \label{sec-2-W-covariance}
        \mathbb{E} \left[ W_{t}(f) W_{r}(g) \right] = 2\pi (t \wedge r) \left((- \Delta)^{-1} f,g \right)_{L^{2}(\Lambda)}, \quad \forall f, g \in L^{2}(\Lambda).
    \end{equation}
    $W$ is also called a cylindrical Wiener process on $L^{2}(\Lambda)$.
\end{proposition}

\begin{proof}
$W$ will be constructed in a larger Hilbert space than $L^{2}(\Lambda)$. 
To this end, for any $s < 0$ the embedding of $H_{0}^{1}(\Lambda)$ into $H_{0}^{s}(\Lambda)$ is Hilbert-Schmidt since by (\ref{sec-2-Weyl-law}) we deduce
\[
    \sum\limits _{k=1}^{\infty} \left\Vert \lambda_{k}^{\frac{1}{2}} e_{k} \right\Vert _{H_{0}^{s}(\Lambda)}^{2} = \sum\limits _{k=1}^{\infty} \lambda_{k} \left\Vert e_{k} \right\Vert _{H_{0}^{s}(\Lambda)}^{2} 
    = \sum\limits _{k=1}^{\infty} \lambda_{k}^{1 - s} < +\infty.
\]
According to Proposition 4.11 in \cite{Da-Prato-Book}, the Karhunen-Lo\`{e}ve expansion 
\begin{equation}
    \label{sec-2-W-partial-sum}
    W_{t} := \sum\limits _{k=1}^{\infty} (2\pi \lambda_{k})^{\frac{1}{2}} B_{t}^{k} e_{k}, \quad t \in [0,1]
\end{equation}
defines a $H_{0}^{s}(\Lambda)$-valued Wiener process such that for any $f \in L^{2}(\Lambda)$ the process $W(f)$ defined in (\ref{sec-2-W-L2-def}) is a centered Gaussian 
process with the covariance (\ref{sec-2-W-covariance}).
\end{proof}

By Proposition \ref{sec-2-proposition-1}, unless indicated, we regard $W$ being constructed as a $H_{0}^{-1}(\Lambda)$-valued process in what follows.
Clearly, it follows from (\ref{sec-2-W-L2-def}) that
\[
\mathcal{F}_{t}^{W} := \sigma\left( W_{s}(f): f \in L^{2}(\Lambda), 0 \leqslant s \leqslant t \right) \vee \mathcal{N} = \mathcal{F}_{t}^{B}, \quad t \in [0,1].
\]
Hence we have $\mathbb{F}^{W} := \{ \mathcal{F}_{t}^{W} : t \in [0,1]\} = \mathbb{F}^{B}$.

Following the route proposed in monograph \cite{Berestycki-GFF-2024}, Chapter 2, we introduce a function $G_{\Lambda} : \Lambda \times \Lambda \mapsto \mathbb{R}$,
such that $W_{1}$ is identically distributed to a Dirichlet Gaussian free field (GFF in short) on $\Lambda$ equipped with $G_{\Lambda}$ as its covariance kernel.
Let $(\tilde{B}, \{\mathbb{P}_{x}\}_{x \in \mathbb{R}^{2}})$ be a two-dimensional Brownian family (see \cite{K-S-Book}, Definition 5.8) with speed ``two'' 
(that is, under each probability measure $\mathbb{P}_{x}$, the process $(\tilde{B}_{t})_{t \geqslant 0}$ is a two-dimensional Brownian motion starting at $x$ with quadratic variation $\langle \tilde{B} \rangle _{t} = 2t \mathrm{I}_{2 \times 2}$ for $t \geqslant 0$),
which is supposed to be independent of $\{(B_{t}^{k})_{0 \leqslant t \leqslant 1}: k \in \mathbb{N}_{+} \}$.
The transition probability of $\tilde{B}$ under $\mathbb{P}_{x}$ is given by
\[
    \mathbb{P}_{x}(\tilde{B}_{t} \in \mathrm{d}y) = \frac{1}{4 \pi t} \exp \left\{ - \frac{\left\vert x - y \right\vert ^{2}}{4t} \right\} \mathrm{d}y, \quad \forall t \geqslant 0, \ x,y \in \mathbb{R}^{2}.
\]
Fixing $x \in \Lambda$, the law of $\tilde{B}$ (starting from $x$) killed when leaving $\Lambda$ is described by
\begin{equation}
    \label{sec-2-law-killed-B}
    \begin{array}[c]{rl}
        \mathbb{E}_{x} \left[ \int_{0}^{\tau_{\Lambda}} f(\tilde{B}_{t}) \mathrm{d}t \right]
        = & \int_{0}^{+\infty} \mathbb{E}_{x} \left[ f(\tilde{B}_{t}) \mathbf{1}_{ \{\tau_{\Lambda} > t \} } \right] \mathrm{d}t \\
        = & \int_{0}^{+\infty} \int_{\mathbb{R}^{2}} f(y) p_{t}^{\Lambda}(x,y) \mathrm{d}y \mathrm{d}t \\
        = & \int_{\mathbb{R}^{2}} f(z) \left( \int_{0}^{+\infty} p_{t}^{\Lambda}(x,y) \mathrm{d}t \right) \mathrm{d}y
    \end{array}
\end{equation}
for any non-negative Borel function $f$ on $\mathbb{R}^{2}$, where the exiting time $\tau_{\Lambda} := \inf \{ t > 0: \tilde{B}_{t} \notin \Lambda \}$ and
\[
p_{t}^{\Lambda}(x,y) := \frac{1}{4 \pi t} \exp \left\{- \frac{\left\vert x - y \right\vert ^{2}}{4t} \right\} \mathbb{P}_{x \rightarrow y;t}(\tau_{\Lambda} > t)
\] 
with $\mathbb{P}_{x \rightarrow y;t}$ denoting the law of $\tilde{B}$ conditionally given $\tilde{B}_{0} = x$ and $\tilde{B}_{t} = y$ (actually a ``speed two'' Brownian bridge of duration $t$ from $x$ to $y$).
The so-called Green function $G_{\Lambda}$ is defined by the time integral in brackets in the last line of (\ref{sec-2-law-killed-B}).

\begin{definition}
\label{sec-2-def-Green-func}
The Green function $G_{\Lambda}(x,y)$ is defined by
\[
    G_{\Lambda}(x,y) = 2\pi \int_{0}^{+\infty} p_{t}^{\Lambda}(x,y) \mathrm{d}t, \quad \forall x,y \in \Lambda, \ x \neq y.
\]
\end{definition}

It should be pointed out that $G_{\Lambda}(x,y) < +\infty$ as long as $x \neq y$ and $\Lambda$ is regular, that is, $\partial \Lambda \neq \varnothing$ and $\mathbb{P}_{x}(\tau_{\Lambda} = 0) = 1$ for all $x \in \partial \Lambda$;
see, for example Lemma 2.32 in \cite{Book-Lawler-2005}.

\begin{remark}
    Corresponding to the sine-Gordon model with a mass $m > 0$, one can define the massive Green function $G_{\Lambda}^{(m)}$ by replacing $p_{t}^{\Lambda}(x,y)$ with $p_{t}^{\Lambda}(x,y) e^{-m^{2}t}$
    which is the transition density of $\tilde{B}$ killed when leaving $\Lambda$ at the rate $m^{2}$.
\end{remark}

\begin{remark}
    We call the attention of readers to the fact that the normalization of the Green function is a little arbitrary. Here we have chosen to normalize it so that $G_{\Lambda}$ is the inverse of (minus) Dirichlet Laplacian on $\Lambda$ multiplied by $2\pi$,
    which is consistent with the standard set up for (real and imaginary) Gaussian multiplicative chaos or in papers on Liouville quantum gravity where the Green function is often normalized so that it blows up like $-\log |x-y|$
    (see \cite{Berestycki-GFF-2024,Webb-AAP-2019,Lacoin-AAP-2022,Rhodes-Vargas-survey-2014} and references therein).
\end{remark}

\begin{remark}
Let $\Lambda^{\prime} \subset \mathbb{C}$ be a domain (not necessarily bounded and simply connected).
According to Proposition 1.14 in \cite{Berestycki-GFF-2024}, if $\Lambda$, $\Lambda^{\prime}$ are regular, then one can prove that the conformal invariance
\[
    G_{\Lambda^{\prime}}(\mathrm{T}(x), \mathrm{T}(y)) = G_{\Lambda}(x,y), \quad \forall x,y \in \Lambda
\]
holds for any conformal isomorphism $\mathrm{T}: \Lambda \mapsto \Lambda^{\prime}$.
In 2D-case, the conformal invariance of $G_{\Lambda}$ is crucial
since it often suffices to prove some desired property in a concrete domain (where we have explicit formulae, such as the upper half plane),
and use conformal invariance to deduce the desired property in an arbitrary regular and simply connected domain (see Example 1.17 in \cite{Berestycki-GFF-2024} for the case $\Lambda$ being the unit disc).
\end{remark}

Some basic and fundamental properties of $G_{\Lambda}(x,y)$ are stated below.

\begin{proposition}[\cite{Berestycki-GFF-2024}, Proposition 1.18, Theorem 1.23]
    \label{sec-2-proposition-Green-func}
    Suppose that $\Lambda$ is regular. Then, for any $x \in \Lambda$,
\begin{enumerate}
\item $G_{\Lambda}(x,y) \rightarrow 0$ as $y \rightarrow y_{0} \in \partial \Lambda$;

\item $G_{\Lambda}(x,y) = - \log |x - y| + g_{\Lambda}(x,y)$, for some $g_{\Lambda}$ being continuous and bounded in $\Lambda \times \Lambda$;

\item $G_{\Lambda}(x,\cdot)$ is harmonic in $\Lambda \setminus \{x\}$; and as a distribution $\Delta G_{\Lambda}(x,\cdot) = -\delta_{x}(\cdot)$.
\end{enumerate}
\end{proposition}

By the first property in the Proposition \ref{sec-2-proposition-Green-func}, we find it convenient to extend $G_{\Lambda}$ to $\mathbb{R}^{2} \times \mathbb{R}^{2}$ by setting $G_{\Lambda}(x,y) = 0$ whenever $(x,y) \notin \Lambda \times \Lambda$.
The proposition below indicates that $W_{1}$ can be viewed as a GFF with the Dirichlet boundary condition, whose law is supported in $H_{0}^{-1}(\Lambda)$, with the covariance kernel $G_{\Lambda}$.
\begin{proposition}
    \label{sec-2-prop-GFF-Green}
    Suppose that $\Lambda$ is regular. Then, for any $f, g \in H_{0}^{1}(\Lambda)$, we have
    \[
        \mathbb{E} \left[ W_{1}(f) W_{1}(g) \right] = \int_{\Lambda \times \Lambda} G_{\Lambda}(x,y) f(x) g(y) \mathrm{d}x \mathrm{d}y.
    \]
\end{proposition}

\begin{proof}
Let $\mu$ be a signed Radon measure supported in $\Lambda$ and $W_{1}^{n}$ be the $n$-th the partial sum of (\ref{sec-2-W-partial-sum}) when $t = 1$. 
By Theorem 1.44 in \cite{Berestycki-GFF-2024}, the integration of $W_{1}^{n}$ with respect to $\mu$ defined by
\[
\mu(W_{1}^{n}) := \sum\limits _{k=1}^{n} (2\pi \lambda_{k})^{\frac{1}{2}} B_{1}^{k} \int_{\Lambda} e_{k}(x) \mu(\mathrm{d}x)
\]
converges in $L^{2}(\mathcal{F}_{1}^{B};\mathbb{R})$ to a Gaussian random variable with variance
$\int_{\Lambda \times \Lambda} G_{\Lambda}(x,y) \mu(\mathrm{d}x)\mu(\mathrm{d}y)$ as $n \rightarrow \infty$.
It follows from (\ref{sec-2-W-L2-def}) that the law of this limit is uniquely defined, and coincides with the Dirichlet GFF when its index set is restricted to the collection of measures whose elements are absolutely continuous with respect to 
the Lebesgue measure such that the density functions belong to $H_{0}^{1}(\Lambda)$.
\end{proof}

\begin{remark}
    One could as well consider different versions of GFF and other types of Green functions with various boundary conditions, say, GFF with Neumann boundary conditions and Neumann Green functions (see \cite{Berestycki-GFF-2024}, Chapter 6 for a detailed discussion)
    leading to log-gases with various boundary conditions.
\end{remark}

Let $L_{2}\left( H_{0}^{1}(\Lambda); \mathbb{R} \right)$ denote the space of Hilbert-Schmidt operators from $H_{0}^{1}(\Lambda)$ to $\mathbb{R}$ endowed with the Hilbert-Schmidt norm,
which is also a separable Hilbert space.
According to the definition of Hilbert-Schmidt norm $\left\Vert \cdot \right\Vert _{L_{2}\left( H_{0}^{1}(\Lambda); \mathbb{R} \right)}$, we observe that $L_{2}\left( H_{0}^{1}(\Lambda); \mathbb{R} \right)$ can be identified with $H_{0}^{-1}(\Lambda)$ the dual space of $H_{0}^{1}(\Lambda)$.
Given an element $\Phi \in \mathcal{H}_{\mathbb{F}^{W}}^{2}([0,1];H_{0}^{-1}(\Lambda))$, one can define the It\^{o}-type stochastic integral for $\Phi$ with respect to $W$
which is denoted by $\{ (\Phi \cdot W)_{t}, t \in [0,1] \} $ or $\{ \int_{0}^{t} \Phi_{s} \mathrm{d} W_{s}, t \in [0,1] \} $.

\begin{proposition}
    \label{sec-2-proposition-3}
Suppose that $\Lambda$ is regular. Then, for any $\Phi \in \mathcal{H}_{\mathbb{F}^{W}}^{2}([0,1];H_{0}^{-1}(\Lambda))$,
the process $\Phi \cdot W$ is a continuous, square integrable, real-valued martingale
on $[0,1]$ and we have the following It\^{o} isometry
\[
    \mathbb{E} \left[ \left\vert \int_{0}^{t} \Phi_{s} \mathrm{d} W_{s} \right\vert ^{2} \right] =
    \mathbb{E} \left[ \int_{0}^{t} \left\Vert \Phi_{s} \right\Vert _{H_{0}^{-1}(\Lambda)}^{2} \mathrm{d}s  \right], \ \forall t \in [0,1].
\]
\end{proposition}

The proof of this proposition is an application of Proposition 4.15 and Proposition 4.17 in \cite{Da-Prato-Book}
together with a continuous extension argument. The following $\mathbb{F}^{W}$-martingale representation theorem can be proved through Theorem 2.5 in \cite{Book-Infinite-SDEs} directly.

\begin{proposition}
    \label{sec-2-proposition-2}
    Suppose that $\Lambda$ is regular. Let $M$ be a real-valued continuous $\mathbb{F}^{W}$-martingale on $[0,1]$ such that $\mathbb{E}\left[ M_{1}^{2} \right] < +\infty$.
Then there exists a unique process $\Phi \in \mathcal{H}_{\mathbb{F}^{W}}^{2}([0,1];H_{0}^{-1}(\Lambda))$ satisfying
\[
M_{t} = M_{0} + \int_{0}^{t} \Phi_{s} \mathrm{d} W_{s}, \quad t \in [0,1].
\]
\end{proposition}

\section{The sine-Gordon measure}
\label{section-3}

For the remainder of the paper, $\Lambda \subset \mathbb{C}$ is always supposed to be a regular bounded simply connected domain and we will not repeat it.
With the above preliminaries and tools in hand, now we can state our main results in this section. 
It should be noted that the constant $C$ appearing the following proofs may change, and we will indicate its dependence on the parameters in our settings
if necessary.

\subsection{Quadratic BSDEs associated with the approximate sine-Gordon measure}
Consider $\eta$ a nonnegative $C^{\infty}$ function being supported in the centered Euclidean unit ball.
We define for $\varepsilon \in (0,1]$ the $\varepsilon$-mollifier $\eta_{\varepsilon}(\cdot) := \varepsilon^{-2} \eta(\varepsilon^{-1} \cdot )$.
For any $\varepsilon, \varepsilon^{\prime} \in (0,1]$, the mollified Green function is defined by 
\begin{equation}
    \label{sec-3-def-mollified-Green}
    G_{\Lambda}^{\varepsilon, \varepsilon^{\prime}}(x,y) := \int_{\mathbb{R}^{2} \times \mathbb{R}^{2}} \eta_{\varepsilon}(x-z_{1}) \eta_{\varepsilon^{\prime}}(y-z_{2}) G_{\Lambda} (z_{1},z_{2}) \mathrm{d}z_{1} \mathrm{d}z_{2}, \quad (x,y) \in \mathbb{R}^{2} \times \mathbb{R}^{2}.
\end{equation}
We simply write $G_{\Lambda}^{\varepsilon}$ when $\varepsilon = \varepsilon^{\prime}$, and $G_{\Lambda}^{\varepsilon}(x)$ when $x = y$.
Introduce the set
\[
\Lambda_{\varepsilon}:=\left\{ x \in \Lambda: \min_{y \in \mathbb{C} \setminus \Lambda} \left\vert y - x \right\vert > 2 \varepsilon \right\}, \quad \varepsilon \in (0,1].
\]
For any $(\varepsilon,x) \in (0,1] \times \Lambda_{\varepsilon}$, the convolution between a distribution $\phi \in \mathcal{D}^{\prime}(\Lambda)$ and $\eta_{\varepsilon}$ is defined by $(\phi \ast \eta_{\varepsilon})(x) := \langle \phi, \eta_{\varepsilon}(x - \cdot) \mathbf{1}_{\Lambda}(\cdot) \rangle$, which does make sense as $\eta_{\varepsilon}(x - \cdot) \mathbf{1}_{\Lambda}(\cdot)$ is actually in $\mathcal{D}_{0}(\Lambda)$.
Then, for each $x \in \Lambda_{\varepsilon}$, the convoluted field of $W_{1}$ is defined by setting $W_{1}^{\varepsilon}(x) := (W_{1} \ast \eta_{\varepsilon})(x)$ since the $H_{0}^{-1}(\Lambda)$-valued process $W_{1}$ can be viewed as an element in $\mathcal{D}^{\prime}(\Lambda)$ almost surely.
Noting that $W_{1}^{\varepsilon}(x) = W_{1} ( \eta_{\varepsilon}(x - \cdot) \mathbf{1}_{\Lambda}(\cdot) )$, it follows from Proposition \ref{sec-2-prop-GFF-Green} that it is a centered Gaussian field indexed by $\mathbb{I} := \bigcup_{\varepsilon \in (0,1]} (\{ \varepsilon \} \times \Lambda_{\varepsilon})$ with covariance function $\mathbb{E} \left[ W_{1}^{\varepsilon}(x) W_{1}^{\varepsilon^{\prime}}(y) \right] = G_{\Lambda}^{\varepsilon, \varepsilon^{\prime}}(x,y)$
whenever $(x,y) \in \Lambda_{\varepsilon} \times \Lambda_{\varepsilon^{\prime}}$.
As $G_{\Lambda}^{\varepsilon, \varepsilon^{\prime}}(x,y)$ is sufficiently regular (that is, both H\"{o}lder continuous in $\varepsilon$ and $x$) to apply Kolmogorov-$\check{\mathrm{C}}$entsov criterion (see, e.g., \cite{Book-Le-Gall}, Theorem 2.9),
there exists a version of $\{ W_{1}^{\varepsilon}(x) \}_{(\varepsilon,x) \in \mathbb{I}}$ which is jointly continuous in $\varepsilon$ and $x$. 
In what follows, we will always be considering this continuous version of the field.

Let $\alpha, \beta$ be two real numbers and $\mu$ be a locally finite Borel (signed) measure on $\Lambda$ of the form $\mu(\mathrm{d}x) = \psi(x) \mathrm{d}x$ for some bounded measurable function $\psi$ whose bound is denoted by $b_{\psi}$. 
For any $\rho \in \mathcal{D}_{0}(\Lambda)$, $F \in \mathcal{B}_{b}(H_{0}^{-1}(\Lambda))$, the bound of $\rho$ (resp. F) is denoted by $b_{\rho}$ (resp. $b_{F}$).
Given $\rho \in \mathcal{D}_{0}(\Lambda)$ arbitrarily, we normalize (\ref{introduction-SG-measure-formal}) by introducing the approximate sine-Gordon measure with ultraviolet cutoffs:
\[
    \mu_{\mathrm{SG}}^{\rho, \varepsilon}( \mathrm{d}\phi) := \Xi_{\rho, \varepsilon}^{-1} \exp\left\{ \alpha \int_{\Lambda_{\varepsilon}} e^{ \frac{\beta^{2}}{2} G_{\Lambda}^{\varepsilon}(x) } \cos(\beta \phi^{\varepsilon}(x)) \rho(x) \mu(\mathrm{d}x) \right\} \mathrm{Law}_{W_{1}}(\mathrm{d} \phi),
\]
with the notation $\phi^{\varepsilon} := \phi \ast \eta_{\varepsilon}$ for $\phi \in H_{0}^{-1}(\Lambda)$, where the normalized partition function
\[
    \Xi_{\rho, \varepsilon} := \int_{H_{0}^{-1}(\Lambda)} \exp\left\{ \alpha \int_{\Lambda_{\varepsilon}} e^{ \frac{\beta^{2}}{2} G_{\Lambda}^{\varepsilon}(x) } \cos(\beta \phi^{\varepsilon}(x)) \rho(x) \mu(\mathrm{d}x) \right\} \mathrm{Law}_{W_{1}}(\mathrm{d} \phi).
\]
Below we will see from (\ref{sec-3-terminal-absolute}) and (\ref{sec-3-terminal-L1}) that $\mu_{\mathrm{SG}}^{\rho, \varepsilon}( e^{p |F|} ) < +\infty$ for all $p > 0$ and all $F \in \mathcal{B}_{b}(H_{0}^{-1}(\Lambda))$.
Thus the following nonstandard log-Laplace transform
\begin{align*}
    \delta & \mapsto \frac{1}{\alpha} \log \int_{H_{0}^{-1}(\Lambda)} \exp\left\{ \alpha \left[ \delta F(\phi) + \left( e^{ \frac{\beta^{2}}{2} G_{\Lambda}^{\varepsilon} } \cos(\beta \phi^{\varepsilon}), \rho \right)_{L^{2}(\Lambda_{\varepsilon};\mu)} \right] \right\} \mathrm{Law}_{W_{1}}(\mathrm{d}\phi) \\
   & = \frac{1}{\alpha} \log \left[ \Xi_{\rho, \varepsilon} \mu_{\mathrm{SG}}^{\rho, \varepsilon}( e^{\alpha \delta F} ) \right]
\end{align*}
is analytic (at least) in a neighborhood of zero. Using the change of variables, $\mu_{\mathrm{SG}}^{\rho, \varepsilon}$ can be described by
\begin{equation}
    \label{sec-3-renormalized-SG-1}
    \begin{array}
    [c]{rl}%
    \mu_{\mathrm{SG}}^{\rho, \varepsilon}( F )
    = & \displaystyle \left. \frac{\mathrm{d}}{\mathrm{d}\delta} \right\vert _{\delta=0} \left( \frac{1}{\alpha} \log \left[ \Xi_{\rho, \varepsilon} \mu_{\mathrm{SG}}^{\rho, \varepsilon}( e^{\alpha \delta F} ) \right] \right) \\
    = & \displaystyle \left.\frac{\mathrm{d}}{\mathrm{d}\delta} \right\vert _{\delta=0} \left( \alpha^{-1} \log \mathbb{E} \left[ e^{ \alpha \left( \delta F(W_{1}) + \left( \left[\left[\cos(\beta W_{1}^{\varepsilon})\right]\right], \rho \right) _{L^{2}(\Lambda_{\varepsilon};\mu)} \right) } \right] \right) \\
    = & \displaystyle \lim\limits_{\delta \rightarrow 0} \frac{\alpha^{-1}}{\delta} \log \left( \frac{\mathbb{E} \left[  e^{ \alpha \left( \delta F(W_{1}) + \left( \left[\left[\cos(\beta W_{1}^{\varepsilon})\right]\right], \rho \right) _{L^{2}(\Lambda_{\varepsilon};\mu)} \right) } \right]} {\mathbb{E} \left[  e^{ \alpha \left( \left[\left[\cos(\beta W_{1}^{\varepsilon})\right]\right], \rho \right) _{L^{2}(\Lambda_{\varepsilon};\mu)} } \right]} \right),
    \end{array}
\end{equation}
where the Wick-ordered cosine of $W_{1}^{\varepsilon}$ is defined by setting
\[
\left[\left[\cos( \beta W_{1}^{\varepsilon}(x) )\right]\right] := \exp \left\{ \frac{\beta^{2}}{2} G_{\Lambda}^{\varepsilon}(x) \right\} \cos( \beta W_{1}^{\varepsilon}(x) ), \quad \forall x \in \Lambda_{\varepsilon}.
\]
By the way, to be used later, the Wick-ordered sine of $W_{1}^{\varepsilon}$ is defined by setting
\[
    \left[\left[\sin( \beta W_{1}^{\varepsilon}(x) )\right]\right] := \exp \left\{ \frac{\beta^{2}}{2} G_{\Lambda}^{\varepsilon}(x) \right\} \sin( \beta W_{1}^{\varepsilon}(x) ), \quad \forall x \in \Lambda_{\varepsilon}.
\]
(\ref{sec-3-renormalized-SG-1}) does make sense because we can show that $\left( \left[\left[\cos(\beta W_{1}^{\varepsilon})\right]\right], \rho \right) _{L^{2}(\Lambda_{\varepsilon};\mu)} $ is in $ L^{\infty}(\mathcal{F}_{1}^{W};\mathbb{R})$.
Indeed, note that
\begin{equation}
    \label{sec-3-terminal-absolute}
    \left\vert \int_{\Lambda_{\varepsilon}} \exp \left\{ \frac{\beta^{2}}{2} G_{\Lambda}^{\varepsilon}(x) \right\} \cos(\beta W_{1}^{\varepsilon}(x)) \rho(x) \mu(\mathrm{d}x) \right\vert \leqslant b_{\rho} b_{\psi} \int_{\Lambda_{\varepsilon}} \exp \left\{ \frac{\beta^{2}}{2} G_{\Lambda}^{\varepsilon}(x) \right\} \mathrm{d}x.
\end{equation}
As $G_{\Lambda}^{\varepsilon}(x) \geqslant 0$ for $x \in \mathbb{R}^{2}$, we have, for any $k \in \mathbb{N}$,
\begin{align*}
    0 \leqslant \int_{\Lambda_{\varepsilon}} \left[G_{\Lambda}^{\varepsilon}(x)\right]^{k} \mathrm{d}x 
        = & \int_{\Lambda_{\varepsilon} \times \Lambda_{\varepsilon}} \left[G_{\Lambda}^{\varepsilon}(x,y)\right]^{k} \mathbf{1}_{\{x=y\}}(x,y) \mathrm{d}x \mathrm{d}y \\
        \leqslant & \int_{\mathbb{R}^{2} \times \mathbb{R}^{2}} \left[ G_{\Lambda}^{\varepsilon}(x,y) \right]^{k} \mathrm{d}x \mathrm{d}y\\
        \leqslant & \int_{\mathbb{R}^{2} \times \mathbb{R}^{2}} \left[ G_{\Lambda}(x,y) \right]^{k} \mathrm{d}x \mathrm{d}y \\
        = & \int_{\Lambda \times \Lambda} \left[ G_{\Lambda}(x,y) \right]^{k} \mathrm{d}x \mathrm{d}y,
\end{align*}
where the last inequality results from Young's convolution inequality (\cite{Book-Fourier-PDE}, Lemma 1.4), and the last line is due to the extension that $G_{\Lambda}(x,y) = 0$ whenever $(x,y) \notin \Lambda \times \Lambda$.
For $\beta^{2} \in [0,2)$, using the monotonicity convergence and Proposition \ref{sec-2-proposition-Green-func} yields that
\begin{equation}
    \label{sec-3-terminal-L1}
    \begin{array}[c]{rl}
        \int_{\Lambda_{\varepsilon}} \exp \left\{ \frac{\beta^{2}}{2} G_{\Lambda}^{\varepsilon}(x) \right\} \mathrm{d}x
        = & \sum\limits _{k=1}^{\infty} \frac{\beta^{2k}}{2^{k}k!} \int_{\Lambda_{\varepsilon}} \left[G_{\Lambda}^{\varepsilon}(x)\right]^{k} \mathrm{d}x \\
        \leqslant & \sum\limits _{k=1}^{\infty} \frac{\beta^{2k}}{2^{k}k!} \int_{\Lambda \times \Lambda} \left[ G_{\Lambda}(x,y) \right]^{k} \mathrm{d}x \mathrm{d}y \\
        = & \int_{\Lambda \times \Lambda} \exp \left\{ \frac{\beta^{2}}{2} G_{\Lambda}(x,y) \right\} \mathrm{d}x \mathrm{d}y \\
        \leqslant & e^{ \frac{\beta^{2}}{2} \sup\limits _{(x,y) \in \Lambda \times \Lambda}  g_{\Lambda}(x,y) } \int_{\Lambda \times \Lambda} \left\vert x - y \right\vert ^{-\frac{\beta^{2}}{2}} \mathbf{1}_{\{ |x - y| \leqslant 1 \}} \mathrm{d}x \mathrm{d}y \\
        \leqslant & C_{\beta} < +\infty,
    \end{array}
\end{equation}
where $C_{\beta}$ is a positive constant depending on $\beta$ and $\Lambda$.
Combing (\ref{sec-3-terminal-absolute}) with (\ref{sec-3-terminal-L1}) we actually obtain that $\sup_{\varepsilon \in (0,1]} \left\Vert \left( \left[\left[\cos(\beta W_{1}^{\varepsilon})\right]\right], \rho \right) _{L^{2}(\Lambda_{\varepsilon};\mu)} \right\Vert _{L^{\infty}(\mathcal{F}_{1}^{W};\mathbb{R})} \leqslant b_{\rho} b_{\psi} C_{\beta}$.

For any $\rho \in \mathcal{D}_{0}(\Lambda)$, $F \in \mathcal{B}_{b}(H_{0}^{-1}(\Lambda))$, we introduce the following quadratic BSDE associated with the approximate sine-Gordon measure, that is,
\begin{equation}
    \label{sec-3-quadratic-BSDE}
    \left\{
    \begin{array}[c]{rl}
       \mathrm{d} Y_{t}^{\varepsilon, \delta F}(\rho) = & \displaystyle - \frac{\alpha}{2} \left\Vert Z_{t}^{\varepsilon, \delta F}(\rho)  \right\Vert _{H_{0}^{-1}(\Lambda)}^{2} \mathrm{d}t +  Z_{t}^{\varepsilon, \delta F}(\rho) \mathrm{d} W_{t}, \quad t \in [0,1),\\
       Y_{1}^{\varepsilon, \delta F}(\rho) = & \delta F(W_{1}) + \left( \left[\left[\cos(\beta W_{1}^{\varepsilon})\right]\right], \rho \right) _{L^{2}(\Lambda_{\varepsilon};\mu)},
    \end{array}
    \right.
\end{equation}
where $(\varepsilon, \delta) \in (0,1] \times [-1,0) \cup (0,1]$. By a solution to (\ref{sec-3-quadratic-BSDE}) we mean a pair $( Y^{\varepsilon, \delta F}(\rho),Z^{\varepsilon, \delta F}(\rho) ) = \{ ( Y_{t}^{\varepsilon, \delta F}(\rho),Z_{t}^{\varepsilon, \delta F}(\rho) ) \}_{t \in [0,1]}$ of predictable processes with values
in $\mathbb{R} \times H_{0}^{-1}(\Lambda)$ satisfying (\ref{sec-3-quadratic-BSDE}) $\mathbb{P}$-a.s. The definition of solutions to other BSDEs in the rest of this paper is similar so we will not repeat it below.

\begin{lemma}
    \label{sec-3-lemma-1}
    Let $\rho \in \mathcal{D}_{0}(\Lambda)$, $F \in \mathcal{B}_{b}(H_{0}^{-1}(\Lambda))$ be given arbitrarily. If $\beta^{2} \in [0,2)$, then (\ref{sec-3-quadratic-BSDE})
    admits a unique solution $\left(Y^{\varepsilon, \delta F}(\rho),Z^{\varepsilon, \delta F}(\rho) \right)$ belonging to $\mathcal{S}_{\mathbb{F}^{W}}^{\infty}([0,1];\mathbb{R}) \times \mathcal{H}_{\mathbb{F}^{W}}^{2}([0,1];H_{0}^{-1}(\Lambda))$ 
    such that $|| Z^{\varepsilon, \delta F}(\rho) \cdot W ||_{\mathrm{BMO}} \leqslant \sqrt{2} |\alpha|^{-1} e^{\frac{|\alpha|}{2} (b_{F} + b_{\rho} b_{\psi} C_{\beta})}$.
\end{lemma}

\begin{proof}
    As $\left\vert  \delta F(W_{1}) + \left( \left[\left[\cos(\beta W_{1}^{\varepsilon})\right]\right], \rho \right) _{L^{2}(\Lambda_{\varepsilon};\mu)} \right\vert \leqslant b_{F} + b_{\rho} b_{\psi} C_{\beta}$, we put
    \[
    \mathcal{Y}_{t} := \mathbb{E} \left[ \exp \left\{ \alpha \left( \delta F(W_{1}) + \left( \left[\left[\cos(\beta W_{1}^{\varepsilon})\right]\right], \rho \right) _{L^{2}(\Lambda_{\varepsilon};\mu)} \right) \right\} \mid \mathcal{F}_{t}^{W} \right], \quad t \in [0,1].
    \]
    Then $\mathcal{Y}$ is a real-valued martingale with $\mathbb{E} [\mathcal{Y}_{1}^{2}] < +\infty$, and it admits a continuous modification due to the continuity of $\mathbb{F}^{W}$ (In what follows, $\mathcal{Y}$ will always be identified with this continuous modification).
    According to Proposition \ref{sec-2-proposition-2}, there exists a unique process $\mathcal{Z} \in \mathcal{H}_{\mathbb{F}^{W}}^{2}([0,1];H_{0}^{-1}(\Lambda))$ such that
    \[
    \mathcal{Y}_{t} = \exp \left\{ \alpha \left( \delta F(W_{1}) + \left( \left[\left[\cos(\beta W_{1}^{\varepsilon})\right]\right], \rho \right) _{L^{2}(\Lambda_{\varepsilon};\mu)} \right) \right\} - \int_{t}^{1} \mathcal{Z}_{s} \mathrm{d} W_{s}, \quad t \in [0,1].
    \]
    Employing It\^{o}'s formula (\cite{Da-Prato-Book}, Theorem 4.17), one can check that $\left( \alpha^{-1} \log (\mathcal{Y}), \mathcal{Z}/ (\alpha \mathcal{Y}) \right)$ satisfies (\ref{sec-3-quadratic-BSDE}),
    which lives in $\mathcal{S}_{\mathbb{F}^{W}}^{\infty}([0,1];\mathbb{R}) \times \mathcal{H}_{\mathbb{F}^{W}}^{2}([0,1];H_{0}^{-1}(\Lambda))$ since $\mathcal{Y}$ is uniformly bounded with a strictly positive lower bound.
    Next, for any solution $(Y, Z) \in \mathcal{S}_{\mathbb{F}^{W}}^{\infty}([0,1];\mathbb{R}) \times \mathcal{H}_{\mathbb{F}^{W}}^{2}([0,1];H_{0}^{-1}(\Lambda))$, we will prove that $Z \cdot W$ is a $\mathrm{BMO}$ martingale by
    considering the function $u(x) = |\alpha|^{-2} (e^{ |\alpha| |x|} - |\alpha| |x| - 1), x \in \mathbb{R}$. Then it can be checked that $u \in \mathcal{C}^{2}(\mathbb{R})$ and for $x \in \mathbb{R}$,
    \[
    \begin{array}[c]{rcl}
        u^{\prime}(x) = |\alpha|^{-1} (e^{|\alpha| |x|} - 1)\mathrm{sgn}(x), & u^{\prime \prime}(x) = e^{|\alpha| |x|}, & u^{\prime \prime}(x) - |\alpha| |u^{\prime}(x)| = 1.
    \end{array}
    \]
    Employing It\^{o}'s formula to compute $u(Y_{t})$, we have
    \begin{align*}
        & u(Y_{t}) +  \frac{1}{2} \int_{t}^{1} \left\Vert Z_{s} \right\Vert _{H_{0}^{-1}(\Lambda)}^{2} \mathrm{d}s \\
        \leqslant & \ u\left( \alpha \delta F(W_{1}) + \alpha \left( \left[\left[\cos(\beta W_{1}^{\varepsilon})\right]\right], \rho \right) _{L^{2}(\Lambda_{\varepsilon};\mu)} \right) - \int_{t}^{1} u^{\prime}(Y_{s}) Z_{s} \mathrm{d} W_{s}.
    \end{align*}
    As $(Y, Z) \in \mathcal{S}_{\mathbb{F}^{W}}^{\infty}([0,1];\mathbb{R}) \times \mathcal{H}_{\mathbb{F}^{W}}^{2}([0,1];H_{0}^{-1}(\Lambda))$, the stochastic integral in the above BSDE is a true martingale with the mean zero (see Proposition 4.13 in \cite{Da-Prato-Book}).
    Taking $\mathbb{E}[ \cdot \mid \mathcal{F}_{t}^{W} ]$ on both sides yields that
    \begin{equation}
        \label{sec-3-QBSDE-Z-BMO-uniform-bounded}
        \mathbb{E}\left[ \int_{t}^{1} \left\Vert Z_{s} \right\Vert _{H_{0}^{-1}(\Lambda)}^{2} \mathrm{d}s \mid \mathcal{F}_{t}^{W} \right] \leqslant 2 |\alpha|^{-2} e^{|\alpha| (b_{F} + b_{\rho} b_{\psi} C_{\beta})} < +\infty, \quad \forall t \in [0,1],
    \end{equation}
    $\mathbb{P}$-almost surely, which implies that $Z \cdot W$ is indeed a $\mathrm{BMO}$ martingale.

    To prove the uniqueness, for $k = 1,2$, let $(Y^{(k)},Z^{(k)}) \in \mathcal{S}_{\mathbb{F}^{W}}^{\infty}([0,1];\mathbb{R}) \times \mathcal{H}_{\mathbb{F}^{W}}^{2}([0,1];H_{0}^{-1}(\Lambda))$ be two pairs of solution to (\ref{sec-3-quadratic-BSDE}).
    Subtracting $(Y^{(2)},Z^{(2)})$ from $(Y^{(1)},Z^{(1)})$, we deduce from Theorem 4.12 and Corollary 4.14 in \cite{Da-Prato-Book} that
    \begin{equation}
        \label{sec-3-BSDE-unique}
         \begin{array}[c]{rl}
        Y_{t}^{(1)} - Y_{t}^{(2)} = & \frac{1}{2} \int_{t}^{1} \left( Z_{s}^{(1)} - Z_{s}^{(2)}, Z_{s}^{(1)} + Z_{s}^{(2)} \right) _{H_{0}^{-1}(\Lambda)} \mathrm{d}s - \int_{t}^{1} (Z_{s}^{(1)} - Z_{s}^{(2)}) \mathrm{d} W_{s} \\
        = & \int_{t}^{1} \mathrm{d} \ll (Z^{(1)} - Z^{(2)}) \cdot W, \frac{1}{2}(Z^{(1)} + Z^{(2)}) \cdot W \gg _{s} \\
        & - \int_{t}^{1} \mathrm{d} \left( (Z^{(1)} - Z^{(2)}) \cdot W \right)_{s} \\
        = & - \int_{t}^{1} \mathrm{d} \left( M - \ll M,N \gg \right) _{s},
    \end{array}
    \end{equation}
    where $M := (Z^{(1)} - Z^{(2)}) \cdot W$, $N := \frac{1}{2}(Z^{(1)} + Z^{(2)}) \cdot W$, and $\ll M,N \gg$ denotes the predictable quadratic covariance of real-valued martingales $M$ and $N$.
    Since $Z^{(k)} \cdot W$ is a $\mathrm{BMO}$ martingale for $k=1,2$, it can be verified that both $M$ and $N$ are $\mathrm{BMO}$ martingales. By Lemma A.4 in \cite{Hu-Tang-2016}, there is a new probability measure $\mathbb{Q}$ defined by $\mathrm{d} \mathbb{Q} := \mathcal{E}(N_{1}) \mathrm{d} \mathbb{P}$
    and two constants $c_{1}>0$ and $c_{2}>0$ depending only on $\alpha$, $\beta$, $b_{F}$, $b_{\rho}$, and $b_{\psi}$, such that
    \[
    c_{1} \left\Vert M \right\Vert _{\mathrm{BMO}(\mathbb{P})} \leqslant \left\Vert \tilde{M} \right\Vert _{\mathrm{BMO}(\mathbb{Q})} \leqslant c_{2} \left\Vert M \right\Vert _{\mathrm{BMO}(\mathbb{P})},
    \]
    where $\tilde{M} := M - \ll M,N \gg $. Thus we can take $\mathbb{E}_{\mathbb{Q}} [ \cdot \mid \mathcal{F}_{t}^{W}]$
    on both sides of (\ref{sec-3-BSDE-unique}) and obtain $Y_{t}^{(1)} = Y_{t}^{(2)}$ for each $t \in [0,1]$, $\mathbb{Q}$-a.s. (of course, $\mathbb{P}$-a.s.). Then $\left\Vert Z_{t}^{(1)}(\omega) - Z_{t}^{(2)}(\omega) \right\Vert _{H_{0}^{-1}(\Lambda)} = 0$ for a.e $(t,\omega) \in [0,1] \times \Omega$ 
    thanks to (\ref{sec-3-BSDE-unique}) and Proposition \ref{sec-2-proposition-3}, since any continuous local martingale with the initial value zero and finite variation paths equals zero up to an evanescent set.
\end{proof}

By Lemma \ref{sec-3-lemma-1}, the unique solution $Y^{\varepsilon, \delta F}(\rho)$ of (\ref{sec-3-quadratic-BSDE}) has the form
\[
    Y_{t}^{\varepsilon, \delta F}(\rho) = \frac{1}{\alpha} \log \left( \mathbb{E} \left[ \exp \left\{ \alpha \left( \delta F(W_{1}) + \left( \left[\left[\cos(\beta W_{1}^{\varepsilon})\right]\right], \rho \right) _{L^{2}(\Lambda_{\varepsilon};\mu)} \right) \right\} \mid \mathcal{F}_{t}^{W} \right] \right)
\]
for $t \in [0,1]$, whence (\ref{sec-3-renormalized-SG-1}) can be further expressed by
\[
    \int_{H_{0}^{-1}(\Lambda)} F(\phi) \mu_{\mathrm{SG}}^{\rho, \varepsilon}( \mathrm{d}\phi) = \lim\limits_{\delta \rightarrow 0} \delta^{-1} \left[ Y_{0}^{\varepsilon, \delta F}(\rho) - Y_{0}^{\varepsilon,0}(\rho) \right].
\]
Here, $\left(Y^{\varepsilon,0}(\rho), Z^{\varepsilon,0}(\rho)\right)$ denotes the unique solution of (\ref{sec-3-quadratic-BSDE}) corresponding to $F = 0$.
To compute the limit on the right-hand side of the above equation, we resort to a variational approach originating from the control theory. To be more precise, we introduce the following variational equation
\begin{equation}
    \label{sec-3-variational-eq}
    \left\{
    \begin{array}[c]{rl}
       \mathrm{d} \widehat{Y}_{t}^{\varepsilon, F}(\rho) = & - \alpha \left( Z_{t}^{\varepsilon, 0}(\rho) , \widehat{Z}_{t}^{\varepsilon, F}(\rho) \right) _{H_{0}^{-1}(\Lambda)} \mathrm{d}t +  \widehat{Z}_{t}^{\varepsilon, F}(\rho) \mathrm{d} W_{t}, \quad t \in [0,1),\\
       \widehat{Y}_{1}^{\varepsilon, F}(\rho) = & F(W_{1}),
    \end{array}
    \right.
\end{equation}
and conjecture that $Y_{0}^{\varepsilon, \delta F}(\rho)$ admits a first-order Taylor expansion with respect to $\delta$ around $0$, that is,
\begin{equation}
    \label{sec-3-Taylor-expansion}
    Y_{0}^{\varepsilon, \delta F}(\rho) = Y_{0}^{\varepsilon,0}(\rho) + \delta \widehat{Y}_{0}^{\varepsilon, F}(\rho) + o(\delta), \quad \text{as} \ \delta \rightarrow 0.
\end{equation}
If BSDE (\ref{sec-3-variational-eq}) admits a unique solution, then we will deduce from (\ref{sec-3-Taylor-expansion}) that
\begin{equation}
    \label{sec-3-SG-variational-expression}
    \mu_{\mathrm{SG}}^{\rho, \varepsilon}( F ) 
    = \lim\limits_{\delta \rightarrow 0} \delta^{-1} \left[ \left( Y_{0}^{\varepsilon,0}(\rho) + \delta \widehat{Y}_{0}^{\varepsilon, F}(\rho) + o(\delta) \right) - Y_{0}^{\varepsilon,0}(\rho) \right] = \widehat{Y}_{0}^{\varepsilon,F}(\rho).
\end{equation}
To prove the well-posedness of (\ref{sec-3-variational-eq}) and Taylor expansion (\ref{sec-3-Taylor-expansion}), we need the following error estimate.

\begin{proposition}
    \label{sec-3-prop-1}
    Let $\rho \in \mathcal{D}_{0}(\Lambda)$, $F \in \mathcal{B}_{b}(H_{0}^{-1}(\Lambda))$ be given arbitrarily. If $\beta^{2} \in [0,2)$, then
    \begin{equation}
        \label{sec-3-variation-error-estimate}
        \left\Vert Y^{\varepsilon, \delta F}(\rho) - Y^{\varepsilon,0}(\rho) \right\Vert _{\mathcal{S}_{\mathbb{F}^{W}}^{p}([0,1];\mathbb{R})} + \left\Vert Z^{\varepsilon, \delta F}(\rho) - Z^{\varepsilon,0}(\rho) \right\Vert _{\mathcal{H}_{\mathbb{F}^{W}}^{p}([0,1];H_{0}^{-1}(\Lambda))} \leqslant C |\delta|
    \end{equation}
    for any $(\varepsilon, \delta) \in (0,1] \times [-1,0) \cup (0,1]$ and $p \in [1, +\infty)$, where the constant $C>0$ depends only on $p$, $\alpha$, $\beta$, $b_{\rho}$, $b_{F}$, and $b_{\psi}$.
\end{proposition}

\begin{proof}
For every $(\varepsilon, \delta) \in (0,1] \times [-1,1]$, consider a $\mathfrak{B}([0,1]) \otimes \mathcal{F}_{1}^{W} \otimes \mathfrak{B}(H_{0}^{-1}(\Lambda))$ jointly measurable function $f^{(\varepsilon, \delta)}: [0,1] \times \Omega \times H_{0}^{-1}(\Lambda) \mapsto \mathbb{R}$ which is defined by setting,
for any $(t,z) \in [0,1] \times H_{0}^{-1}(\Lambda)$,
\[
    f^{(\varepsilon, \delta)}(t,z) :=
    \left\{
\begin{array}{ll}
    \alpha \left( \frac{Z_{t}^{\varepsilon, \delta F}(\rho) + Z_{t}^{\varepsilon, 0}(\rho)}{2} , z \right) _{H_{0}^{-1}(\Lambda)} & , \quad \delta \in [-1,0) \cup (0,1]; \\
    \alpha \left( Z_{t}^{\varepsilon, 0}(\rho), z \right) _{H_{0}^{-1}(\Lambda)} & , \quad \delta = 0
\end{array}
    \right.
\]
(here, $f^{(\varepsilon, \delta)}$ relies on $\omega$ through $Z^{\varepsilon, \delta F}(\rho)$ and $Z^{\varepsilon, 0}(\rho)$ but we do not write the dependence explicitly). 
Then, for any $t \in [0,1]$ and $(z_{1}, z_{2}) \in H_{0}^{-1}(\Lambda) \times H_{0}^{-1}(\Lambda)$, through the Cauchy-Schwartz inequality we have
\[
    \left\vert f^{(\varepsilon, \delta)}(t,z_{1}) - f^{(\varepsilon, \delta)}(t,z_{2}) \right\vert \leqslant \Phi_{t}^{(\varepsilon, \delta)} \left\Vert z_{1} - z_{2} \right\Vert _{H_{0}^{-1}(\Lambda)},
\]
with the real-valued process
\[
\Phi_{t}^{(\varepsilon, \delta)} :=
\left\{
\begin{array}{ll}
    \displaystyle \frac{\alpha}{2} \left\Vert Z_{t}^{\varepsilon, \delta F}(\rho) + Z_{t}^{\varepsilon, 0}(\rho)  \right\Vert _{H_{0}^{-1}(\Lambda)} & , \quad \delta \in [-1,0) \cup (0,1]; \\
    \displaystyle \alpha \left\Vert Z_{t}^{\varepsilon, 0}(\rho) \right\Vert _{H_{0}^{-1}(\Lambda)} & , \quad \delta = 0.
\end{array}
\right.
\]
Recall that $\{ B_{t}^{(1)} = B_{t} \circ \Pi^{1}: t \in [0,1] \}$ is a one-dimensional, $\mathbb{F}^{W}$-adapted, standard Wiener process under $\mathbb{P}$.
The ordinary stochastic integral $\Phi^{(\varepsilon, \delta)} \cdot B^{(1)}$ is in $\mathrm{BMO}$ due to Lemma \ref{sec-3-lemma-1} and the definition of $\Phi^{(\varepsilon, \delta)}$.
In indeed, it follows from the definition of BMO martingales, the triangle inequality and estimate (\ref{sec-3-QBSDE-Z-BMO-uniform-bounded}) that 
\begin{equation}
    \label{sec-3-BMO-Lip-uniform-bound}
    \begin{array}[c]{rl}
    \left\Vert \Phi^{(\varepsilon, \delta)} \cdot B^{(1)} \right\Vert _{\mathrm{BMO}} \leqslant &  |\alpha| \cdot \left( \left\Vert \frac{Z^{\varepsilon, \delta F}(\rho) + Z^{\varepsilon, 0}(\rho)}{2}\cdot W \right\Vert _{\mathrm{BMO}} \vee \left\Vert Z^{\varepsilon, 0}(\rho) \cdot W \right\Vert _{\mathrm{BMO}} \right) \\
    \leqslant & \sqrt{2} e^{\frac{|\alpha|}{2} (b_{F} + b_{\rho} b_{\psi} C_{\beta})}, \quad \forall (\varepsilon, \delta) \in (0,1] \times [-1,1].
    \end{array}
\end{equation}
Recall the function $\kappa$ defined by (\ref{sec-2-func-Psi}). Let $\overline{p}$ be uniquely determined by $\kappa(\overline{p}) = \sqrt{2} e^{\frac{|\alpha|}{2} (b_{F} + b_{\rho} b_{\psi} C_{\beta})}$ and $\overline{p}^{\ast}$ be the conjugate exponent of $\overline{p}$.
Then, thanks to (\ref{sec-2-reverse-Holder-coeff}), for all $(\varepsilon, \delta) \in (0,1] \times [-1,1]$, we have
\begin{equation}
    \label{sec-3-constant-K-compare}
    K \left(q, \left\Vert \Phi^{(\varepsilon, \delta)} \cdot B^{(1)} \right\Vert _{\mathrm{BMO}} \right) \leqslant K \left(q, \sqrt{2} e^{\frac{|\alpha|}{2} (b_{F} + b_{\rho} b_{\psi} C_{\beta})} \right) < +\infty, \quad \forall \ q \in (1, \overline{p}).
\end{equation}
When $(\varepsilon, \delta) \in (0,1] \times [-1,0) \cup (0,1]$, consider the BSDE
\[
    \begin{array}[c]{rl}
        Y_{t}^{\varepsilon, \delta F}(\rho) - Y_{t}^{\varepsilon,0}(\rho) = & \displaystyle \delta F(W_{1}) + \int_{t}^{1} f^{(\varepsilon, \delta)}\left( s,Z_{s}^{\varepsilon, \delta F}(\rho) - Z_{s}^{\varepsilon,0}(\rho) \right) \mathrm{d}s \\
        & \displaystyle - \int_{t}^{1} \left( Z_{s}^{\varepsilon, \delta F}(\rho) - Z_{s}^{\varepsilon,0}(\rho) \right) \mathrm{d} W_{s}.
    \end{array}
\]
One can check that the terminal condition $\delta F$ and the generator $f^{(\varepsilon, \delta)}$ verify the Assumptions A1--A4 in \cite{Briand-Confortola-2008}.
Since $Y^{\varepsilon, \delta F}(\rho)$ and $Y^{\varepsilon,0}(\rho) $ are in $\mathcal{S}_{\mathbb{F}^{W}}^{\infty}([0,1];\mathbb{R})$, we can use Corollary 9 in \cite{Briand-Confortola-2008}
to deduce that, for any $p > \overline{p}^{\ast}$ with $\lfloor p \rfloor$ denoting the maximal integer no more than $p$, there exists a constant $C_{p}>0$ depending only on $p$ such that
\[
    \begin{array}[c]{rl}
        & \left\Vert Y^{\varepsilon, \delta F}(\rho) - Y^{\varepsilon,0}(\rho) \right\Vert _{\mathcal{S}_{\mathbb{F}^{W}}^{p}([0,1];\mathbb{R})} + \left\Vert Z^{\varepsilon, \delta F}(\rho) - Z^{\varepsilon,0}(\rho) \right\Vert _{\mathcal{H}_{\mathbb{F}^{W}}^{p}([0,1];H_{0}^{-1}(\Lambda))} \\
        \leqslant & C_{p} K \left( \frac{p}{p-1}, \left\Vert M^{(\varepsilon, \delta)} \right\Vert _{\mathrm{BMO}} \right) ^{\frac{p-1}{p}} \left( 1 + \left\Vert \Phi^{(\varepsilon, \delta)} \right\Vert _{\mathcal{H}_{\mathbb{F}^{W}}^{2p}([0,1];\mathbb{R})} \right) \left\Vert \delta F(W_{1}) \right\Vert _{L^{3p}(\mathcal{F}_{1}^{W}; \mathbb{R})} \\
        \leqslant & C_{p} b_{F} K \left( \frac{p}{p-1}, \left\Vert M^{(\varepsilon, \delta)} \right\Vert _{\mathrm{BMO}} \right) ^{\frac{p-1}{p}} \left\{ 1 + \left[ (\lfloor p \rfloor + 1)! \left\Vert M^{(\varepsilon, \delta)}\right\Vert _{\mathrm{BMO}}^{2(\lfloor p \rfloor + 1)} \right] ^{\frac{1}{2(\lfloor p \rfloor + 1)}} \right\} |\delta|\\
        \leqslant & C_{p} b_{F} K \left( \frac{p}{p-1}, \sqrt{2} e^{\frac{|\alpha|}{2} (b_{F} + b_{\rho} b_{\psi} C_{\beta})} \right) ^{\frac{p-1}{p}} \left[ 1 + (\lfloor p \rfloor + 1)! \sqrt{2} e^{\frac{|\alpha|}{2} (b_{F} + b_{\rho} b_{\psi} C_{\beta})} \right] |\delta|,
    \end{array}
\]
with the notation $M^{(\varepsilon, \delta)} := \Phi^{(\varepsilon, \delta)} \cdot B^{(1)}$,
where the second inequality is due to the energy inequality for $\mathrm{BMO}$ martingales (\cite{BMO-Book}, p.26) and H\"{o}lder's inequality,
and the last inequality follows from (\ref{sec-3-BMO-Lip-uniform-bound}) and (\ref{sec-3-constant-K-compare}). Thus we have proved (\ref{sec-3-variation-error-estimate}) for $p \in (\overline{p}^{\ast}, +\infty)$.
The estimate for the case $p \in [1, \overline{p}^{\ast}]$ follows from H\"{o}lder's inequality immediately.
\end{proof}

\begin{remark}
In the proof of Proposition \ref{sec-3-prop-1}, $\Phi^{(\varepsilon, \delta)} \cdot B^{(1)}$ is well defined since $\Phi^{(\varepsilon, \delta)}$ is $\mathbb{F}^{W}$-predictable.
In addition, the choice of such a Wiener process is not unique so that one can replace $B^{(1)}$ with any other $\mathbb{F}^{W}$-adapted standard Wiener process under $\mathbb{P}$, say, $B^{(k)}$ for some $k \in \mathbb{N}_{+}$.
\end{remark}

Now we can prove the well-posedness of (\ref{sec-3-variational-eq}) and Taylor expansion (\ref{sec-3-Taylor-expansion}) by the following lemma.

\begin{lemma}
    \label{sec-3-lemma-2}
    Let $\rho \in \mathcal{D}_{0}(\Lambda)$, $F \in \mathcal{B}_{b}(H_{0}^{-1}(\Lambda))$ be given arbitrarily. If $\beta^{2} \in [0,2)$, then, for any $\varepsilon \in (0,1]$, (\ref{sec-3-variational-eq}) admits a unique solution 
    $$
    \left( \widehat{Y}^{\varepsilon, F}(\rho), \widehat{Z}^{\varepsilon, F}(\rho) \right) \in \mathcal{S}_{\mathbb{F}^{W}}^{\infty}([0,1];\mathbb{R}) \times \bigcap_{p>1} \mathcal{H}_{\mathbb{F}^{W}}^{p}([0,1];H_{0}^{-1}(\Lambda))
    $$
    such that (\ref{sec-3-Taylor-expansion}) hold.
\end{lemma}

\begin{proof}
    Adopting the notation used in the proof of Proposition \ref{sec-3-prop-1}, BSDE (\ref{sec-3-variational-eq}) can be rewritten as
    \[
        \widehat{Y}_{t}^{\varepsilon, F}(\rho) = F(W_{1}) + \int_{t}^{1} f^{(\varepsilon, 0)}\left( s, \widehat{Z}_{s}^{\varepsilon, F}(\rho) \right) \mathrm{d}s - \int_{t}^{1} \widehat{Z}_{s}^{\varepsilon, F}(\rho) \mathrm{d} W_{s}.
    \]
    As $F$ is bounded, there is a couple $\left( \widehat{Y}^{\varepsilon, F}(\rho), \widehat{Z}^{\varepsilon, F}(\rho) \right) \in \mathcal{S}_{\mathbb{F}^{W}}^{\infty}([0,1];\mathbb{R}) \times \bigcap_{p>1} \mathcal{H}_{\mathbb{F}^{W}}^{p}([0,1];H_{0}^{-1}(\Lambda))$ 
    which uniquely solves BSDE (\ref{sec-3-variational-eq}) according to Theorem 10 in \cite{Briand-Confortola-2008}.
    Combing this with Proposition \ref{sec-3-prop-1}, we further deduce that $(\mathcal{Y}^{(\delta)},\mathcal{Z}^{(\delta)})$ is also in 
    $\mathcal{S}_{\mathbb{F}^{W}}^{\infty}([0,1];\mathbb{R}) \times \bigcap_{p>1} \mathcal{H}_{\mathbb{F}^{W}}^{p}([0,1];H_{0}^{-1}(\Lambda))$,
    where
    \begin{align*}
        \mathcal{Y}^{(\delta)} := Y^{\varepsilon, \delta F}(\rho) - Y^{\varepsilon,0}(\rho) - \delta \widehat{Y}^{\varepsilon, F}(\rho); \\
        \mathcal{Z}^{(\delta)} := Z^{\varepsilon, \delta F}(\rho) - Z^{\varepsilon,0}(\rho) - \delta \widehat{Z}^{\varepsilon, F}(\rho).
    \end{align*}
    Hence proving (\ref{sec-3-Taylor-expansion}) is equivalent to proving $\mathcal{Y}_{0}^{(\delta)} = o(\delta)$. 
    
    For any $(\varepsilon, \delta) \in (0,1] \times [-1,0) \cup (0,1]$, it follows from (\ref{sec-3-quadratic-BSDE}) and (\ref{sec-3-variational-eq}) that $(\mathcal{Y}^{(\delta)},\mathcal{Z}^{(\delta)})$ satisfies the BSDE
    \[
        \mathcal{Y}_{t}^{(\delta)} = \int_{t}^{1} g^{(\varepsilon, \delta)}\left( s, \mathcal{Z}_{s}^{(\delta)} \right) \mathrm{d}s  - \int_{t}^{1} \mathcal{Z}_{s}^{(\delta)} \mathrm{d} W_{s}, \quad t \in [0,1],
    \]
    where, for any $(t,z) \in [0,1] \times H_{0}^{-1}(\Lambda)$,
    \[
    g^{(\varepsilon, \delta)}(t,z) := f^{(\varepsilon, 0)}\left( t, z \right) + \frac{\alpha}{2} \left\Vert Z_{t}^{\varepsilon, \delta F}(\rho) - Z_{t}^{\varepsilon,0}(\rho) \right\Vert _{H_{0}^{-1}(\Lambda)}^{2}.
    \]
    Recall the constants $\overline{p}$ and its conjugate exponent $\overline{p}^{\ast}$ that is defined in the proof of Proposition \ref{sec-3-prop-1}.
    Since $\mathcal{Z}^{(\delta)} \in \bigcap_{p>1} \mathcal{H}_{\mathbb{F}^{W}}^{p}([0,1];H_{0}^{-1}(\Lambda))$, we deduce from Lemma 7 in \cite{Briand-Confortola-2008}, Proposition \ref{sec-3-prop-1}, (\ref{sec-3-BMO-Lip-uniform-bound}), (\ref{sec-3-constant-K-compare}) that,
    for any $p > \overline{p}^{\ast}$,
    \[
    \begin{array}[c]{rl}
        & \left\Vert \mathcal{Y}^{(\delta)} \right\Vert _{\mathcal{S}_{\mathbb{F}^{W}}^{2p}([0,1];\mathbb{R})} \\
        \leqslant & \displaystyle C_{p} K \left( \frac{p}{p-1}, \left\Vert \Phi^{(\varepsilon, 0)} \cdot B^{(1)} \right\Vert _{\mathrm{BMO}} \right) ^{\frac{p-1}{p}} \left\Vert \int_{0}^{1} \left\vert g^{(\varepsilon, \delta)}(s,0) \right\vert \mathrm{d}s \right\Vert _{L^{3p}(\mathcal{F}_{1}^{W}; \mathbb{R})}\\
        \leqslant & \displaystyle C_{p} \mathcal{K} ^{\frac{p-1}{p}} \left\Vert \int_{0}^{1} \left\Vert Z_{s}^{\varepsilon, \delta F}(\rho) - Z_{s}^{\varepsilon,0}(\rho) \right\Vert _{H_{0}^{-1}(\Lambda)}^{2} \mathrm{d}s \right\Vert _{L^{3p}(\mathcal{F}_{1}^{W}; \mathbb{R})}\\
        = & \displaystyle C_{p} \mathcal{K} ^{\frac{p-1}{p}}  \left\Vert Z^{\varepsilon, \delta F}(\rho) - Z^{\varepsilon,0}(\rho) \right\Vert _{\mathcal{H}_{\mathbb{F}^{W}}^{6p}([0,1];H_{0}^{-1}(\Lambda))}^{2} \\
        \leqslant & C \left\vert \delta \right\vert ^{2},
    \end{array}
    \]
    with the notation $\mathcal{K} := K \left( \frac{p}{p-1}, \sqrt{2} e^{\frac{|\alpha|}{2} (b_{F} + b_{\rho} b_{\psi} C_{\beta})} \right)$, where the constant $C_{p}>0$ depends only on $p$, and the constant $C$ depends on $p$, $\alpha$, $\beta$, $b_{\rho}$, $b_{F}$, and $b_{\psi}$.
    Therefore we have $\left\vert \mathcal{Y}_{0}^{(\delta)} \right\vert \leqslant \left\Vert \mathcal{Y}^{(\delta)} \right\Vert _{\mathcal{S}_{\mathbb{F}^{W}}^{2p}([0,1];\mathbb{R})} \leqslant C \left\vert \delta \right\vert ^{2}$ and hence $\mathcal{Y}_{0}^{(\delta)} = o(\delta)$ as $\delta \rightarrow 0$.
\end{proof}

For any $\varepsilon \in (0,1]$ and $\rho \in \mathcal{D}_{0}(\Lambda)$, define a stochastic exponential on $[0,1]$:
\[
    \Gamma_{t}^{\varepsilon}(\rho) = \mathcal{E}\left( \alpha Z^{(\varepsilon,0)}(\rho) \cdot W \right)_{t} = e^{ \alpha \int_{0}^{t} Z_{s}^{(\varepsilon,0)}(\rho) \mathrm{d} W_{s} - \frac{\alpha^{2}}{2} \int_{0}^{t} \left\Vert Z_{s}^{(\varepsilon,0)}(\rho) \right\Vert _{H_{0}^{-1}(\Lambda)}^{2} \mathrm{d}s }.
\]
Due to (\ref{sec-3-BMO-Lip-uniform-bound}), $\Gamma^{\varepsilon}(\rho)$ is a uniform integrable martingale. Then, by the reverse H\"{o}lder inequality, 
\begin{equation}
    \label{sec-3-Gamma-varepsilon-reverse-Holder}
    \sup_{\varepsilon \in (0,1]} \mathbb{E}\left[ \left( \frac{\Gamma_{1}^{\varepsilon}(\rho)}{\Gamma_{\tau}^{\varepsilon}(\rho)} \right) ^{p} \mid \mathcal{F}_{\tau}^{W} \right] \leqslant K\left( p, \sqrt{2} e^{\frac{|\alpha|}{2} b_{\rho} b_{\psi} C_{\beta}} \right) < +\infty, \quad \mathbb{P}-a.s.
\end{equation}
for any $p \in [1, \overline{p})$ and all $\mathbb{F}^{W}$-stopping times $\tau$ taking values in $[0,1]$. Moreover, using It\^{o}'s formula one can check that $\Gamma^{\varepsilon}(\rho)$ satisfies the (forward) scalar stochastic differential equation:
\begin{equation}
    \label{sec-3-SDE-Gamma-varepsilon}
    \Gamma_{t}^{\varepsilon}(\rho) = 1 + \alpha \int_{0}^{t} \Gamma_{s}^{\varepsilon}(\rho) Z_{s}^{(\varepsilon,0)}(\rho) \mathrm{d} W_{s}, \quad t \in [0,1].
\end{equation}
The following lemma shows that $\widehat{Y}^{\varepsilon, F}(\rho)$ can be expressed explicitly by $\Gamma^{\varepsilon}(\rho)$.

\begin{lemma}
    \label{sec-3-lemma-3}
    Let $\rho \in \mathcal{D}_{0}(\Lambda)$, $F \in \mathcal{B}_{b}(H_{0}^{-1}(\Lambda))$ be given arbitrarily. If $\beta^{2} \in [0,2)$, then, for any $\varepsilon \in (0,1]$,
    \begin{equation}
        \label{sec-3-hatY-explicitly-express}
        \widehat{Y}_{t}^{\varepsilon, F}(\rho) = \mathbb{E} \left[ \left( \frac{\Gamma_{1}^{\varepsilon}(\rho)}{\Gamma_{t}^{\varepsilon}(\rho)} \right) F(W_{1}) \mid \mathcal{F}_{t}^{W} \right], \quad \forall t \in [0,1],
    \end{equation}
    and particularly $\mu_{\mathrm{SG}}^{\rho, \varepsilon}( F ) = \mathbb{E} \left[ \Gamma_{1}^{\varepsilon}(\rho) F(W_{1}) \right]$.
\end{lemma}

\begin{proof}
    Using It\^{o}'s formula yields that
    \begin{align}
        \label{sec-3-Gamma-hatY-eq}
        \Gamma_{t}^{\varepsilon}(\rho) \widehat{Y}_{t}^{\varepsilon, F}(\rho) = & \ \Gamma_{1}^{\varepsilon}(\rho) F(W_{1}) - \alpha \int_{t}^{1} \Gamma_{s}^{\varepsilon}(\rho) \widehat{Y}_{s}^{\varepsilon, F}(\rho) Z_{s}^{(\varepsilon,0)}(\rho) \mathrm{d} W_{s} \\
       \nonumber & - \int_{t}^{1}  \Gamma_{s}^{\varepsilon}(\rho) \widehat{Z}_{s}^{\varepsilon, F}(\rho) \mathrm{d} W_{s}.
    \end{align}
    The stochastic integrals in (\ref{sec-3-Gamma-hatY-eq}) are true martingales. Due to the B-D-G inequality, H\"{o}lder's inequality and Doob's maximal inequality, it suffices to prove that, say, for the first stochastic integral in (\ref{sec-3-Gamma-hatY-eq}),
    \[
    \begin{array}[c]{rl}
        & \left\Vert \Gamma^{\varepsilon}(\rho) \widehat{Y}^{\varepsilon, F}(\rho) Z^{(\varepsilon,0)}(\rho) \right\Vert _{\mathcal{H}_{\mathbb{F}^{W}}^{2}([0,1];H_{0}^{-1}(\Lambda))} \\
        \leqslant & \left\Vert \Gamma^{\varepsilon}(\rho) \right\Vert _{\mathcal{S}_{\mathbb{F}^{W}}^{p}([0,1];\mathbb{R})} \left\Vert \widehat{Y}^{\varepsilon, F}(\rho) Z^{(\varepsilon,0)}(\rho) \right\Vert _{\mathcal{H}_{\mathbb{F}^{W}}^{p^{\ast}}([0,1];H_{0}^{-1}(\Lambda))} \\
        \leqslant & p^{\ast} \left\Vert \Gamma_{1}^{\varepsilon}(\rho) \right\Vert _{L^{p}(\mathcal{F}_{1}^{W};\mathbb{R})} \left\Vert \widehat{Y}^{\varepsilon, F}(\rho) \right\Vert _{\mathcal{S}_{\mathbb{F}^{W}}^{2p^{\ast}}([0,1];\mathbb{R})} \left\Vert Z^{(\varepsilon,0)}(\rho) \right\Vert _{\mathcal{H}_{\mathbb{F}^{W}}^{2p^{\ast}}([0,1];H_{0}^{-1}(\Lambda))}
    \end{array} 
    \]
    with $p = (1 + \overline{p})/2$, which is finite due to Lemma \ref{sec-3-lemma-2} and (\ref{sec-3-Gamma-varepsilon-reverse-Holder}). The same approach can be applied to the second stochastic integral in (\ref{sec-3-Gamma-hatY-eq}).
    As $\Gamma_{t}^{\varepsilon}(\rho) > 0$ for every $t \in [0,1]$, $\mathbb{P}$-a.s., (\ref{sec-3-hatY-explicitly-express}) is obtained by dividing the both sides of (\ref{sec-3-Gamma-hatY-eq}) by $\Gamma_{t}^{\varepsilon}(\rho)$ and then taking $\mathbb{E} [\cdot \mid \mathcal{F}_{t}^{W}]$.
    The last claim follows from (\ref{sec-3-SG-variational-expression}) immediately.
\end{proof}

\subsection{Weak convergence of the approximate sine-Gordon measure}
As mentioned before, our goal is to sending $\varepsilon$ to $0$ to see whether
there exists a nontrivial limit of $\mu_{\mathrm{SG}}^{\rho, \varepsilon}$ (in the sense of weak convergence).
Since Lemma \ref{sec-3-lemma-2} guarantees the validity of (\ref{sec-3-SG-variational-expression}), telling us that $\mu_{\mathrm{SG}}^{\rho, \varepsilon}$ can be characterized by 
the initial value of the solution $\widehat{Y}^{\varepsilon, F}(\rho)$ to the BSDE (\ref{sec-3-variational-eq}), it is equivalent to studying the asymptotic behavior of $\widehat{Y}^{\varepsilon, F}(\rho)$ as $\varepsilon \rightarrow 0^{+}$ in suitable sense.

To this end, it is necessary to obtain the convergence result of $Z^{\varepsilon, 0}(\rho)$ as $\varepsilon \rightarrow 0^{+}$,
because $\widehat{Y}^{\varepsilon, F}(\rho)$ relies on $Z^{\varepsilon, 0}(\rho)$ through BSDE (\ref{sec-3-variational-eq}) for each $\varepsilon \in (0,1]$.
So we have to focus on investigating the asymptotic property of the solution (when $F = 0$) to quadratic BSDE (\ref{sec-3-quadratic-BSDE}) as $\varepsilon \rightarrow 0^{+}$,
which is completely determined by the asymptotic behavior of terminal condition $\left( \left[\left[\cos(\beta W_{1}^{\varepsilon})\right]\right], \rho \right) _{L^{2}(\Lambda_{\varepsilon};\mu)}$.
It should be emphasized that $\left( \left[\left[\cos(\beta W_{1}^{\varepsilon})\right]\right], \rho \right) _{L^{2}(\Lambda_{\varepsilon};\mu)}$ is actually the real part of the $\varepsilon$-mollified complex (more precisely, imaginary)
multiplicative chaos associated with $W_{1}$ (recall that with our notation $G_{\Lambda}^{\varepsilon}(x) = \mathbb{E} \left[ (W_{1}^{\varepsilon}(x))^{2} \right]$):
\[
M_{\varepsilon}^{(\beta)}(\rho) := \int_{\Lambda_{\varepsilon}} \exp \left\{ i \beta W_{1}^{\varepsilon}(x) + \frac{\beta^{2}}{2} G_{\Lambda}^{\varepsilon}(x) \right\} \rho(x) \mu(\mathrm{d}x), \quad \forall \rho \in \mathcal{D}_{0}(\Lambda),
\]
where we recall that $\mu$ is of the form $\mu(\mathrm{d}x) = \psi(x) \mathrm{d}x$ with $\psi$ being bounded and measurable on $\Lambda$ (see. e.g. \cite{Webb-AAP-2019,Lacoin-AAP-2022} for more details about complex multiplicative chaos).
Theorem 2.1 in \cite{Lacoin-AAP-2022} indicates that, if $\beta^{2} \in [0,2)$, then the limit $\lim_{\varepsilon \rightarrow 0^{+}} M_{\varepsilon}^{(\beta)}(\rho) = M_{0}^{(\beta)}(\rho)$ exists in $\mathbb{P}$,
and $M_{0}^{(\beta)}(\rho)$ does not depend on the choice of the smoothing kernel $\eta$. 
Thus the cosine (resp. sine) of $W_{1}$ tested against a given $\rho \in \mathcal{D}_{0}(\Lambda)$ (simply denoted by ``$\langle \cos(\beta W_{1}), \rho \rangle _{\mu}$'' (resp. ``$\langle \sin(\beta W_{1}), \rho \rangle _{\mu}$'')) is defined as the real (resp. imaginary) part of $M_{0}^{(\beta)}(\rho)$, or in other words, one has
\begin{align}
    \label{sec-3-cosine-W-converge-prob}
    \langle \cos(\beta W_{1}), \rho \rangle _{\mu} := & \lim\limits _{\varepsilon \rightarrow 0^{+}} \int_{\Lambda_{\varepsilon}} e^{ \frac{\beta^{2}}{2} G_{\Lambda}^{\varepsilon}(x)} \cos(\beta W_{1}^{\varepsilon}(x)) \rho(x) \mu(\mathrm{d}x) \\
    \nonumber = & \lim\limits _{\varepsilon \rightarrow 0^{+}} \left( \left[\left[\cos(\beta W_{1}^{\varepsilon})\right]\right], \rho \right) _{L^{2}(\Lambda_{\varepsilon};\mu)},
\end{align}
\begin{align}
    \label{sec-3-sine-W-converge-prob}
    \langle \sin(\beta W_{1}), \rho \rangle _{\mu} := & \lim\limits _{\varepsilon \rightarrow 0^{+}} \int_{\Lambda_{\varepsilon}} e^{ \frac{\beta^{2}}{2} G_{\Lambda}^{\varepsilon}(x)} \sin(\beta W_{1}^{\varepsilon}(x)) \rho(x) \mu(\mathrm{d}x) \\
    \nonumber = & \lim\limits _{\varepsilon \rightarrow 0^{+}} \left( \left[\left[\sin(\beta W_{1}^{\varepsilon})\right]\right], \rho \right) _{L^{2}(\Lambda_{\varepsilon};\mu)},
\end{align}
where the convergence is in $\mathbb{P}$. It should be noted that the subscript ``$\mu$'' in (\ref{sec-3-cosine-W-converge-prob}) and (\ref{sec-3-sine-W-converge-prob}) highlights the dependence on $\mu$ of the above two limits.

\begin{remark}
    \label{sec-3-remark-1}
    When $\beta^{2} \in [2,4)$, although $M_{\varepsilon}^{(\beta)}(\rho)$ does not converge as a random variable,
    it can be renormalized by subtracting an ever-increasing number of field independent counter-terms \cite{Gallavotti-1985}.
    However it would require some refinement for BSDE (\ref{sec-3-quadratic-BSDE}), which we leave for future work.
\end{remark}

\begin{proposition}
    \label{sec-3-proposition-2}
    Let $\rho \in \mathcal{D}_{0}(\Lambda)$, $F \in \mathcal{B}_{b}(H_{0}^{-1}(\Lambda))$ be given arbitrarily. If $\beta^{2} \in [0,2)$, then for any $p \in [1, +\infty)$,
    \begin{align}
        \label{sec-3-lem-4-eq1}
     \left\Vert Y^{\varepsilon,0}(\rho) - Y^{\varepsilon^{\prime},0}(\rho) \right\Vert _{\mathcal{S}_{\mathbb{F}^{W}}^{p}([0,1];\mathbb{R})} + \left\Vert Z^{\varepsilon,0}(\rho) - Z^{\varepsilon^{\prime},0}(\rho) \right\Vert _{\mathcal{H}_{\mathbb{F}^{W}}^{p}([0,1];H_{0}^{-1}(\Lambda))} \rightarrow 0
    \end{align}
    as $(\varepsilon, \varepsilon^{\prime}) \rightarrow (0,0)^{+}$. In particular, there exists a couple $(\overline{Y}(\rho), \overline{Z}(\rho)) $ belonging to $ \mathcal{S}_{\mathbb{F}^{W}}^{\infty}([0,1];\mathbb{R}) \times \mathcal{H}_{\mathbb{F}^{W}}^{2}([0,1];H_{0}^{-1}(\Lambda))$ such that
    \begin{equation}
        \label{sec-3-lem-4-eq2}
        \lim\limits _{\varepsilon \rightarrow 0^{+}} \left( \left\Vert Y^{\varepsilon,0}(\rho) - \overline{Y}(\rho) \right\Vert _{\mathcal{S}_{\mathbb{F}^{W}}^{2}([0,1];\mathbb{R})} + \left\Vert Z^{\varepsilon,0}(\rho) - \overline{Z}(\rho) \right\Vert _{\mathcal{H}_{\mathbb{F}^{W}}^{2}([0,1];H_{0}^{-1}(\Lambda))} \right) = 0
    \end{equation} 
    and $(\overline{Y}(\rho), \overline{Z}(\rho))$ uniquely solves the quadratic BSDE
    \begin{equation}
        \label{sec-3-quadratic-BSDE-bar}
        \left\{
        \begin{array}[c]{rl}
           \mathrm{d} \overline{Y}_{t}(\rho) = & \displaystyle - \frac{\alpha}{2} \left\Vert \overline{Z}_{t}(\rho) \right\Vert _{H_{0}^{-1}(\Lambda)}^{2} \mathrm{d}t + \overline{Z}_{t}(\rho) \mathrm{d} W_{t}, \quad t \in [0,1),\\
           \overline{Y}_{1}(\rho) = & \langle \cos(\beta W_{1}), \rho \rangle _{\mu}.
        \end{array}
        \right.
    \end{equation}
\end{proposition}

\begin{proof}
By Lemma \ref{sec-3-lemma-1} and the energy inequality for $\mathrm{BMO}$ martingales, we have for any $p \in [1, +\infty)$,
\begin{equation}
    \label{sec-3-Yk-Zk-uniform-bounded}
    \sup\limits _{\varepsilon \in (0,1]} \left( \left\Vert Y^{\varepsilon,0}(\rho) \right\Vert _{\mathcal{S}_{\mathbb{F}^{W}}^{p}([0,1];\mathbb{R})} + \left\Vert Z^{\varepsilon,0}(\rho) \right\Vert _{\mathcal{H}_{\mathbb{F}^{W}}^{p}([0,1];H_{0}^{-1}(\Lambda))} \right) \leqslant C < +\infty,
\end{equation}
where the constant $C>0$ depends only on $\alpha$, $\beta$, $p$, $b_{\rho}$, and $b_{\psi}$.
According to the Vitali convergence theorem, it is thus enough to prove that, as $(\varepsilon, \varepsilon^{\prime}) \rightarrow (0,0)^{+}$,
\[
        \sup_{t \in [0,1]} \left\vert Y_{t}^{\varepsilon,0}(\rho)-Y_{t}^{\varepsilon^{\prime},0}(\rho) \right\vert + \int_{0}^{1} \left\Vert Z_{s}^{\varepsilon,0}(\rho)-Z_{s}^{\varepsilon^{\prime},0}(\rho) \right\Vert _{H_{0}^{-1}(\Lambda)}^{2} \mathrm{d}s \xrightarrow{\mathbb{P}} 0
\]
to get the desired result. To this end, for any $\varepsilon \in (0,1]$, set $\xi_{\varepsilon} := \left( \left[\left[\cos(\beta W_{1}^{\varepsilon})\right]\right], \rho \right) _{L^{2}(\Lambda_{\varepsilon};\mu)}$, $(Y^{\varepsilon},Z^{\varepsilon}):=(Y^{\varepsilon,0}(\rho),Z^{\varepsilon,0}(\rho))$.
Then, for any $(\varepsilon, \varepsilon^{\prime}) \in (0,1] \times (0,1]$ and $t \in [0,1]$, we have
\[
    \widehat{Y}_{t}^{\varepsilon, \varepsilon^{\prime}} = \widehat{\xi}_{\varepsilon, \varepsilon^{\prime}} + \alpha \int_{t}^{1} \left[ \left( Z_{s}^{\varepsilon^{\prime}} , \widehat{Z}_{s}^{\varepsilon, \varepsilon^{\prime}} \right) _{H_{0}^{-1}(\Lambda)} + \frac{1}{2} \left\Vert \widehat{Z}_{s}^{\varepsilon, \varepsilon^{\prime}}  \right\Vert _{H_{0}^{-1}(\Lambda)}^{2} \right] \mathrm{d}s 
    - \displaystyle \int_{t}^{1} \widehat{Z}_{s}^{\varepsilon, \varepsilon^{\prime}} \mathrm{d} W_{s},
\]
where $\widehat{\xi}_{\varepsilon, \varepsilon^{\prime}}:=\xi_{\varepsilon} - \xi_{\varepsilon^{\prime}}$, $\widehat{Y}^{\varepsilon, \varepsilon^{\prime}}:=Y^{\varepsilon}-Y^{\varepsilon^{\prime}}$, $\widehat{Z}^{\varepsilon, \varepsilon^{\prime}}:=Z^{\varepsilon}-Z^{\varepsilon^{\prime}}$.
When $\alpha>0$, it follows from It\^{o}'s formula that
\[
        \begin{array}[c]{rl}
           e^{\alpha \gamma \widehat{Y}_{t}^{\varepsilon, \varepsilon^{\prime}}} = & e^{\alpha \gamma \widehat{\xi}_{\varepsilon, \varepsilon^{\prime}}} + \alpha^{2} \gamma \int_{t}^{1} e^{\alpha \gamma \widehat{Y}_{s}^{\varepsilon, \varepsilon^{\prime}}} \left( Z_{s}^{\varepsilon^{\prime}} , \widehat{Z}_{s}^{\varepsilon, \varepsilon^{\prime}} \right) _{H_{0}^{-1}(\Lambda)} \mathrm{d}s \\
           & + \frac{\alpha^{2} \gamma (1 - \gamma)}{2} \int_{t}^{1} e^{\alpha \gamma \widehat{Y}_{s}^{\varepsilon, \varepsilon^{\prime}}} \left\Vert \widehat{Z}_{s}^{\varepsilon, \varepsilon^{\prime}}  \right\Vert _{H_{0}^{-1}(\Lambda)}^{2} \mathrm{d}s - \alpha \gamma \int_{t}^{1} e^{\alpha \gamma \widehat{Y}_{s}^{\varepsilon, \varepsilon^{\prime}}} \widehat{Z}_{s}^{\varepsilon, \varepsilon^{\prime}} \mathrm{d} W_{s} \\
           \leqslant & \displaystyle e^{\alpha \gamma \widehat{\xi}_{\varepsilon, \varepsilon^{\prime}}} + \alpha^{2} \gamma \int_{t}^{1} e^{\alpha \gamma \widehat{Y}_{s}^{\varepsilon, \varepsilon^{\prime}}} \left( Z_{s}^{\varepsilon^{\prime}} , \widehat{Z}_{s}^{\varepsilon, \varepsilon^{\prime}} \right) _{H_{0}^{-1}(\Lambda)} \mathrm{d}s \\
           & - \alpha \gamma \int_{t}^{1} e^{\alpha \gamma \widehat{Y}_{s}^{\varepsilon, \varepsilon^{\prime}}} \widehat{Z}_{s}^{\varepsilon, \varepsilon^{\prime}} \mathrm{d} W_{s}
        \end{array}
\]
for any $\gamma \in [1,+\infty)$. Hence we further deduce from (\ref{sec-3-SDE-Gamma-varepsilon}) and It\^{o}'s formula that
\[
    \Gamma_{t}^{\varepsilon^{\prime}}(\rho) e^{\alpha \gamma \widehat{Y}_{t}^{\varepsilon, \varepsilon^{\prime}}} \leqslant \Gamma_{1}^{\varepsilon^{\prime}}(\rho) e^{\alpha \gamma \widehat{\xi}_{\varepsilon, \varepsilon^{\prime}}} - \alpha \int_{t}^{1} \Gamma_{s}^{\varepsilon^{\prime}}(\rho) e^{\alpha \gamma \widehat{Y}_{s}^{\varepsilon, \varepsilon^{\prime}}} \left( Z_{s}^{\varepsilon^{\prime}} + \gamma \widehat{Z}_{s}^{\varepsilon, \varepsilon^{\prime}} \right) \mathrm{d} W_{s}.
\]
Recall the constant $\overline{p}$ defined in the proof of Proposition \ref{sec-3-prop-1} ($b_{F} = 0$ because of $F = 0$).
Dividing both sides of the above inequality by $\Gamma_{t}^{\varepsilon^{\prime}}(\rho)$ and then taking $\mathbb{E} [\cdot \mid \mathcal{F}_{t}^{W}]$,
we apply the fundamental inequality $\log(x) \leqslant x$, H\"{o}lder's inequality and the reverse H\"{o}lder inequality to $\mathcal{E}\left( \alpha Z^{\varepsilon^{\prime}} \cdot W \right)$ to deduce
\begin{align*}
    Y_{t}^{\varepsilon}-Y_{t}^{\varepsilon^{\prime}} \leqslant & \frac{1}{\alpha \gamma} \log \mathbb{E} \left[ e^{\alpha \gamma \left\vert \xi_{\varepsilon} - \xi_{\varepsilon^{\prime}} \right\vert} \left( \frac{\Gamma_{1}^{\varepsilon^{\prime}}(\rho)}{\Gamma_{t}^{\varepsilon^{\prime}}(\rho)} \right) \mid \mathcal{F}_{t}^{W} \right] \\
    \leqslant & \frac{1}{\alpha \gamma} \left( \mathbb{E} \left[ \left( \frac{\Gamma_{1}^{\varepsilon^{\prime}}(\rho)}{\Gamma_{t}^{\varepsilon^{\prime}}(\rho)} \right) ^{\frac{\overline{p} + 1}{2}} \mid \mathcal{F}_{1}^{W} \right] \right) ^{\frac{2}{\overline{p} + 1}} 
    \left( \mathbb{E} \left[ e^{ \frac{\alpha \gamma (\overline{p}+1)}{\overline{p}-1} \left\vert \xi_{\varepsilon} - \xi_{\varepsilon^{\prime}} \right\vert } \mid \mathcal{F}_{t}^{W} \right] \right) ^{\frac{\overline{p} - 1}{\overline{p} + 1}} \\
    \leqslant & \frac{1}{\alpha \gamma} K\left( \frac{1 + \overline{p}}{2}, \sqrt{2} e^{\frac{\alpha}{2} b_{\rho} b_{\psi} C_{\beta}} \right) \left( \mathbb{E} \left[ e^{ \frac{\alpha \gamma (\overline{p}+1)}{\overline{p}-1} \left\vert \xi_{\varepsilon} - \xi_{\varepsilon^{\prime}} \right\vert } \mid \mathcal{F}_{t}^{W} \right] \right) ^{\frac{\overline{p} - 1}{\overline{p} + 1}}
\end{align*}
(reminding that $\alpha>0$), where the last inequality is due to (\ref{sec-3-BMO-Lip-uniform-bound}), (\ref{sec-3-constant-K-compare}), and 
\begin{equation}
    \label{sec-3-Zk-BMO-uniform-bounded}
    \left\Vert \alpha Z^{\varepsilon^{\prime}} \cdot W \right\Vert _{\mathrm{BMO}} = \left\Vert \Phi^{(\varepsilon^{\prime},0)} \cdot B^{(1)} \right\Vert _{\mathrm{BMO}} \leqslant \sqrt{2} e^{\frac{\alpha}{2} b_{\rho} b_{\psi} C_{\beta}}, \quad \forall \varepsilon \in (0,1]. 
\end{equation}
The same analysis can be applied to $Y^{\varepsilon^{\prime}}-Y^{\varepsilon}$, whence we obtain
\begin{equation}
    \label{sec-3-Yn-Ym-abs-dominate}
    \left\vert Y_{t}^{\varepsilon}-Y_{t}^{\varepsilon^{\prime}} \right\vert \leqslant \frac{C}{\alpha \gamma} \mathbb{E} \left[ e^{\alpha \gamma (\bar{p}+1)(\bar{p}-1)^{-1} \left\vert \xi_{\varepsilon} - \xi_{\varepsilon^{\prime}} \right\vert} \mid \mathcal{F}_{t}^{W} \right]
\end{equation}
for $\alpha > 0$ with $C>0$ a constant depending only on $\alpha$, $\beta$, $b_{\rho}$, and $b_{\psi}$.
When $\alpha<0$, taking the minus sign on both sides of BSDE (\ref{sec-3-quadratic-BSDE}) and treating $(-Y^{\varepsilon,0}(\rho), -Z^{\varepsilon,0}(\rho))$ as a new pair of solutions, we can prove (\ref{sec-3-Yn-Ym-abs-dominate}) by repeating the previous procedures in which $\alpha$ is substituted by $-\alpha$.
Consequently, the right-hand side of (\ref{sec-3-Yn-Ym-abs-dominate}), which is a uniform integrable martingale since $\sup_{\varepsilon \in (0,1]} \left\Vert \xi_{\varepsilon} \right\Vert _{L^{\infty}(\mathcal{F}_{1}^{W};\mathbb{R})} < +\infty$,
admits a continuous modification denoted by $(M_{t}^{\varepsilon, \varepsilon^{\prime}})_{t \in [0,1]}$ due to the continuity of $\mathbb{F}^{W}$.
Then, for any fixed $\delta > 0$, we deduce from (\ref{sec-3-Yn-Ym-abs-dominate}) and Doob's maximal inequality that
\begin{align*}
    \mathbb{P} \left( \sup\limits _{t \in [0,1]} \left\vert Y_{t}^{\varepsilon}-Y_{t}^{\varepsilon^{\prime}} \right\vert \geqslant \delta \right) \leqslant & \mathbb{P} \left( \sup\limits _{t \in [0,1]} M_{t}^{\varepsilon, \varepsilon^{\prime}} \geqslant \delta \right) \\
    \leqslant & \frac{C}{|\alpha| \delta \gamma} \mathbb{E} \left[ e^{|\alpha| \gamma (\bar{p}+1)(\bar{p}-1)^{-1} \left\vert \xi_{\varepsilon} - \xi_{\varepsilon^{\prime}} \right\vert} \right],
\end{align*}
where the constant $C>0$ depends only on $\alpha$, $\beta$, $b_{\rho}$, and $b_{\psi}$. By $\sup_{\varepsilon \in (0,1]} \left\Vert \xi_{\varepsilon} \right\Vert _{L^{\infty}(\mathcal{F}_{1}^{W};\mathbb{R})} < +\infty$ and (\ref{sec-3-cosine-W-converge-prob}),
applying the dominated convergence to the above inequality yields
\[
\limsup\limits _{(\varepsilon, \varepsilon^{\prime}) \rightarrow (0,0)^{+}} \mathbb{P} \left( \sup\limits _{t \in [0,1]} \left\vert Y_{t}^{\varepsilon}-Y_{t}^{\varepsilon^{\prime}} \right\vert \geqslant \delta \right) \leqslant \frac{C}{|\alpha| \delta \gamma}, \quad \forall \gamma \in [1, +\infty),
\]
and then we obtain $\sup_{t \in [0,1]} \left\vert Y_{t}^{\varepsilon}-Y_{t}^{\varepsilon^{\prime}} \right\vert \xrightarrow{\mathbb{P}} 0$ by sending $\gamma$ to $+\infty$.
To prove $\int_{0}^{1} \left\Vert Z_{s}^{\varepsilon}-Z_{s}^{\varepsilon^{\prime}} \right\Vert _{H_{0}^{-1}(\Lambda)}^{2} \mathrm{d}s \xrightarrow{\mathbb{P}} 0$, using It\^{o}'s formula and taking expectation to get
\begin{align*}
    & \mathbb{E} \left[ \int_{0}^{1} \left\Vert Z_{s}^{\varepsilon}-Z_{s}^{\varepsilon^{\prime}} \right\Vert _{H_{0}^{-1}(\Lambda)}^{2} \mathrm{d}s \right] \\
    \leqslant & \mathbb{E} \left[ |\alpha| \sup\limits _{t \in [0,1]} \left\vert Y_{t}^{\varepsilon}-Y_{t}^{\varepsilon^{\prime}} \right\vert \int_{0}^{1} \left( \left\Vert Z_{s}^{\varepsilon} \right\Vert _{H_{0}^{-1}(\Lambda)}^{2} + \left\Vert Z_{s}^{\varepsilon^{\prime}} \right\Vert _{H_{0}^{-1}(\Lambda)}^{2} \right) \mathrm{d}s + \left\vert \xi_{\varepsilon}-\xi_{\varepsilon^{\prime}} \right\vert ^{2} \right]\\
    \leqslant & |\alpha| \left( \mathbb{E} \left[ \sup\limits _{t \in [0,1]} \left\vert Y_{t}^{\varepsilon}-Y_{t}^{\varepsilon^{\prime}} \right\vert ^{2} \right] \right) ^{\frac{1}{2}} \left( \mathbb{E} \left[ \left( \int_{0}^{1} \left( \left\Vert Z_{s}^{\varepsilon} \right\Vert _{H_{0}^{-1}(\Lambda)}^{2} + \left\Vert Z_{s}^{\varepsilon^{\prime}} \right\Vert _{H_{0}^{-1}(\Lambda)}^{2} \right) \mathrm{d}s \right) ^{2} \right] \right) ^{\frac{1}{2}} \\
    & + \mathbb{E} \left[ \left\vert \xi_{\varepsilon}-\xi_{\varepsilon^{\prime}} \right\vert ^{2} \right],
\end{align*}
where the last inequality is due to the Cauchy-Schwartz. 
As we have proved $\sup_{t \in [0,1]} | Y_{t}^{\varepsilon}-Y_{t}^{\varepsilon^{\prime}} | \xrightarrow{\mathbb{P}} 0$, by (\ref{sec-3-Yk-Zk-uniform-bounded}) and the dominated convergence we obtain $\lim_{(\varepsilon, \varepsilon^{\prime}) \rightarrow (0,0)^{+}} \mathbb{E} \left[ \int_{0}^{1} || Z_{s}^{\varepsilon}-Z_{s}^{\varepsilon^{\prime}} || _{H_{0}^{-1}(\Lambda)}^{2} \mathrm{d}s \right] = 0$.
Consequently, (\ref{sec-3-lem-4-eq1}) follows immediately from the Markov inequality directly.

To prove (\ref{sec-3-lem-4-eq2}), on the one hand, we introduce a sequence of solutions to (\ref{sec-3-quadratic-BSDE}) along $\{ \varepsilon_{n} \} _{n \in \mathbb{N}_{+}}$ such that $\varepsilon_{n} = 2^{-n}$.
It follows from (\ref{sec-3-lem-4-eq1}) with $p=2$ that $\{ ( Y^{\varepsilon_{n}}, Z^{\varepsilon_{n}} ) \} _{n \in \mathbb{N}_{+}}$ is Cauchy in the corresponding spaces,
whence we define $(\overline{Y}(\rho), \overline{Z}(\rho))$ as the limit of $\{ ( Y^{\varepsilon_{n}}, Z^{\varepsilon_{n}} ) \} _{n \in \mathbb{N}_{+}}$ in $\mathcal{S}_{\mathbb{F}^{W}}^{2}([0,1];\mathbb{R}) \times \mathcal{H}_{\mathbb{F}^{W}}^{2}([0,1];H_{0}^{-1}(\Lambda))$.
On the other hand, let $\{ \varepsilon_{n}^{\prime} \} _{n \in \mathbb{N}_{+}} \subset (0,1]$ be another sequence tending to $0$ as $n \rightarrow \infty$, and $\{ ( Y^{\varepsilon_{n}^{\prime}}, Z^{\varepsilon_{n}^{\prime}} ) \} _{n \in \mathbb{N}_{+}}$ be the corresponding sequence of solutions to (\ref{sec-3-quadratic-BSDE}) along $\{ \varepsilon_{n}^{\prime} \} _{n \in \mathbb{N}_{+}}$.
By the triangle inequality and (\ref{sec-3-lem-4-eq1}) with $p=2$, we deduce that
\begin{align*}
    & \limsup\limits _{n \rightarrow \infty} \left( \left\Vert Y^{\varepsilon_{n}^{\prime}} - \overline{Y}(\rho) \right\Vert _{\mathcal{S}_{\mathbb{F}^{W}}^{2}([0,1];\mathbb{R})} + \left\Vert Z^{\varepsilon_{n}^{\prime}} - \overline{Z}(\rho) \right\Vert _{\mathcal{H}_{\mathbb{F}^{W}}^{2}([0,1];H_{0}^{-1}(\Lambda))} \right) \\
        \leqslant & \lim\limits _{n \rightarrow \infty} \left( \left\Vert Y^{\varepsilon_{n}^{\prime}} - Y^{\varepsilon_{n}} \right\Vert _{\mathcal{S}_{\mathbb{F}^{W}}^{2}([0,1];\mathbb{R})} + \left\Vert Z^{\varepsilon_{n}^{\prime}} - Z^{\varepsilon_{n}} \right\Vert _{\mathcal{H}_{\mathbb{F}^{W}}^{2}([0,1];H_{0}^{-1}(\Lambda))} \right) \\
        & + \lim\limits _{n \rightarrow \infty} \left( \left\Vert Y^{\varepsilon_{n}} - \overline{Y}(\rho) \right\Vert _{\mathcal{S}_{\mathbb{F}^{W}}^{2}([0,1];\mathbb{R})} + \left\Vert Z^{\varepsilon_{n}} - \overline{Z}(\rho) \right\Vert _{\mathcal{H}_{\mathbb{F}^{W}}^{2}([0,1];H_{0}^{-1}(\Lambda))} \right)
        = 0,
\end{align*}
from which (\ref{sec-3-lem-4-eq2}) follows immediately due to the arbitrariness of $\{ \varepsilon_{n}^{\prime} \} _{n \in \mathbb{N}_{+}}$.

To prove the last claim, we first deduce from Lemma \ref{sec-3-lemma-1} that 
$$
\sup_{\varepsilon \in (0,1]} || Y^{\varepsilon} || _{\mathcal{S}_{\mathbb{F}^{W}}^{\infty}([0,1];\mathbb{R})} \leqslant b_{\rho} b_{\psi} C_{\beta}.
$$
Applying (\ref{sec-3-lem-4-eq2}) along $\{ \varepsilon_{n}\}_ {n \in \mathbb{N}_{+}}$, there is a subsequence $\{ \varepsilon_{n_{j}} \}_{j=1}^{\infty}$ such that 
$$
\lim\limits _{j \rightarrow \infty} \sup_{t \in [0,1]} \left\vert Y_{t}^{\varepsilon_{n_{j}}} - \overline{Y}_{t}(\rho) \right\vert = 0, \quad \mathbb{P}-a.s.,
$$
which implies that $\overline{Y}(\rho) \in \mathcal{S}_{\mathbb{F}^{W}}^{\infty}([0,1];\mathbb{R})$.
Now we show that $(\overline{Y}(\rho), \overline{Z}(\rho))$ uniquely solves (\ref{sec-3-quadratic-BSDE-bar}). On the one hand, we have
\begin{align}
        \label{sec-3-limit-LN-marting-part}
        \lim\limits _{\varepsilon \rightarrow 0^{+}} \left\Vert \int_{0}^{1} \left( Z_{s}^{\varepsilon} - \overline{Z}_{s}(\rho) \right) \mathrm{d} W_{s} \right\Vert _{L^{2}(\mathcal{F}_{1}^{W};\mathbb{R})}^{2}
        = & \lim\limits _{\varepsilon \rightarrow 0^{+}}\left\Vert Z^{\varepsilon} - \overline{Z}(\rho) \right\Vert _{\mathcal{H}_{\mathbb{F}^{W}}^{2}([0,1];H_{0}^{-1}(\Lambda))}^{2} \\
       \nonumber = & \ 0
\end{align}
    due to the It\^{o} isometry (Proposition \ref{sec-2-proposition-3}).
    On the other hand, the Cauchy-Schwartz inequality results to
\begin{align*}
    & \left\vert \int_{0}^{1}  \left\Vert Z_{s}^{\varepsilon} \right\Vert _{H_{0}^{-1}(\Lambda)}^{2} \mathrm{d}s - \int_{0}^{1} \left\Vert \overline{Z}_{s}(\rho) \right\Vert _{H_{0}^{-1}(\Lambda)}^{2} \mathrm{d}s \right\vert \\
    \leqslant & 2 \int_{0}^{1} \left\Vert \overline{Z}_{s}(\rho) \right\Vert _{H_{0}^{-1}(\Lambda)} \left\Vert Z_{s}^{\varepsilon} - \overline{Z}_{s}(\rho) \right\Vert _{H_{0}^{-1}(\Lambda)} \mathrm{d}s 
     + \int_{0}^{1} \left\Vert Z_{s}^{\varepsilon} - \overline{Z}_{s}(\rho) \right\Vert _{H_{0}^{-1}(\Lambda)}^{2} \mathrm{d}s \\
    \leqslant & 2 \left( \int_{0}^{1} \left\Vert \overline{Z}_{s}(\rho) \right\Vert _{H_{0}^{-1}(\Lambda)}^{2} \right) ^{\frac{1}{2}} \left( \int_{0}^{1} \left\Vert Z_{s}^{\varepsilon} - \overline{Z}_{s}(\rho) \right\Vert _{H_{0}^{-1}(\Lambda)}^{2} \mathrm{d}s \right) ^{\frac{1}{2}} \\
    & + \int_{0}^{1} \left\Vert Z_{s}^{\varepsilon} - \overline{Z}_{s}(\rho) \right\Vert _{H_{0}^{-1}(\Lambda)}^{2} \mathrm{d}s.
\end{align*}
Taking expectation and sending $\varepsilon$ to $0^{+}$, it follows from the Cauchy-Schwartz inequality again that
\begin{equation}
    \label{sec-3-limit-LN-quad-part}
    \lim\limits _{\varepsilon \rightarrow 0^{+}} \mathbb{E} \left[ \left\vert \int_{0}^{1}  \left\Vert Z_{s}^{\varepsilon,0}(\rho) \right\Vert _{H_{0}^{-1}(\Lambda)}^{2} \mathrm{d}s - \int_{0}^{1} \left\Vert \overline{Z}_{s}(\rho) \right\Vert _{H_{0}^{-1}(\Lambda)}^{2} \mathrm{d}s \right\vert \right] = 0.
\end{equation}
Applying (\ref{sec-3-limit-LN-marting-part}) and (\ref{sec-3-limit-LN-quad-part}) along $\{ \varepsilon_{n}\}_ {n \in \mathbb{N}_{+}}$ and extracting a subsequence $\{ \varepsilon_{n_{j}} \}_{j=1}^{\infty}$ if necessary, 
we put $F=0$ and replace $\varepsilon$ with $\varepsilon_{n_{j}}$ in (\ref{sec-3-quadratic-BSDE}) to pass to the limit as $j \rightarrow \infty$ to deduce that $(\overline{Y}(\rho), \overline{Z}(\rho))$ satisfies (\ref{sec-3-quadratic-BSDE-bar}),
where the convergence of $( [[\cos(\beta W_{1}^{\varepsilon_{n_{j}}})]], \rho) _{L^{2}(\Lambda_{\varepsilon};\mu)}$ attributes to (\ref{sec-3-cosine-W-converge-prob}).
Ultimately, the uniqueness in $\mathcal{S}_{\mathbb{F}^{W}}^{\infty}([0,1];\mathbb{R}) \times \mathcal{H}_{\mathbb{F}^{W}}^{2}([0,1];H_{0}^{-1}(\Lambda))$ follows from the same argument as that in the proof of Lemma \ref{sec-3-lemma-1}.
\end{proof}

The theorem below is the main result of this section.

\begin{theorem}
    \label{sec-3-theorem-1}
    Let $\rho \in \mathcal{D}_{0}(\Lambda)$ be given arbitrarily and $\{ \xi_{\varepsilon} \}_{\varepsilon \in (0,1]}$ be a family of $\sigma(W_{1})$-measurable bounded random variables that converge to a $\sigma(W_{1})$-measurable bounded random variable $\xi$ in $\mathbb{P}$.
    If $\beta^{2} \in [0,2)$, then we have
    \begin{equation}
        \label{sec-3-thm-1-general-convergence}
        \lim\limits _{\varepsilon \rightarrow 0^{+}} \mathbb{E} \left[ \Gamma_{1}^{\varepsilon}(\rho) \cdot \xi_{\varepsilon} \right] = \mathbb{E} \left[ \Gamma(\rho) \cdot \xi \right],
    \end{equation}
    where $\Gamma(\rho) := \exp \left\{ \alpha \int_{0}^{1} \overline{Z}_{s}(\rho) \mathrm{d} W_{s} - \frac{\alpha^{2}}{2} \int_{0}^{1} \left\Vert \overline{Z}_{s}(\rho) \right\Vert _{H_{0}^{-1}(\Lambda)}^{2} \mathrm{d}s \right\}$ belongs to $L^{p}(\mathcal{F}_{1}^{W};\mathbb{R})$
    for any $p \in [1, \overline{p})$ with $\overline{p}$ being such that $\kappa(\overline{p}) = \sqrt{2} e^{\frac{|\alpha|}{2} b_{\rho} b_{\psi} C_{\beta}}$, and $\overline{Z}(\rho)$ satisfies (\ref{sec-3-quadratic-BSDE-bar}).
    Moreover, there is a probability measure $\mu_{\mathrm{SG}}^{\rho}$ on $(H_{0}^{-1}(\Lambda), \mathfrak{B}(H_{0}^{-1}(\Lambda)))$ such that $\mu_{\mathrm{SG}}^{\rho,\varepsilon}$ converges weakly to $\mu_{\mathrm{SG}}^{\rho}$ as $\varepsilon \rightarrow 0^{+}$.
\end{theorem}

\begin{proof}

According to the definition of $\Gamma^{\varepsilon}(\rho)$, it follows from (\ref{sec-3-limit-LN-marting-part}), (\ref{sec-3-limit-LN-quad-part}), and the continuity of $x \mapsto e^{x}$ that
$\lim_{\varepsilon \rightarrow 0^{+}}\Gamma_{1}^{\varepsilon}(\rho) = \Gamma(\rho)$ in $\mathbb{P}$ since the convergence in measure is preserved under continuous mappings.
By (\ref{sec-3-Zk-BMO-uniform-bounded}), employing the reverse H\"{o}lder inequality and Fatou's lemma yields that, for any $p \in (1, \bar{p})$,
\begin{equation}
    \label{sec-3-Gamma-rho-p-moments}
    \mathbb{E} \left[ \Gamma(\rho)^{p} \right] \leqslant \liminf\limits _{\varepsilon \rightarrow 0^{+}} \mathbb{E} \left[ \left( \Gamma_{1}^{\varepsilon}(\rho) \right) ^{p}  \right] \leqslant \sup_{\varepsilon \in (0,1]} \mathbb{E}\left[ \left( \Gamma_{1}^{\varepsilon}(\rho) \right) ^{p} \right]
    \leqslant K\left( p, \sqrt{2} e^{\frac{|\alpha|}{2} b_{\rho} b_{\psi} C_{\beta}} \right).
\end{equation}
Then, due to the boundedness of $\{ \xi \}_{\varepsilon \in (0,1]}$ and the convergence $\lim_{\varepsilon \rightarrow 0^{+}} \xi_{\varepsilon} = \xi$ in probability, we deduce (\ref{sec-3-thm-1-general-convergence}) from the Vitali convergence theorem,
and particularly we also have $\lim_{\varepsilon \rightarrow 0^{+}} \mathbb{E} \left[ \left\vert \Gamma_{1}^{\varepsilon}(\rho) - \Gamma(\rho) \right\vert \right] = 0$, from which $\mathbb{E} [\Gamma(\rho)] = \lim_{\varepsilon \rightarrow 0^{+}} \mathbb{E} [\Gamma_{1}^{\varepsilon}(\rho)] = 1$ follows immediately.
Therefore, a new probability measure $\mathbb{P}_{\rho}$ on $(\Omega, \mathcal{F})$ can be defined by $\mathbb{P}_{\rho}(A) := \mathbb{E} \left[ \Gamma(\rho) \mathbf{1}_{A} \right], \ \forall A \in \mathcal{F}$. 
For any $F \in \mathcal{C}_{b}(H_{0}^{-1}(\Lambda))$, according to Lemma \ref{sec-3-lemma-3} and taking $\xi_{\varepsilon} = \xi = F(W_{1})$ for all $\varepsilon \in (0,1]$, we obtain that
\begin{equation}
    \label{sec-3-SG-weak-convergence}
    \lim\limits _{\varepsilon \rightarrow 0^{+}} \int_{H_{0}^{-1}(\Lambda)} F(\phi) \mu_{\mathrm{SG}}^{\rho, \varepsilon}( \mathrm{d}\phi) = \lim\limits _{\varepsilon \rightarrow 0^{+}} \mathbb{E} \left[ \Gamma_{1}^{\varepsilon}(\rho) F(W_{1}) \right] = \mathbb{E} \left[ \Gamma(\rho) F(W_{1}) \right].
\end{equation}
The weak convergence of $\mu_{\mathrm{SG}}^{\rho, \varepsilon}$ follows from (\ref{sec-3-SG-weak-convergence}) if we put $\mu_{\mathrm{SG}}^{\rho} := \mathbb{P}_{\rho} \circ (W_{1})^{-1}$.
\end{proof}

\begin{remark}
    Theorem \ref{sec-3-theorem-1} tells us that our sine-Gordon measure $\mu_{\mathrm{SG}}^{\rho}$ is absolutely continuous with respect to the law of GFF when $\beta^{2} \in [0,2)$, in which case only the Wick-ordered cosine is enough to normalize
    the terminal condition of quadratic BSDE (\ref{sec-3-quadratic-BSDE}). This observation is consistent with the results concerning the finite ultraviolet regime $\beta^{2} \in [0,2)$ in the existing literature (see, e.g., \cite{Barashkov-2022} and references therein).
\end{remark}

\section{Applications in the 2D critical XOR-Ising model and 2D log-gases}
\label{section-4}

In this section, we will discuss the connection between BSDE (\ref{sec-3-quadratic-BSDE-bar}) and the critical planar XOR-Ising model and use Theorem \ref{sec-3-theorem-1} to provide the limit of the sine-Gordon representation for two-dimensional log-gases.
Let $\alpha, \beta$ be two positive parameters and $\mu$ be a locally finite Borel (positive) measure on $\Lambda$ of the form $\mu(\mathrm{d}x) = \psi(x) \mathrm{d}x$ for some nonnegative bounded measurable function $\psi$.
In what follows, $\alpha, \beta$, $\psi$ may appear in some concrete form that will be indicated if necessary.
It should be noted that the constant $C$ appearing in the following proofs may change, and we will indicate its dependence on the parameters if necessary.

\subsection{Partition functions and normalization}

What we are concerned with is the following type of partition functions being (at least formally) given by
\[
    \mathcal{Q}(\alpha, \beta) := \sum\limits _{n=0}^{\infty} \frac{\alpha^{n}}{2^{n}n!} \sum\limits_{\{ \gamma_{k} \}_{k=1}^{n} \in \{ \pm 1 \}^{n}} \int_{\Lambda^{n}} \exp \left\{ - \beta^{2} \sum_{1 \leqslant k < l \leqslant n} \gamma_{k} \gamma_{l} G_{\Lambda}(x_{k}, x_{l}) \right\} \prod\limits _{k = 1}^{n} \mathrm{d}x_{k},
\]
where the interaction is described by the potential in the exponential---$G_{\Lambda}$ in Definition \ref{sec-2-def-Green-func}.
As mentioned in \cite{Lacoin-Rhodes-Vargas-2023-PTRF}, $\mathcal{Q}(\alpha, \beta)$ can characterize the gas of interacting charged particles with potentials of log-type,
for example, the two-dimensional Coulomb gas (For the Yukawa-like potentials, we need to consider the massive Green function $G_{\Lambda}^{(m)}$ with a mass $m > 0$ instead of $G_{\Lambda}$.)
living in a box $\Lambda$ \cite{Gallavotti-Nicolo-1985}, where $\alpha$ stands for particle activity and $\beta^{2}$ for inverse temperature.
By (ii) in Proposition \ref{sec-2-proposition-Green-func}, it should be noted that the above integral over $\Lambda^{n}$ diverges when $\beta^{2} \geqslant 2$, which implies that $\mathcal{Q}(\alpha, \beta)$ also diverges in this case.
For the case $\beta^{2} \in [0,2)$, we will show that the following slight modification of $\mathcal{Q}(\alpha, \beta)$ converges, that is, for any $\rho \in \mathcal{D}_{0}(\Lambda)$ being nonnegative,
\begin{equation}
    \label{sec-4-partition-function-formal}
    \mathcal{Q}_{\rho, \mu}(\alpha, \beta) = \sum\limits _{n=0}^{\infty} \frac{\alpha^{n}}{2^{n}n!} \sum\limits_{\{ \gamma_{k} \}_{k=1}^{n} \in \{ \pm 1 \}^{n}} \int_{\Lambda^{n}} e^{ - \beta^{2} \sum\limits_{1 \leqslant k < l \leqslant n} \gamma_{k} \gamma_{l} G_{\Lambda}(x_{k}, x_{l}) } \prod\limits _{k = 1}^{n} \rho(x_{k}) \mu(\mathrm{d}x_{k})
\end{equation}
is well defined. To this end, applying the global Onsager-type inequality in the case of the GFF on a bounded simply connected domain (see Proposition 3.9 in \cite{Webb-AAP-2019}), there exists a constant $C$ depending only on $\Lambda$ such that, for any $n \in \mathbb{N}_{+}$, $(\gamma_{1},\ldots,\gamma_{n}) \in \{ \pm 1 \}^{n}$, and arbitrary $x_{1},\ldots,x_{n} \in \Lambda$ being distinct,
\[
- \sum_{1 \leqslant k < l \leqslant n} \gamma_{k} \gamma_{l} G_{\Lambda}(x_{k}, x_{l}) \leqslant \frac{1}{2} \sum_{k=1}^{n} \log \left( \frac{1}{\frac{1}{2} \min\limits _{k \neq l} \left\vert x_{k} - x_{l} \right\vert }  \right) + C n.
\]
Then, for any nonnegative $\rho \in \mathcal{D}_{0}(\Lambda)$ and for each $n \in \mathbb{N}_{+}$, we may by using compactness first cover $\mathrm{supp}(\rho)$ with a finite number of balls and then deduce from the above inequality with Lemma 3.10 in \cite{Webb-AAP-2019} that
\begin{align}
    \label{sec-4-partit-summand-dominate}
        0 \leqslant & \sum\limits_{\{ \gamma_{k} \}_{k=1}^{n} \in \{ \pm 1 \}^{n}} \int_{\Lambda^{n}} e^{ - \beta^{2} \sum\limits_{1 \leqslant k < l \leqslant n} \gamma_{k} \gamma_{l} G_{\Lambda}(x_{k}, x_{l}) } \prod\limits _{k = 1}^{n} \rho(x_{k}) \mu(\mathrm{d}x_{k}) \\
        \nonumber \leqslant & \ (2 b_{\rho} b_{\psi})^{n} \int_{\left( \mathrm{supp}(\rho) \right)^{n}} \exp \left\{ \frac{\beta^{2}}{2} \sum\limits_{k=1}^{n} \log \left( \frac{1}{\frac{1}{2} \min\limits _{k \neq l} \left\vert x_{k} - x_{l} \right\vert }  \right) + \beta^{2} C n \right\} \prod\limits _{k = 1}^{n} \mathrm{d}x_{k} \\
        \nonumber \leqslant & \ C^{n} n^{\frac{\beta^{2}}{4}n},
\end{align}
where the constant $C>0$ depends on $\beta$, $b_{\rho}$, $b_{\psi}$, and $\mathrm{supp}(\rho)$, but is independent of $n$.
Thus it follows from a straightforward Stirling estimate that $\mathcal{Q}_{\rho, \mu}(\alpha, \beta)$ converges absolutely for $\beta^{2} \in [0,2)$.

To have an interpretation of (\ref{sec-4-partition-function-formal}) in terms of the sine-Gordon measure constructed in the last section, we need to replace $G_{\Lambda}$ with its regular approximation $G_{\Lambda}^{\varepsilon}$ defined by (\ref{sec-3-def-mollified-Green})
to consider the corresponding limit as $\varepsilon \rightarrow 0^{+}$.
More precisely, it is defined as
\begin{equation}
    \label{sec-4-renormalized-part-func}
    \mathcal{Q}_{\rho, \mu}^{(\varepsilon)}(\alpha,\beta) = \sum\limits _{n=0}^{\infty} \frac{\alpha^{n}}{2^{n}n!} \sum\limits_{\{ \gamma_{k} \}_{k=1}^{n} \in \{ \pm 1 \}^{n}} \int_{\Lambda_{\varepsilon}^{n}} e^{ - \beta^{2} \sum\limits_{1 \leqslant k < l \leqslant n} \gamma_{k} \gamma_{l} G_{\Lambda}^{\varepsilon}(x_{k}, x_{l}) } \prod\limits _{k = 1}^{n} \rho(x_{k}) \mu(\mathrm{d}x_{k}).
\end{equation}
We follow here the nice exposition of \cite{Lacoin-Rhodes-Vargas-2023-PTRF} by using our notation.
Consider $n$ interacting charged particles or $n$ magnetic dipoles with interactions that constrain the spins to point parallel or anti-parallel along a given direction. For any charge (resp. spin) distribution $\{ \gamma_{k} \}_{k=1}^{n}$ taking values in $\{-1,1\}^{n}$ with $k$ positive charges (resp. spin-up states) and $n - k$ negative (resp. spin-down) ones, we use the Fubini theorem to obtain the following correspondence for $\varepsilon \in (0,1]$,
\begin{align}
    \label{sec-4-Qn-varepsilon-express}
        Q_{n}^{(\varepsilon)} := & \sum\limits_{\{ \gamma_{k} \}_{k=1}^{n} \in \{ \pm 1 \}^{n}} \int_{\Lambda_{\varepsilon}^{n}} e^{ - \beta^{2} \sum\limits_{1 \leqslant k < l \leqslant n} \gamma_{k} \gamma_{l} G_{\Lambda}^{\varepsilon}(x_{k}, x_{l}) } \prod\limits _{k = 1}^{n} \rho(x_{k}) \mu(\mathrm{d}x_{k}) \\
       \nonumber = & \ \mathbb{E} \left[ \sum\limits_{\{ \gamma_{k} \}_{k=1}^{n} \in \{ \pm 1 \}^{n}} \int_{\Lambda_{\varepsilon}^{n}} e^{ \sum\limits_{k=1}^{n} \left( i \beta \gamma_{k} W_{1}^{\varepsilon}(x_{k}) + \frac{\beta^{2}}{2} G_{\Lambda}^{\varepsilon}(x_{k}) \right) } \prod\limits _{k = 1}^{n} \rho(x_{k}) \mu(\mathrm{d}x_{k}) \right] \\
       \nonumber = & \ \mathbb{E} \left[ \sum\limits _{k=0}^{n} \binom{n}{k} \left( \int_{\Lambda_{\varepsilon}} e^{i \beta W_{1}^{\varepsilon}(x) + \frac{\beta^{2}}{2} G_{\Lambda}^{\varepsilon}(x)} \rho(x) \mu(\mathrm{d}x) \right) ^{k} \right. \\
       \nonumber & \left. \quad \quad \quad \quad \times \left( \int_{\Lambda_{\varepsilon}} e^{-i \beta W_{1}^{\varepsilon}(x) + \frac{\beta^{2}}{2} G_{\Lambda}^{\varepsilon}(x)} \rho(x) \mu(\mathrm{d}x) \right) ^{n - k} \right] \\
       \nonumber = & \ \mathbb{E} \left[ \left\{ \int_{\Lambda_{\varepsilon}} \left( e^{i \beta W_{1}^{\varepsilon}(x) + \frac{\beta^{2}}{2} G_{\Lambda}^{\varepsilon}(x)} + e^{-i \beta W_{1}^{\varepsilon}(x) + \frac{\beta^{2}}{2} G_{\Lambda}^{\varepsilon}(x)} \right) \rho(x) \mu(\mathrm{d}x) \right\} ^{n} \right]
\end{align}
since $\{ W_{1}^{\varepsilon}(x) \}_{(\varepsilon,x) \in \mathbb{I}}$ is a centered Gaussian field with covariance function $G_{\Lambda}^{\varepsilon}$.
According to this correspondence and noting that for each $n \in \mathbb{N}_{+}$
\[
    \begin{array}[c]{rl}
        & \displaystyle \sum\limits _{k=0}^{n} \frac{\alpha^{k}}{2^{k} k!} \left\vert \int_{\Lambda_{\varepsilon}} \left( e^{i \beta W_{1}^{\varepsilon}(x) + \frac{\beta^{2}}{2} G_{\Lambda}^{\varepsilon}(x)} + e^{-i \beta W_{1}^{\varepsilon}(x) + \frac{\beta^{2}}{2} G_{\Lambda}^{\varepsilon}(x)} \right) \rho(x) \mu(\mathrm{d}x) \right\vert ^{k} \\
        = & \displaystyle \sum\limits _{k=0}^{n} \frac{\alpha^{k}}{2^{k} k!} \left\vert \left( \left[\left[ \cos(\beta W_{1}^{\varepsilon}) \right]\right], \rho \right)_{L^{2}(\Lambda_{\varepsilon};\mu)} \right\vert ^{k}
        \leqslant e^{ \frac{\alpha}{2} b_{\rho} b_{\psi} C_{\beta} }
    \end{array}
\]
because of $\sup_{\varepsilon \in (0,1]} \left\Vert \left( \left[\left[\cos(\beta W_{1}^{\varepsilon})\right]\right], \rho \right) _{L^{2}(\Lambda_{\varepsilon};\mu)} \right\Vert _{L^{\infty}(\mathcal{F}_{1}^{W};\mathbb{R})} \leqslant b_{\rho} b_{\psi} C_{\beta}$, we may apply the dominated convergence theorem to obtain
\begin{align}
    \label{sec-4-charge-partition-GFF-express}
        \mathcal{Q}_{\rho, \mu}^{(\varepsilon)}(\alpha,\beta) = \sum\limits _{n=0}^{\infty} \frac{\alpha^{n}}{2^{n}n!} Q_{n}^{(\varepsilon)} = & \ \mathbb{E} \left[ \exp\left\{ \alpha \left( \left[\left[ \cos(\beta W_{1}^{\varepsilon}) \right]\right], \rho \right)_{L^{2}(\Lambda_{\varepsilon};\mu)} \right\} \right] \\
        \nonumber = & \ \Xi_{\rho, \varepsilon} = e^{\alpha Y_{0}^{\varepsilon, 0}(\rho)},
\end{align}
where the last equation results from Lemma \ref{sec-3-lemma-1}.

\begin{proposition}
    \label{sec-4-proposition-3}
    Let $\alpha > 0$ and $\rho \in \mathcal{D}_{0}(\Lambda)$ be nonnegative. If $\beta^{2} \in [0,2)$, then 
    \[
        \mathcal{Q}_{\rho, \mu}(\alpha,\beta) = \lim\limits _{\varepsilon \rightarrow 0^{+}} \mathcal{Q}_{\rho, \mu}^{(\varepsilon)}(\alpha,\beta) = e^{\alpha \overline{Y}_{0}(\rho) },
    \]
    where $\overline{Y}(\rho)$ satisfies quadratic BSDE (\ref{sec-3-quadratic-BSDE-bar}).
\end{proposition}

\begin{proof}
At first, for any $n \in \mathbb{N}_{+}$, $(\gamma_{1},\ldots,\gamma_{n}) \in \{ \pm 1 \}^{n}$, and arbitrary $x_{1},\ldots,x_{n} \in \Lambda$ being distinct, there is a constant $C>0$ relying on $\mathrm{supp}(\rho)$ and the bound of $g_{\Lambda}$, but being independent of $\varepsilon$ such that
\begin{align*}
    0 \leqslant & \ \mathbf{1}_{\Lambda_{\varepsilon}^{n}}(x_{1},\ldots,x_{n}) \cdot \prod_{k=1}^{n} \rho(x_{k})\psi (x_{k}) \cdot e^{ -\beta^{2} \sum\limits_{1 \leqslant k < l \leqslant n} \gamma_{k} \gamma_{l} G_{\Lambda}^{\varepsilon}(x_{k}, x_{l}) } \\
    \leqslant & \ (b_{\rho} b_{\psi})^{n} \mathbf{1}_{\Lambda_{\varepsilon}^{n} \cap (\mathrm{supp}(\rho))^{n}}(x_{1},\ldots,x_{n}) \cdot e^{ -\beta^{2} \sum\limits_{1 \leqslant k < l \leqslant n} \gamma_{k} \gamma_{l} G_{\Lambda}^{\varepsilon}(x_{k}, x_{l}) } \\
    \leqslant & \ (b_{\rho} b_{\psi})^{n} \mathbf{1}_{(\mathrm{supp}(\rho))^{n}}(x_{1},\ldots,x_{n}) \cdot e^{ \beta^{2} \left[ \sum\limits_{k=1}^{n} \log\left(\frac{1}{\frac{1}{2} \min _{k \neq l} \left\vert x_{k} - x_{l} \right\vert } \right) + C n^{2} \right] }
\end{align*}
holds for all $\varepsilon \in (0,1]$, where the last inequality results from part (ii) of Proposition 3.6 in \cite{Webb-AAP-2019}.
According to Lemma 3.10 in \cite{Webb-AAP-2019}, the quantity in the last line is integrable over $\Lambda^{n}$ with respect to the Lebesgue measure, whence we obtain
\[
    Q_{n}^{(0)} := \lim\limits _{\varepsilon \rightarrow 0^{+}} Q_{n}^{(\varepsilon)}
    = \sum\limits_{\{ \gamma_{k} \}_{k=1}^{n} \in \{ \pm 1 \}^{n}} \int_{\Lambda^{n}} e^{ - \beta^{2} \sum\limits_{1 \leqslant k < l \leqslant n} \gamma_{k} \gamma_{l} G_{\Lambda}(x_{k}, x_{l}) } \prod\limits _{k = 1}^{n} \rho(x_{k}) \mu(\mathrm{d}x_{k})
\]
through the definition of $G_{\Lambda}^{\varepsilon}$ and employing the dominated convergence theorem. On the one hand, it follows from 
$$
\sup_{\varepsilon \in (0,1]} || \left( \left[\left[\cos(\beta W_{1}^{\varepsilon})\right]\right], \rho \right) _{L^{2}(\Lambda_{\varepsilon};\mu)} || _{L^{\infty}(\mathcal{F}_{1}^{W};\mathbb{R})} \leqslant b_{\rho} b_{\psi} C_{\beta}
$$
and (\ref{sec-4-Qn-varepsilon-express}) that $\sup_{\varepsilon \in (0,1]} Q_{n}^{(\varepsilon)} \leqslant (b_{\rho} b_{\psi} C_{\beta})^{n}$ for each $n \in \mathbb{N}_{+}$. On the other hand, by the definition of $Q_{n}^{(0)}$, (\ref{sec-4-partit-summand-dominate}) indicates that $Q_{n}^{(0)} \leqslant C^{n} n^{\frac{\beta^{2}}{4} n}$.
Consequently, the series $\sum_{n=0}^{\infty} \left( \alpha^{n}Q_{n}^{(\varepsilon)}/2^{n}n! \right)$ converges uniformly in $\varepsilon \in [0,1]$. So the limitation commutes with the summation and then we deduce from (\ref{sec-4-charge-partition-GFF-express}) with (\ref{sec-3-lem-4-eq2}) that
\begin{align*}
    \mathcal{Q}_{\rho, \mu}(\alpha,\beta) = \sum\limits _{n=0}^{\infty} \frac{\alpha^{n}}{2^{n}n!} Q_{n}^{(0)} = \lim\limits _{\varepsilon \rightarrow 0^{+}} \sum\limits _{n=0}^{\infty} \frac{\alpha^{n}}{2^{n}n!} Q_{n}^{(\varepsilon)} = & \lim\limits _{\varepsilon \rightarrow 0^{+}} \mathcal{Q}_{\rho, \mu}^{(\varepsilon)}(\alpha,\beta) \\
    = & \lim\limits _{\varepsilon \rightarrow 0^{+}} e^{\alpha Y_{0}^{\varepsilon,0}(\rho)} = e^{\alpha \overline{Y}_{0}(\rho)},
\end{align*}
which accomplishes the proof.
\end{proof}

\begin{remark}
    Proposition \ref{sec-4-proposition-3} provides a probabilistic representation for partition function $\mathcal{Q}_{\rho, \mu}(\alpha,\beta)$ by the solution to quadratic BSDE (\ref{sec-3-quadratic-BSDE-bar}) 
    whose terminal condition is the real part of imaginary multiplicative chaos tested against the nonnegative given $\rho \in \mathcal{D}_{0}(\Lambda)$.
\end{remark}

It is well known that the Ising model plays a central role in equilibrium statistical mechanics, being a standard example of an order-disorder phase transition in dimensions two and above. 
Besides pure mathematical interest, it has found successful applications in several fields in theoretical physics and computer sciences. (See \cite{Book-Baxter-1982} for an extensive introduction to the Ising model.)
To illustrate that (\ref{sec-4-partition-function-formal}) is closely related to the scaling limit of correlation functions of the spin field (By a spin field, we mean a function defined on $\Lambda$ which is constant on these lattice faces, and in each face, it agrees with the value of the spin configuration on that face.)
for the critical planar XOR-Ising model, we first recall the definition of the Ising model (with $+$ boundary conditions), which is the law of a random assignment of $\pm 1$ spins to the dual graph of lattice approximation $\Lambda_{\delta} := \Lambda \cap \delta \mathbb{Z}^{2}$ to $\Lambda$ for $\delta>0$.
Let $\mathbf{F}_{\delta}$ be the set of faces of square lattice $\delta \mathbb{Z}^{2}$ contained in $\Lambda$, and $\partial \mathbf{F}_{\delta}$ be the set of faces in $\delta \mathbb{Z}^{2}$ adjacent to a face in $\mathbf{F}_{\delta}$ but not in $\mathbf{F}_{\delta}$ themselves.
To avoid overlap, we make the convention that the faces are half-open, that is, of the form $\delta ([j,j+1) \times [k,k+1))$ for some $j,k \in \mathbb{Z}$.
Following \cite{Chelkak-Hongler-Izyurov-2015,Webb-AAP-2019}, we define a spin configuration as a function $\sigma_{\delta}: \mathbf{F}_{\delta} \cup \partial \mathbf{F}_{\delta} \mapsto \{ \pm 1 \}$, and the critical planar Ising model with $+$ boundary conditions, the critical inverse temperature $\beta_{c} := \log(1 + \sqrt{2})/2$,
and zero magnetic field to be the following \textit{Ising Boltzmann measure} that describes the probability of occurrence of a spin configuration $\sigma_{\delta}$:
\[
\mathbb{P}_{\Lambda_{\delta}}^{+}(\sigma_{\delta}) := \mathcal{Z}_{\beta_{c}}^{-1} \exp \left\{ \beta_{c} \sum_{\substack{a,b \in \mathbf{F}_{\delta} \cup \partial \mathbf{F}_{\delta}, \\ a \sim b}} \sigma_{\delta}(a) \sigma_{\delta}(b) \right\} \mathbf{1}_{\{ \sigma_{\delta}^{\prime}: \ \sigma_{\delta}^{\prime}(a) = 1, \forall a \in \partial \mathbf{F}_{\delta} \} } (\sigma_{\delta}),
\]
where by $a \sim b$ we mean that $a,b \in \mathbf{F}_{\delta} \cup \partial \mathbf{F}_{\delta}$ are neighboring faces, and $\mathcal{Z}_{\beta_{c}}$ is a normalizing constant. We count each pair $a,b$ of nearest neighbor faces only once.
To consider the spin at an arbitrary point $x \in \Lambda$, by overloading the notation slightly, we define a function $\sigma_{\delta}(x) = \sigma_{\delta}(a)$ if $x \in a$ for all $a \in \mathbf{F}_{\delta}$, and $\sigma_{\delta}(x) = 1$ otherwise. 
We will also write from now on $\mathbb{P}_{\Lambda_{\delta}}^{+}$ for the law of the critical planar Ising model as well as the law of the induced spin field, and by $\mathbb{E}_{\Lambda_{\delta}}^{+}$ we denote the expectation corresponding to $\mathbb{P}_{\Lambda_{\delta}}^{+}$.

Let $\varphi: \Lambda \mapsto \mathbb{H}^{+}$ be a conformal bijection and $\mathcal{C} := 2^{\frac{5}{48}} e^{\frac{3}{2} \zeta^{\prime}(-1)}$ be a lattice-dependent constant
with $\zeta^{\prime}$ denoting the derivative of Riemann's zeta function. By taking $\alpha = 2^{-\frac{1}{2}} \mathcal{C}^{2}$, $\beta = 2^{-\frac{1}{2}}$, and $\psi(x) = \left( \left\vert \varphi^{\prime}(x) \right\vert / 2 \mathrm{Im}\{ \varphi(x) \} \right) ^{\frac{1}{4}}$, we find that (\ref{sec-4-partition-function-formal})
admits an expression that involves the scaling limit of correlation functions of the Ising model as $\delta \rightarrow 0^{+}$.

\begin{remark}
    It should be emphasized that rotating the lattice may change the value of $\mathcal{C}$, say, $\mathcal{C} = 2^{\frac{1}{6}} e^{\frac{3}{2} \zeta^{\prime}(-1)}$ in \cite{Chelkak-Hongler-Izyurov-2015}, since the authors consider actually the square lattice rotated by $\pi/4$ and with diagonal mesh $2\delta$ in which case the lattice spacing is $\sqrt{2} \delta$ instead of $\delta$ in our framework.
\end{remark}

\begin{proposition}
    \label{sec-4-proposition-1}
Let $n \in \mathbb{N}_{+}$ and $x_{1},\ldots,x_{n} \in \Lambda$ be distinct. Then $\mathcal{Q}_{\rho, \mu_{\varphi}}\left( \frac{\mathcal{C}^{2}}{\sqrt{2}}, \frac{1}{\sqrt{2}} \right)$ equals to
\[
\sum\limits_{n=0}^{\infty} \frac{1}{2^{n} n!} \int_{\Lambda^{n}} \left( \mathcal{C}^{n} \cdot \left\langle \sigma_{\varphi(x_{1})}, \cdots, \sigma_{\varphi(x_{n})} \right\rangle _{\mathbb{H}^{+}}^{+} \cdot \prod_{k=1}^{n} \left\vert \varphi^{\prime}(x_{k}) \right\vert ^{\frac{1}{8}} \right)^{2} \prod_{k=1}^{n} \rho(x_{k}) \mathrm{d}x_{k},
\]
where $\mu_{\varphi}( \mathrm{d}x ) := \left( \left\vert \varphi^{\prime}(x) \right\vert / 2 \mathrm{Im}\{ \varphi(x) \} \right) ^{\frac{1}{4}} \mathrm{d}x$, and for any $z_{1},\ldots,z_{n} \in \mathbb{H}^{+}$,
\[
    \left\langle \sigma_{z_{1}}, \cdots, \sigma_{z_{n}} \right\rangle _{\mathbb{H}^{+}}^{+} := \prod_{k=1}^{n} \frac{1}{\left( 2 \mathrm{Im}\{ z_{k} \} \right) ^{\frac{1}{8}} } \cdot \left( 2^{-\frac{n}{2}} \sum\limits_{\{ \gamma_{k} \}_{k=1}^{n} \in \{ \pm 1 \}^{n}} \prod_{1 \leqslant k < l \leqslant n} \left\vert \frac{z_{k} - z_{l}}{z_{k} - \overline{z_{l}}} \right\vert ^{\frac{\gamma_{k} \gamma_{l}}{2}} \right) ^{\frac{1}{2}}
\]
is the scaling limit of $n$-point correlation function on the upper half-plane $\mathbb{H}^{+}$.
\end{proposition}

\begin{proof}
The conformal invariance of Green functions implies that $G_{\Lambda}(x,y) = G_{\mathbb{H}^{+}}(\varphi(x), \varphi(y))$ for all $(x,y) \in \Lambda \times \Lambda$.
Since $G_{\mathbb{H}^{+}}$ admits the explicit formula
\[
    G_{\mathbb{H}^{+}}(w,z) = \log \left\vert \frac{w - \overline{z}}{w - z} \right\vert, \quad \forall (w,z) \in \mathbb{H}^{+} \times \mathbb{H}^{+},
\]
we may compute
\[
    \begin{array}[c]{rl}
        I_{n}^{(0)} := & \left( \frac{\mathcal{C}^{2}}{\sqrt{2}} \right)^{n} \sum\limits_{\{ \gamma_{k} \}} \int_{\Lambda^{n}} \exp \left\{ - \frac{1}{2} \sum\limits_{k < l} \gamma_{k} \gamma_{l} G_{\Lambda}(x_{k}, x_{l}) \right\} \prod\limits _{k = 1}^{n} \rho(x_{k}) \mu_{\varphi}(\mathrm{d}x_{k}) \\
        = & \left( \frac{\mathcal{C}^{2}}{\sqrt{2}} \right)^{n} \sum\limits_{\{ \gamma_{k} \}} \int_{\Lambda^{n}} \prod\limits_{k < l} \exp \left\{ - \frac{1}{2} \gamma_{k} \gamma_{l} \log \left\vert \frac{\varphi(x_{k}) - \overline{\varphi(x_{l})}}{\varphi(x_{k}) - \varphi(x_{l})} \right\vert \right\} \prod\limits _{k = 1}^{n} \rho(x_{k}) \mu_{\varphi}(\mathrm{d}x_{k}) \\
        = & \int_{\Lambda^{n}} \left[ \left( \frac{\mathcal{C}^{2}}{\sqrt{2}} \right)^{n} \prod\limits_{k=1}^{n} \left( \frac{\left\vert \varphi^{\prime}(x_{k}) \right\vert}{2 \mathrm{Im}\{ \varphi(x_{k}) \}} \right) ^{\frac{1}{4}} \sum\limits_{\{ \gamma_{k} \}} \prod\limits_{k < l} \left\vert \frac{\varphi(x_{k}) - \varphi(x_{l})}{\varphi(x_{k}) - \overline{\varphi(x_{l})}} \right\vert ^{\frac{\gamma_{k} \gamma_{l}}{2}} \right] \prod\limits _{k = 1}^{n} \rho(x_{k}) \mathrm{d}x_{k}.
    \end{array}
\]
(Here, we simply denote the summation over $\{ \gamma_{k} \}_{k=1}^{n} \in \{ \pm 1 \}^{n}$ by $\{ \gamma_{k} \}$.)
Consequently, the desired result follows immediately from the relationship
$
\mathcal{Q}_{\rho, \mu_{\varphi}}\left( 2^{-\frac{1}{2}} \mathcal{C}^{2}, 2^{-\frac{1}{2}} \right) = \sum_{n=0}^{\infty} \left( I_{n}^{(0)}/2^{n}n! \right)
$.
\end{proof}

The (zero-magnetic) XOR-Ising model is the law of a random spin configuration $\mathcal{S}_{\delta}$ on $\mathbf{F}_{\delta}$ given by a point-wise product of two independent Ising spin configurations,
$\sigma_{\delta}, \tilde{\sigma}_{\delta}$, being identically distributed to $\mathbb{P}_{\Lambda_{\delta}}^{+}$, i.e., $\mathcal{S}_{\delta}(a) = \sigma_{\delta}(a) \tilde{\sigma}_{\delta}(a)$ for $a \in \mathbf{F}_{\delta}$.
The corresponding spin field is defined by setting $\mathcal{S}_{\delta}(x) := \sigma_{\delta}(x) \tilde{\sigma}_{\delta}(x)$ for $x \in \Lambda$.
Based on Proposition \ref{sec-4-proposition-1}, $\mathcal{Q}_{\rho, \mu_{\varphi}}\left( 2^{-\frac{1}{2}} \mathcal{C}^{2}, 2^{-\frac{1}{2}} \right)$ can be further expressed by the scaling limit of an exponential moment for $\mathcal{S}_{\delta}$ tested against $\rho$ as $\delta \rightarrow 0^{+}$.

\begin{proposition}
    \label{sec-4-proposition-2}
    Let $\{ \mathcal{S}_{\delta} \}_{\delta>0}$ be the family of XOR-Ising spin fields. Then, for any nonnegative $\rho \in \mathcal{D}_{0}(\Lambda)$,
    \[
        \mathcal{Q}_{\rho, \mu_{\varphi}}\left( \frac{\mathcal{C}^{2}}{\sqrt{2}}, \frac{1}{\sqrt{2}} \right) = \lim\limits _{\delta \rightarrow 0^{+}} \mathbb{E}_{\Lambda_{\delta}}^{+} \left[ \exp \left\{ \frac{\delta^{-\frac{1}{4}}}{2} \int_{\Lambda} \rho(x) \mathcal{S}_{\delta}(x) \mathrm{d}x \right\} \right].
    \]
\end{proposition}

\begin{proof}
According to Theorem 1.2 in \cite{Chelkak-Hongler-Izyurov-2015}, for any conformal bijection $\varphi: \Lambda \mapsto \mathbb{H}^{+}$ and arbitrary $x_{1},\ldots,x_{n} \in \Lambda$ being distinct, we have
\begin{equation}
    \label{sec-4-Ising-corr-scale-lim}
    \lim\limits _{\delta \rightarrow 0^{+}} \delta^{-\frac{n}{8}} \mathbb{E}_{\Lambda_{\delta}}^{+} \left[ \prod_{k=1}^{n} \sigma_{\delta}(x_{k}) \right] = \mathcal{C}^{n} \cdot \left\langle \sigma_{\varphi(x_{1})}, \cdots, \sigma_{\varphi(x_{n})} \right\rangle _{\mathbb{H}^{+}}^{+} \cdot \prod_{k=1}^{n} \left\vert \varphi^{\prime}(x_{k}) \right\vert ^{\frac{1}{8}}.
\end{equation}
Using a variant of the Onsager inequality for the Ising model (\cite{Webb-AAP-2019}, (4.3) on p. 2150) yields that, for any $n \in \mathbb{N}_{+}$ and arbitrary $x_{1},\ldots,x_{n} \in \mathrm{supp}(\rho)$ being distinct,
\begin{equation}
    \label{sec-4-Ising-type-Onsager-ineq}
    0 < \delta^{-\frac{n}{8}} \mathbb{E}_{\Lambda_{\delta}}^{+} \left[ \prod_{k=1}^{n}\sigma_{\delta}(x_{k}) \right] \leqslant C_{1}^{n} \prod_{k=1}^{n} \left( \min\limits _{k \neq l} \left\vert x_{k} - x_{l} \right\vert \right) ^{-\frac{1}{8}}
\end{equation}
for some constant $C_{1}>0$ being independent of $n$, $\delta$, $x_{1},\ldots,x_{n}$, but it may depend on $\mathrm{supp}(\rho)$.
Thus we deduce from Lemma 3.10 in \cite{Webb-AAP-2019} that there is a constant $C_{2}>0$ depending on $b_{\rho}$ and $\mathrm{supp}(\rho)$ such that
\begin{equation}
    \label{sec-4-Ising-corr-summand-In-delta}
    \begin{array}[c]{rl}
        I_{n}^{(\delta)} := & \displaystyle \int_{\Lambda^{n}} \left\{ \prod_{k=1}^{n} \rho(x_{k}) \cdot \left( \delta^{-\frac{n}{8}} \mathbb{E}_{\Lambda_{\delta}}^{+} \left[ \prod_{k=1}^{n} \sigma_{\delta}(x_{k}) \right] \right)^{2} \right\} \prod_{k=1}^{n} \mathrm{d}x_{k} \\
        \leqslant & \displaystyle (C_{1}^{2} b_{\rho})^{n} \int_{\left(\mathrm{supp}(\rho)\right)^{n}} \prod_{k=1}^{n} \left( \min _{k \neq l} \left\vert x_{k} - x_{l} \right\vert \right) ^{-\frac{1}{4}} \prod_{k=1}^{n} \mathrm{d}x_{k} \\
        \leqslant & C_{2}^{n} n^{\frac{n}{8}}
    \end{array}
\end{equation}
holds for all $\delta > 0$, which together with (\ref{sec-4-Ising-corr-scale-lim}), (\ref{sec-4-Ising-type-Onsager-ineq}), and (\ref{sec-4-Ising-corr-summand-In-delta}) implies that $\lim_{\delta \rightarrow 0^{+}} I_{n}^{(\delta)} = I_{n}^{(0)}$ by the dominated convergence theorem.
Furthermore, putting $\beta = 2^{-\frac{1}{2}}$ and $\psi(x) = \left( \left\vert \varphi^{\prime}(x) \right\vert / 2 \mathrm{Im}\{ \varphi(x) \} \right) ^{\frac{1}{4}}$ in (\ref{sec-4-partit-summand-dominate}) yields $I_{n}^{(0)} \leqslant C_{3}^{n} n^{\frac{n}{8}}$ for some constant $C_{3}>0$ depending on $b_{\rho}$, $b_{\psi}$, and $\mathrm{supp}(\rho)$, 
with which we combine (\ref{sec-4-Ising-corr-summand-In-delta}) to ultimately obtain $\sup_{\delta \geqslant 0}  I_{n}^{(\delta)} \leqslant C^{n} n^{\frac{n}{8}}$ for some constant $C>0$ depending on $b_{\rho}$, $b_{\psi}$, and $\mathrm{supp}(\rho)$.
Then it follows from a Stirling estimate that the series $\sum_{n=0}^{\infty} \left( I_{n}^{(\delta)}/2^{n} n! \right)$ converges uniformly in $\delta \geqslant 0$.
Based on the above analysis, the limitation commutes with the summation and then it leads to
\[
    \mathcal{Q}_{\rho, \mu_{\varphi}}\left( \frac{\mathcal{C}^{2}}{\sqrt{2}}, \frac{1}{\sqrt{2}} \right) = \sum\limits_{n=0}^{\infty} \frac{I_{n}^{(0)}}{2^{n} n!} = \sum\limits_{n=0}^{\infty} \lim\limits _{\delta \rightarrow 0^{+}} \frac{I_{n}^{(\delta)}}{2^{n} n!} = \lim\limits _{\delta \rightarrow 0^{+}} \sum\limits_{n=0}^{\infty} \frac{I_{n}^{(\delta)}}{2^{n} n!}.
\]
For $n \in \mathbb{N}_{+}$ and $\delta>0$, by the definition of $\mathcal{S}_{\delta}$ and the Fubini theorem, we have
\begin{align*}
    I_{n}^{(\delta)} = & \int_{\Lambda^{n}} \left\{ \prod_{k=1}^{n} \rho(x_{k}) \cdot \delta^{-\frac{n}{4}} \mathbb{E}_{\Lambda_{\delta}}^{+} \left[ \prod_{k=1}^{n} \mathcal{S}_{\delta}(x_{k}) \right] \right\} \prod_{k=1}^{n} \mathrm{d}x_{k} \\
    = & \ \mathbb{E}_{\Lambda_{\delta}}^{+} \left[ \left( \delta^{-\frac{1}{4}} \int_{\Lambda} \rho(x) \mathcal{S}_{\delta}(x) \mathrm{d}x \right)^{n} \right].
\end{align*}
Finally, applying the Vitali convergence theorem yields 
\begin{align*}
    \mathcal{Q}_{\rho, \mu_{\varphi}}\left( \frac{\mathcal{C}^{2}}{\sqrt{2}}, \frac{1}{\sqrt{2}} \right) 
    = & \lim\limits _{\delta \rightarrow 0^{+}} \sum\limits_{n=0}^{\infty} \frac{1}{n!} \mathbb{E}_{\Lambda_{\delta}}^{+} \left[ \left( \frac{\delta^{-\frac{1}{4}}}{2} \int_{\Lambda} \rho(x) \mathcal{S}_{\delta}(x) \mathrm{d}x \right)^{n} \right] \\
    = & \lim\limits _{\delta \rightarrow 0^{+}} \mathbb{E}_{\Lambda_{\delta}}^{+} \left[ e^{ \frac{\delta^{-\frac{1}{4}}}{2} \int_{\Lambda} \rho(x) \mathcal{S}_{\delta}(x) \mathrm{d}x } \right]
\end{align*}
since we have
\[
\sup\limits _{\delta>0} \mathbb{E}_{\Lambda_{\delta}}^{+} \left[ \exp \left\{ p \left\vert \delta^{-\frac{1}{4}} \int_{\Lambda} \rho(x) \mathcal{S}_{\delta}(x) \mathrm{d}x \right\vert \right\} \right] < +\infty, \quad \forall p \in (0, +\infty)
\]
according to Lemma 4.3 in \cite{Webb-AAP-2019}.
\end{proof}

\subsection{Weak convergence of the normalized charge or spin distributions}

In this part, we focus on the limit of the sine-Gordon representation for distributions of particle/spin configurations under normalization. For simplicity of writing, we put $\psi \equiv 1$ without loss of generality and hence $\mu(\mathrm{d}x) = \mathrm{d}x$.
The state space of particle/spin configurations is $\varPhi_{\Lambda} := \{ (n,\mathbf{x},\Upsilon) : n \in \mathbb{N}, \mathbf{x} \in \Lambda^{n}, \Upsilon \in \{\pm 1\}^{n} \}$ equipped with its canonical $\sigma$-algebra $\mathcal{A}$.
A measurable function $F$ on $\varPhi_{\Lambda}$ is thus a sequence $\{F(n,\cdot,\cdot)\}_{n \in \mathbb{N}}$ of measurable functions on $\Lambda^{n} \times \{-1,1\}^{n} $.
The normalized partition function (\ref{sec-4-renormalized-part-func}) induces a probability measure $\mathbb{P}_{\rho, \varepsilon}$ (with expectation $\mathbb{E}_{\rho, \varepsilon}$) by setting (below $\mathbf{x} = (x_{k})_{k=1}^{n}$, $\Upsilon = ( \gamma_{k} )_{k=1}^{n}$ are vectors with $n$ coordinates)
\[
    \Xi_{\rho, \varepsilon}^{-1} \sum\limits _{n=0}^{\infty} \frac{\alpha^{n}}{2^{n}n!} \sum\limits_{\Upsilon \in \{ \pm 1 \}^{n}} \int_{\Lambda_{\varepsilon}^{n}} F(n,\mathbf{x},\Upsilon) e^{ - \beta^{2} \sum\limits_{1 \leqslant k < l \leqslant n} \gamma_{k} \gamma_{l} G_{\Lambda}^{\varepsilon}(x_{k}, x_{l}) } \prod\limits _{k = 1}^{n} \rho(x_{k}) \mathrm{d}x_{k}
\]
for arbitrary bounded measurable function $F$ on $\varPhi_{\Lambda}$ as its expectation $\mathbb{E}_{\rho, \varepsilon} [F]$.

where we adopt $\Xi_{\rho, \varepsilon}$ to be the normalizing constant as we have shown that $\mathcal{Q}_{\rho, \mu}^{(\varepsilon)}(\alpha, \beta) = \Xi_{\rho, \varepsilon}$ in Proposition \ref{sec-4-proposition-3}.
The physically relevant quantity to be concerned is the charge/spin distribution, being obtained as the push-forward (or image measure) of $\mathbb{P}_{\rho, \varepsilon}$ by a measurable mapping from $(\varPhi_{\Lambda}, \mathcal{A})$ to $(\mathcal{S}^{\prime}(\mathbb{R}^{2}), \mathfrak{B}(\mathcal{S}^{\prime}(\mathbb{R}^{2})))$:
\[
\begin{array}[c]{cccc}
    \mathcal{T}: & \varPhi_{\Lambda} & \longmapsto & \mathcal{S}^{\prime}(\mathbb{R}^{2}) \\
    & (n,\mathbf{x},\Upsilon) & \longmapsto & \sum\limits _{k=1}^{n} \gamma_{k} \delta_{x_{k}},
\end{array}
\]
which sums signed Dirac masses $\pm \delta_{x_{k}}$ corresponding to particles' (resp. spins') locations and charges (resp. states of spins).
(In fact, $\mathcal{T}(\varPhi_{\Lambda}) \subset \mathcal{S}^{\prime}(\mathbb{R}^{2})$ since $\mathrm{supp}(\sum_{k=1}^{n} \gamma_{k} \delta_{x_{k}}) = \{ x_{k} \} _{k=1}^{n} \subset \Lambda$ is compact in $\mathbb{R}^{2}$ for each $(n,\mathbf{x},\Upsilon) \in \varPhi_{\Lambda}$.)
To establish the sine-Gordon representation for the characteristic function of charge/spin distribution, we introduce the Fourier transform of image measure $\mathcal{T}_{\#} \mathbb{P}_{\rho, \varepsilon}$ by setting
\[
    \Psi_{\mathcal{T}_{\#} \mathbb{P}_{\rho, \varepsilon}}(\theta) := \int_{\mathcal{S}^{\prime}(\mathbb{R}^{2})} e^{i \nu(\theta)} \mathcal{T}_{\#} \mathbb{P}_{\rho, \varepsilon} \left( \mathrm{d}\nu \right), \quad \forall \theta \in \mathcal{S}(\mathbb{R}^{2}).
\]
We have the following main theorem of this part.

\begin{theorem}
    \label{sec-4-main-thm}
    Let $\rho \in \mathcal{D}_{0}(\Lambda)$ be nonnegative and $\Gamma(\rho)$ be defined in Theorem \ref{sec-3-theorem-1}. If $\beta^{2} \in [0,2)$, then there is a functional $\Psi_{\rho}: \mathcal{S}(\mathbb{R}^{2}) \mapsto \mathbb{R}$ being defined by
    \[
        \Psi_{\rho}(\theta) := \mathbb{E} \left[ \Gamma(\rho) \exp\left\{ -\alpha \left( \left\langle \sin\left( \beta W_{1} \right), \rho \right\rangle _{\mu_{1}} + \left\langle \cos\left( \beta W_{1} \right), \rho \right\rangle _{\mu_{2}} \right) \right\} \right]
    \]
    with two locally finite Borel measures $\mu_{1}(\mathrm{d}x) := \sin(\theta(x)) \mathrm{d}x$, $\mu_{2}(\mathrm{d}x) := \left( 1 - \cos (\theta(x)) \right) \mathrm{d}x$, such that $\Psi_{\rho}$ is continuous with respect to the topology on $\mathcal{S}(\mathbb{R}^{2})$ and
    \begin{equation}
        \label{sec-4-thm-char-func-limit}
        \lim\limits _{\varepsilon \rightarrow 0^{+}} \Psi_{\mathcal{T}_{\#} \mathbb{P}_{\rho, \varepsilon}}(\theta) = \Psi_{\rho}(\theta), \quad \forall \theta \in \mathcal{S}(\mathbb{R}^{2}).
    \end{equation}
    In particular, the family of probability measures $\{ \mathcal{T}_{\#} \mathbb{P}_{\rho, \varepsilon} \}_{\varepsilon \in (0,1]}$ converge weakly to the law of some $\mathcal{S}^{\prime}(\mathbb{R}^{2})$-valued random variable on $(\varPhi_{\Lambda}, \mathcal{A})$ with $\Psi_{\rho}$ as its characteristic function, 
    which is independent of the choice of smoothing kernel $\eta$.
\end{theorem}

\begin{proof}
For any $\theta \in \mathcal{S}(\mathbb{R}^{2})$, by repeating the computation leading to (\ref{sec-4-charge-partition-GFF-express}), $\Psi_{\mathcal{T}_{\#} \mathbb{P}_{\rho, \varepsilon}}(\theta)$ can be further expressed by
\begin{align*}  
        & \Xi_{\rho, \varepsilon}^{-1} \sum\limits _{n=0}^{\infty} \frac{\alpha^{n}}{2^{n}n!} \sum\limits_{\Upsilon \in \{ \pm 1 \}^{n}} \int_{\Lambda_{\varepsilon}^{n}} e^{ i \sum\limits _{k=1}^{n} \gamma_{k} \theta(x_{k}) - \beta^{2} \sum\limits _{1 \leqslant k < l \leqslant n} \gamma_{k} \gamma_{l} G_{\Lambda}^{\varepsilon}(x_{k}, x_{l}) } \prod\limits _{k = 1}^{n} \rho(x_{k}) \mathrm{d}x_{k} \\
        = & \ \Xi_{\rho, \varepsilon}^{-1} \mathbb{E} \left[ \exp\left\{ \alpha \int_{\Lambda_{\varepsilon}} e^{\frac{\beta^{2}}{2} G_{\Lambda}^{\varepsilon}(x) } \cos(\beta W_{1}^{\varepsilon}(x) + \theta(x)) \rho(x) \mathrm{d}x \right\} \right] \\
        = & \ \int_{H_{0}^{-1}(\Lambda)} e^{ \alpha \int_{\Lambda_{\varepsilon}} e^{\frac{\beta^{2}}{2} G_{\Lambda}^{\varepsilon}(x) } \left( \cos(\beta \phi^{\varepsilon}(x) + \theta(x)) - \cos(\beta \phi^{\varepsilon}(x)) \right) \rho(x) \mathrm{d}x } \mu_{\mathrm{SG}}^{\rho, \varepsilon}(\mathrm{d}\phi),
\end{align*}
where we use the notation $\phi^{\varepsilon} := \phi \ast \eta_{\varepsilon}$ and the last line results from the definition of $\mu_{\mathrm{SG}}^{\rho, \varepsilon}$. By transforming the difference of cosines as follows
\[
    \cos(\beta \phi^{\varepsilon}(x) + \theta(x)) - \cos(\beta \phi^{\varepsilon}(x)) = -2 \sin\left( \beta \phi^{\varepsilon}(x) + \frac{\theta(x)}{2} \right) \sin\left( \frac{\theta(x)}{2} \right)
\]
and defining a function in $\mathcal{C}_{b}(H_{0}^{-1}(\Lambda))$ as
\begin{equation}
    \label{sec-4-def-F-Cb}
    F_{\rho, \theta}^{(\varepsilon)}(\phi) := \int_{\Lambda_{\varepsilon}} e^{ \frac{\beta^{2}}{2} G_{\Lambda}^{\varepsilon}(x) } \sin\left( \beta \phi^{\varepsilon}(x) + \frac{\theta(x)}{2} \right) \sin\left( \frac{\theta(x)}{2} \right) \rho(x) \mathrm{d}x,
\end{equation}
the above Fourier transform can be further written as
\begin{equation}
    \label{sec-4-approxi-char-func}
    \Psi_{\mathcal{T}_{\#} \mathbb{P}_{\rho, \varepsilon}}(\theta)
        = \int_{H_{0}^{-1}(\Lambda)} e^{ -2\alpha F_{\rho, \theta}^{(\varepsilon)}(\phi) } \mu_{\mathrm{SG}}^{\rho, \varepsilon}(\mathrm{d}\phi)
        = \mathbb{E} \left[ \Gamma_{1}^{\varepsilon}(\rho) e^{ -2\alpha F_{\rho, \theta}^{(\varepsilon)}(W_{1}) } \right]
\end{equation}
thanks to Lemma \ref{sec-3-lemma-3}. 
The uniform boundedness of $F_{\rho, \theta}^{(\varepsilon)}$ is easy to verify since we deduce from (\ref{sec-3-terminal-L1}) and (\ref{sec-4-def-F-Cb}) that, for all $\phi \in H_{0}^{-1}(\Lambda)$ and $\varepsilon \in (0,1]$,
\[
    \left\vert F_{\rho, \theta}^{(\varepsilon)}(\phi) \right\vert \leqslant b_{\rho} \int_{\Lambda_{\varepsilon}} e^{ \frac{\beta^{2}}{2} G_{\Lambda}^{\varepsilon}(x) } \mathrm{d}x \leqslant b_{\rho} C_{\beta}.
\]
To prove the continuity, by choosing $\phi_{1}, \phi_{2} \in H_{0}^{-1}(\Lambda)$ arbitrarily, we obtain the following error estimate:
\begin{align}
    \label{sec-4-error-estimate}
        & \left\vert F_{\rho, \theta}^{(\varepsilon)}(\phi_{1}) - F_{\rho, \theta}^{(\varepsilon)}(\phi_{2}) \right\vert \\
        \nonumber \leqslant & \displaystyle b_{\rho} \int_{\Lambda_{\varepsilon}} e^{ \frac{\beta^{2}}{2} G_{\Lambda}^{\varepsilon}(x) } \left\vert \sin\left( \beta \phi_{1}^{\varepsilon}(x) + \frac{\theta(x)}{2} \right) - \sin\left( \beta \phi_{2}^{\varepsilon}(x) + \frac{\theta(x)}{2} \right) \right\vert \mathrm{d}x \\
        \nonumber \leqslant & \displaystyle \beta b_{\rho} \int_{\Lambda_{\varepsilon}} e^{ \frac{\beta^{2}}{2} G_{\Lambda}^{\varepsilon}(x) } \left\vert \phi_{1}^{\varepsilon}(x) - \phi_{2}^{\varepsilon}(x) \right\vert \mathrm{d}x \\
        \nonumber = & \displaystyle \beta b_{\rho} \int_{\Lambda_{\varepsilon}} e^{ \frac{\beta^{2}}{2} G_{\Lambda}^{\varepsilon}(x) } \left\vert \langle \phi_{1} - \phi_{2}, \eta_{\varepsilon}(x - \cdot) \rangle \right\vert \mathrm{d}x,
\end{align}
where the second line is due to the fundamental inequality $\left\vert \sin(a) - \sin(b) \right\vert \leqslant |a - b|$.
To implement estimate (\ref{sec-4-error-estimate}), noting that $\eta_{\varepsilon}(x - \cdot) \in \mathcal{D}_{0}(\Lambda)$ whenever $x \in \Lambda_{\varepsilon}$ and following the argument in Remark \ref{sec-2-remark-2}, we deduce from the definition of $H_{0}^{1}(\Lambda)$-norm and the Cauchy-Schwartz inequality that, for all $x \in \Lambda_{\varepsilon}$,
\[
\begin{array}[c]{rl}
    \left\Vert \eta_{\varepsilon}(x - \cdot) \right\Vert _{H_{0}^{1}(\Lambda)} \leqslant & \displaystyle \left\Vert \varepsilon^{-4} \Delta \eta \left( \frac{x - \cdot}{\varepsilon} \right) \right\Vert _{L^{2}(\Lambda)} \left\Vert \varepsilon^{-2} \eta \left( \frac{x - \cdot}{\varepsilon} \right) \right\Vert _{L^{2}(\Lambda)} \\
    \leqslant & \displaystyle \varepsilon^{-4} \left\Vert \varepsilon^{-1} \Delta \eta \left( \frac{x - \cdot}{\varepsilon} \right) \right\Vert _{L^{2}(\mathbb{R}^{2})} \left\Vert \varepsilon^{-1} \eta \left( \frac{x - \cdot}{\varepsilon} \right) \right\Vert _{L^{2}(\mathbb{R}^{2})} \\
    = & \displaystyle \varepsilon^{-4} \left\Vert \Delta \eta \right\Vert _{L^{2}(\mathbb{R}^{2})} \left\Vert \eta \right\Vert _{L^{2}(\mathbb{R}^{2})},
\end{array}
\]
where the last line attributes to the change of variables and the fact $\eta \in \mathcal{D}_{0}(\mathbb{R}^{2})$.
Thus, thanks to the duality between $H_{0}^{-1}(\Lambda)$ and $H_{0}^{1}(\Lambda)$, there is a constant $C_{\eta}>0$ depending only on $\eta$ such that
\begin{align*}
    \left\vert \langle \phi_{1} - \phi_{2}, \eta_{\varepsilon}(x - \cdot) \rangle \right\vert 
    \leqslant & \left\Vert \phi_{1} - \phi_{2} \right\Vert _{H_{0}^{-1}(\Lambda)} \left\Vert \eta_{\varepsilon}(x - \cdot) \right\Vert _{H_{0}^{1}(\Lambda)} \\
    \leqslant & C_{\eta} \varepsilon^{-4} \left\Vert \phi_{1} - \phi_{2} \right\Vert _{H_{0}^{-1}(\Lambda)}.
\end{align*}
Plugging this into (\ref{sec-4-error-estimate}) and utilizing estimate (\ref{sec-3-terminal-L1}), we obtain
\begin{align*}
    \left\vert F_{\rho, \theta}^{(\varepsilon)}(\phi_{1}) - F_{\rho, \theta}^{(\varepsilon)}(\phi_{2}) \right\vert \leqslant & \beta b_{\rho} C_{\eta} \varepsilon^{-4} \left\Vert \phi_{1} - \phi_{2} \right\Vert _{H_{0}^{-1}(\Lambda)} \int_{\Lambda_{\varepsilon}} e^{ \frac{\beta^{2}}{2} G_{\Lambda}^{\varepsilon}(x) } \mathrm{d}x \\
    \leqslant & C_{\eta, \varepsilon} \left\Vert \phi_{1} - \phi_{2} \right\Vert _{H_{0}^{-1}(\Lambda)},
\end{align*}
for some constant $C_{\eta, \varepsilon}>0$ depending on $\beta$, $b_{\rho}$, $\eta$, and $\varepsilon$.

To prove (\ref{sec-4-thm-char-func-limit}), we first compute $F_{\rho, \theta}^{(\varepsilon)}(W_{1})$ by
\begin{align*}
     & \int_{\Lambda_{\varepsilon}} e^{ \frac{\beta^{2}}{2} G_{\Lambda}^{\varepsilon}(x) } \left[ \sin\left( \beta W_{1}^{\varepsilon}(x) \right) \cos\left( \frac{\theta(x)}{2} \right) \sin\left( \frac{\theta(x)}{2} \right) + \cos\left( \beta W_{1}^{\varepsilon}(x) \right) \sin^{2} \left( \frac{\theta(x)}{2} \right) \right] \rho(x) \mathrm{d}x \\
        = & \ \int_{\Lambda_{\varepsilon}} \left( \left[\left[ \sin\left( \beta W_{1}^{\varepsilon}(x) \right) \right]\right] \frac{\sin\left( \theta(x) \right)}{2} + \left[\left[ \cos\left( \beta W_{1}^{\varepsilon}(x) \right) \right]\right] \frac{1 - \cos\left( \theta(x) \right)}{2} \right) \rho(x) \mathrm{d}x \\
        = & \ \frac{1}{2} \left[ \left( \left[\left[ \sin\left( \beta W_{1}^{\varepsilon} \right) \right]\right], \rho \right) _{L^{2}(\Lambda_{\varepsilon};\mu_{1})} + \left( \left[\left[ \cos\left( \beta W_{1}^{\varepsilon} \right) \right]\right], \rho \right) _{L^{2}(\Lambda_{\varepsilon};\mu_{2})} \right],
\end{align*}
where $\mu_{1}(\mathrm{d}x) := \sin(\theta(x)) \mathrm{d}x$, $\mu_{2}(\mathrm{d}x) := \left( 1 - \cos (\theta(x)) \right) \mathrm{d}x$ are two locally finite signed measures.
Then it follows from (\ref{sec-3-cosine-W-converge-prob}) and (\ref{sec-3-sine-W-converge-prob}) that the following limit
\[
\lim\limits_{\varepsilon \rightarrow 0^{+}} F_{\rho, \theta}^{(\varepsilon)}(W_{1}) = \frac{1}{2} \left( \left\langle \sin\left( \beta W_{1} \right), \rho \right\rangle _{\mu_{1}} + \left\langle \cos\left( \beta W_{1} \right), \rho \right\rangle _{\mu_{2}} \right)
\]
exists in the sense of convergence in $\mathbb{P}$ and it is independent of the choice of smoothing kernel $\eta$.
Since $\{ F_{\rho, \theta}^{(\varepsilon)}(\phi) \}_{\varepsilon \in (0,1]} \subset \mathcal{C}_{b}(H_{0}^{-1}(\Lambda)) \subset \mathcal{B}_{b}(H_{0}^{-1}(\Lambda))$, it follows from Doob's measurability theorem (see, e.g. \cite{Book-He-Wang-Yan-1992}, Theorem 1.5) that $\{ \exp \{ -2\alpha F_{\rho, \theta}^{(\varepsilon)}(W_{1}) \} \}_{\varepsilon \in (0,1]}$ is a family of $\sigma(W_{1})$-measurable bounded random variables.
According to (\ref{sec-3-thm-1-general-convergence}) in Theorem \ref{sec-3-theorem-1}, we deduce (\ref{sec-4-thm-char-func-limit}) from letting $\varepsilon \rightarrow 0^{+}$ in (\ref{sec-4-approxi-char-func}). 

It remains to prove the continuity of $\Psi_{\rho}$. To this end, it suffices to show that, for any fixed $\phi \in H_{0}^{-1}(\Lambda)$, the function $\theta \mapsto F_{\rho, \theta}^{(\varepsilon)}(\phi)$ is continuous with respect to the topology on $\mathcal{S}(\mathbb{R}^{2})$, uniformly in $\varepsilon$.
Actually, choosing $\theta_{1}, \theta_{2} \in \mathcal{S}(\mathbb{R}^{2})$ arbitrarily and observing (\ref{sec-4-def-F-Cb}), it follows from the fundamental inequality $\left\vert \sin(a) - \sin(b) \right\vert \leqslant |a - b|$ and (\ref{sec-3-terminal-L1}) that
\begin{align*}
    \left\vert F_{\rho, \theta_{1}}^{(\varepsilon)}(\phi) - F_{\rho, \theta_{2}}^{(\varepsilon)}(\phi) \right\vert \leqslant & \ 2 b_{\rho} \int_{\Lambda_{\varepsilon}} e^{ \frac{\beta^{2}}{2} G_{\Lambda}^{\varepsilon}(x) } \left\vert \theta_{1}(x) - \theta_{2}(x) \right\vert \mathrm{d}x \\
    \leqslant & \ 2 b_{\rho} C_{\beta} \left\Vert \theta_{1} - \theta_{2} \right\Vert _{0,0},
\end{align*}
where $C_{\beta}$ is the constant appearing in the estimate (\ref{sec-3-terminal-L1}).
Since $\mathbb{E}[\Gamma_{1}^{\varepsilon}(\rho)] = 1$ for all $\varepsilon \in (0,1]$, we deduce from the above inequality that
\begin{align*}
    & \mathbb{E}\left[ \Gamma_{1}^{\varepsilon}(\rho) \left\vert e^{ -2\alpha F_{\rho, \theta_{1}}^{(\varepsilon)}(W_{1}) } - e^{ -2\alpha F_{\rho, \theta_{2}}^{(\varepsilon)}(W_{1}) } \right\vert \right] \\
    \leqslant & 2|\alpha| e^{2|\alpha| b_{\rho} C_{\beta}} \mathbb{E}\left[ \Gamma_{1}^{\varepsilon}(\rho) \left\vert F_{\rho, \theta_{1}}^{(\varepsilon)}(W_{1}) - F_{\rho, \theta_{2}}^{(\varepsilon)}(W_{1}) \right\vert \right] \\
    \leqslant & 2|\alpha| b_{\rho} C_{\beta} e^{2|\alpha| b_{\rho} C_{\beta}} \left\Vert \theta_{1} - \theta_{2} \right\Vert _{0,0}.
\end{align*}
Thus, according to (\ref{sec-4-thm-char-func-limit}) and (\ref{sec-4-approxi-char-func}), we finally have
\[
    \left\vert \Psi_{\rho}(\theta_{1}) - \Psi_{\rho}(\theta_{2}) \right\vert = \lim\limits _{\varepsilon \rightarrow 0^{+}} \left\vert \Psi_{\mathcal{T}_{\#} \mathbb{P}_{\rho, \varepsilon}}(\theta_{1}) - \Psi_{\mathcal{T}_{\#} \mathbb{P}_{\rho, \varepsilon}}(\theta_{2}) \right\vert \leqslant C \left\Vert \theta_{1} - \theta_{2} \right\Vert _{0,0},
\]
where the constant $C$ depends only on $\alpha$, $\beta$, and $\rho$.
The last assertion follows immediately from the continuity of $\Psi_{\rho}$ and a consequence of Levy-Schwartz Theorem (\cite{Bierme-Durieu-Wang-2017}, Theorem 2.3), which accomplishes the proof.
\end{proof}

        





 \ack 
      The research of Shanjian Tang was partially supported by National Science Foundation of China (Grant No.11631004) and National Key R\&D Program of China (Grant No.2018YFA0703903).
      The research of Rundong Xu was partially supported by China Postdoctoral Science Foundation (Grant No.2024M760481) and Shanghai Postdoctoral Excellence Program (Grant No.2023201).


\frenchspacing
\bibliographystyle{cpam}
\bibliography{bibdata}

\end{document}